\documentclass[10pt,twoside]{article}
\usepackage{amsmath,amsthm}
\usepackage[latin1]{inputenc}
\usepackage{picinpar}
\usepackage{longtable}
\usepackage{amsrefs}
\usepackage[dvips]{graphicx}
\usepackage[mathscr]{eucal}
\usepackage{calc}
\usepackage{ifthen}
\usepackage{dcpic,pictexwd}
\usepackage[all]{xy}
\usepackage{enumerate}
\usepackage{amssymb,latexsym}
\usepackage{amsthm}
\usepackage[dvips]{graphics}
\usepackage{color}
\usepackage{tabularx}
\usepackage{makeidx}
\usepackage{tkz-graph}
\usepackage{tikz-cd}
\usetikzlibrary{calc}
\usepackage{tikz}
\usetikzlibrary{arrows}
\usepackage{stackrel}
\usepackage{multicol}
\usepackage{float}
\newtheorem{teorema}{Theorem}

\newtheorem{corolario}{Corollary}

\newtheorem{proposicion}{Proposition}

\theoremstyle{definition}

 \newtheorem{teor}{Theorem}

\newtheorem{prop}{Proposition}
\newtheorem{corol}{Corollary}

\newtheorem{Nota}{Remark}

\begin{document}
\thispagestyle{plain}
\par\bigskip
\begin{centering}

\textbf{Brauer Configuration Algebras and Matrix Problems to Categorify Integer Sequences}

\end{centering}\par\bigskip
\begin{centering}
\footnotesize{Agust\'{\i}n Moreno Ca\~{n}adas (corresponding author)}\\
  \footnotesize{Pedro Fernando Fern\'andez Espinosa}\\
      \footnotesize{Isa\'ias David Mar\'in Gaviria}\\
        \footnotesize{Gabriel Bravo R{i}os}\\
   
\end{centering}
\par\bigskip

\begin{centering}

\textbf{Abstract}\par\bigskip

\end{centering}

\small{Bijections between invariants associated to indecomposable projective modules over some suitable Brauer configuration algebras and invariants associated to solutions of the Kronecker problem and the four subspace problem are used to categorify integer sequences in the sense of Ringel and Fahr. Dimensions of the Brauer configuration algebras and their corresponding centers involved in the different processes are given as well.}


\section{Introduction}

According to Ringel and Fahr \cite{Fahr2} a categorification of a sequence of
numbers means to consider instead of these numbers suitable objects in a category (for instance, representation of quivers) so that the numbers in question occur as invariants of the objects, equality of numbers may be visualized by isomorphisms of objects functional relations by functorial ties.  The notion of this kind of categorification arose from  the use of suitable arrays of numbers to obtain integer partitions of dimensions of indecomposable preprojective modules over the 3-Kronecker algebra (see Figure (\ref{covering}) where it is shown the 3-Kronecker quiver and a piece of the oriented 3-regular tree or universal covering $(T,E,\Omega_{t})$ as described by Ringel and Fahr in \cite{Fahr1}). Firstly they noted that  the vector dimension of these kind of modules consists of even-index Fibonacci numbers (denoted $f_{i}$ and such that $f_{i}=f_{i-1}+f_{i-2}$, for $i\geq2$, $f_{0}=0$, $f_1=1$) then they used results from the universal covering theory developed by Gabriel and his students to identify such Fibonacci numbers with dimensions of representations of the corresponding universal covering. In particular, preinjective and preprojective representations of the 3-Kronecker quiver were used in \cite{Fahr1} by Ringel and Fahr in order to derive a partition formula for even-index Fibonacci numbers $f_{2n}$.

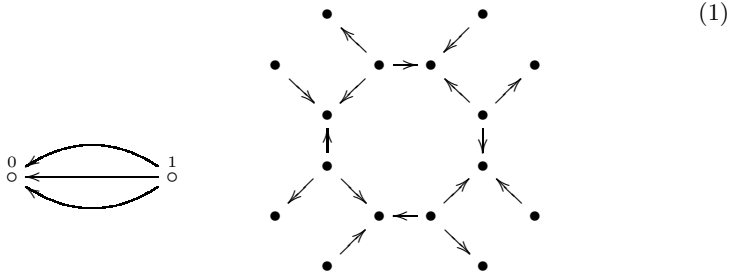
 \begin{equation}\label{covering}
\xymatrix@=50pt{
	&\\
\overset{0}{\circ}&\ar@/_1pc/[l]\ar@/^1pc/[l] \ar[l] \overset{1}{\circ}
}\hspace{10mm}
\xymatrix@=5pt @R=0.75pc @C=0.75pc{
&\bullet&&&\ar[dl]\bullet& \\
\bullet\ar[dr]&&\ar[ul]\ar[dl]\ar[r]\bullet&\bullet&&\bullet \\
&\bullet&&&\ar[ul]\ar[ur]\ar[d]\bullet& \\
&\ar[dl]\ar[dr]\ar[u]\bullet&&&\bullet& \\
\bullet&&\bullet&\ar[l]\ar[ur]\ar[dr]\bullet&&\ar[ul]\bullet \\
&\ar[ur]\bullet&&&\bullet&
}
\end{equation}

\par\bigskip
For the sake of clarity we give here a brief insight into the program of Ringel and Fahr.\par\bigskip

First of all note that the road to a categorification of the Fibonacci numbers has several stops some of them dealing with diagonal (lower) arrays of numbers of the form $D=(d_{i, j})$ with $0\leq j\leq i\leq n$, (columns numbered from right to the left, see array $(5)$) for some $n\geq0$ fixed and such that:

\begin{equation}\label{farray}
\begin{split}
d_{i,i}&=1, \quad\text{for all}\hspace{0.1cm}i\geq0,\\
d_{2, 0}&=2, \\
d_{2t+1,0}&=0,\quad\text{for all}\hspace{0.1cm}t\geq0,\\
d_{i+2t, i-1}&=0, \quad\text{if}\hspace{0.1cm}i\geq 1, \hspace{0.1cm}t\geq0,\\
2d_{i, j}+d_{i, j-2}-d_{i-1, j-1}&=d_{i+1, j-1},\quad i, j\geq2.
\end{split}
\end{equation}
Besides, if $i \geq 4$ then the following identity (hook rule) holds:

\begin{equation}
\begin{split}
\underset{t=0}{\overset{i-2}{\sum}}d_{i+t, i-t}+d_{2i-2, 0}&=d_{2i-1, 1}.
\end{split}
\end{equation}

Note that to each entry $d_{i, i-j}$ it is possible to assign a weight $w_{i, i- j}$ such that:

\[w_{i, i-j}=
\begin{cases}
3.2^{2\lfloor\frac{i-j}{2}\rfloor+a},  & \hspace{0.2cm}\text{if}\hspace{0.2cm}j\hspace{0.2cm}\text{is even},\hspace{0.2cm}i\neq j, \\
0,            &\hspace{0.2cm}\text{if}\hspace{0.2cm}j\hspace{0.2cm}\text{is odd},\hspace{0.2cm}i\neq j, \\
1,&\hspace{0.2cm}\text{if}\hspace{0.2cm}i=j=2h\hspace{0.2cm}\text{for some}\hspace{0.1cm}h\geq0.
\end{cases}\]

Where $\lfloor x\rfloor$ is the greatest integer number less than $x$, $a\in\{0,-1\}$, $a=-1$ if $i$ is even, it is 0 otherwise.\par\bigskip

The first stop consists of defining partitions of the even-index Fibonacci numbers in the following form:
\begin{equation}
\begin{split}
f_{2i+2}&=\underset{j=0}{\overset{i}{\sum}}(w_{i,i-j})(d_{i, i-j}),
\end{split}
\end{equation}

to do that, Ringel and Fahr interpreted weights $w_{i, i-j}$ as distances in a 3-regular tree $(T,E)$ (with $T$ a vertex set and $E$ a set of edges) from a fixed point $x_{0}\in T$ to any point $y\in T$. They define sets $T_{r}$ whose points have distance $r$ to $x_{0}$, in such a case $T_{0}=\{x_{0}\}$, $T_{1}$ are the neighbors of $x_{0}$ and so on. A given vertex $y$ is said to be even or odd according to this parity \cite{Fahr1}.  \par\bigskip

Any vertex $y\in T$ yields a suitable reflection $\sigma_{y}$ on
the set of functions $T\rightarrow \mathbb{Z}$ with finite support, denoted $\mathbb{Z}[T]$, and some reflection products denoted
 $\Phi_{0}$ and $\Phi_{1}$ according to the parity of $y$ are introduced in \cite{Fahr1}. Then some maps
$a_{t}:\mathbb{N}_{0}\rightarrow \mathbb{Z}\in\mathbb{Z}[T]$ are defined in such a way that if $a_{0}$ is
the characteristic function of $T_{0}$ then
\begin{center}
$a_{0}(x)=0$ unless $x=x_{0}$ in which case $a_{0}(x_{0})=1$,\quad and\quad $a_{t}=(\Phi_{0}\Phi_{1})^{t}a_{0}$, for $t\geq0$,\par\bigskip
\end{center}

with $a_{t}[r]=a_{t}(x)$, for $r\in\mathbb{N}_{0}$ and $x\in T_{r}$, these maps $a_{t}$ give the values $d_{i, j}$ of the array (\ref{farray}). The following table (called the even-index Fibonacci partition triangle) is an example of such array with $n=12$. Rows are giving by the values of $t$, $P_{t}$ is a notation for a 3-Kronecker preprojective module with dimension vector $[f_{2t+2}\hspace{0.2cm}f_{2t}]$  (see \cite{Fahr3}).
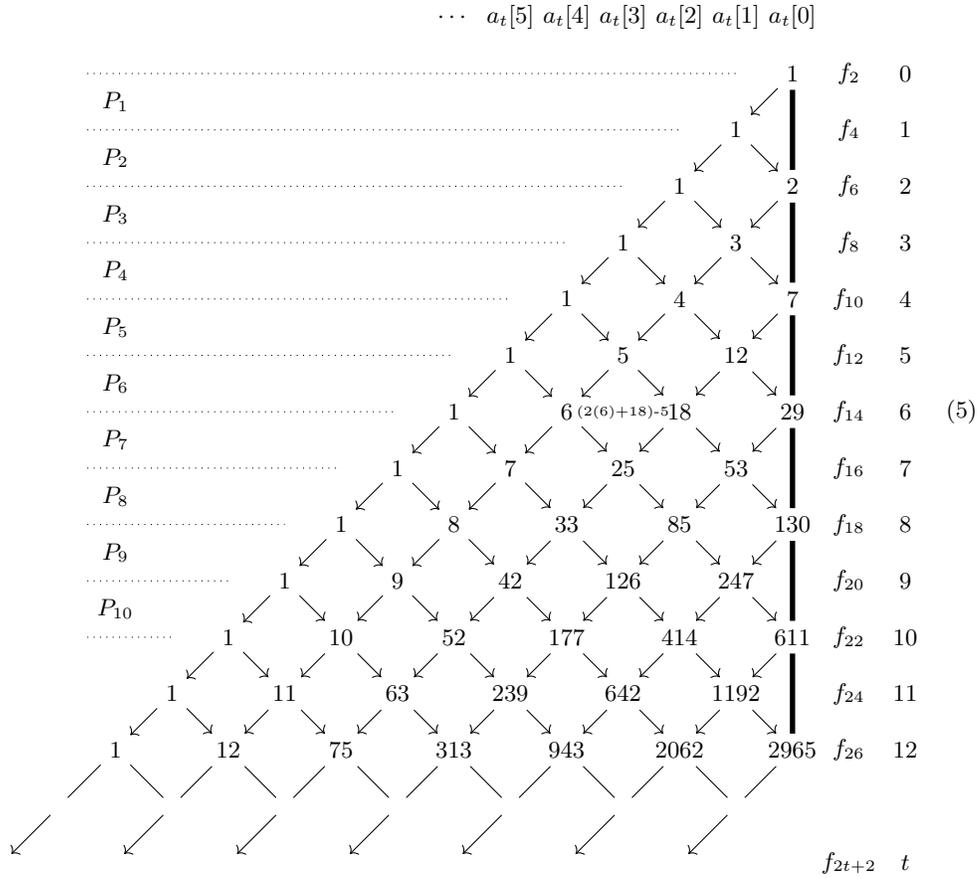
\par\bigskip
\begin{equation*}
\begin{tikzpicture}[scale=0.75]
\node at (11,9)(a){$(5)$};


\node at (10,15)(a){$0$};
\node at (10,14)(a){$1$};
\node at (10,13)(a){$2$};
\node at (10,12)(a){$3$};
\node at (10,11)(a){$4$};
\node at (10,10)(a){$5$};
\node at (10,9)(a){$6$};
\node at (10,8)(a){$7$};
\node at (10,7)(a){$8$};
\node at (10,6)(a){$9$};
\node at (10,5)(a){$10$};
\node at (10,4)(a){$11$};
\node at (10,3)(a){$12$};

\node at (10,1)(a){$t$};
\node at (9,15)(a){$f_2$};
\node at (9,14)(a){$f_4$};
\node at (9,13)(a){$f_6$};
\node at (9,12)(a){$f_8$};
\node at (9,11)(a){$f_{10}$};
\node at (9,10)(a){$f_{12}$};
\node at (9,9)(a){$f_{14}$};
\node at (9,8)(a){$f_{16}$};
\node at (9,7)(a){$f_{18}$};
\node at (9,6)(a){$f_{20}$};
\node at (9,5)(a){$f_{22}$};
\node at (9,4)(a){$f_{24}$};
\node at (9,3)(a){$f_{26}$};

\node at (9,1)(a){$f_{2t+2}$};

\node at (8,15)(1){$1$};
\node at (8,13)(2){$2$};
\node at (8,11)(3){$7$};
\node at (8,9)(4){$29$};
\node at (8,7)(5){$130$};
\node at (8,5)(6){$611$};
\node at (8,3)(7){$2965$};
\node at (7,14)(8){$1$};
\node at (7,12)(9){$3$};
\node at (7,10)(10){$12$};
\node at (7,8)(11){$53$};
\node at (7,6)(12){$247$};
\node at (7,4)(13){$1192$};
\node at (6,13)(14){$1$};
\node at (6,11)(15){$4$};
\node at (6,9)(16){$18$};
\node at (6,7)(17){$85$};
\node at (6,5)(18){$414$};
\node at (6,3)(19){$2062$};

\node at (5,12)(20){$1$};
\node at (5,10)(21){$5$};
\node at (5,8)(22){$25$};
\node at (5,6)(23){$126$};
\node at (5,4)(24){$642$};

\node at (4,11)(25){$1$};
\node at (4,9)(26){$6$};
\node at (5,9){\tiny{(2(6)+18)-5}};
\node at (4,7)(27){$33$};
\node at (4,5)(28){$177$};
\node at (4,3)(29){$943$};

\node at (3,10)(30){$1$};
\node at (3,8)(31){$7$};
\node at (3,6)(32){$42$};
\node at (3,4)(33){$239$};

\node at (2,9)(34){$1$};
\node at (2,7)(35){$8$};
\node at (2,5)(36){$52$};
\node at (2,3)(37){$313$};

\node at (1,8)(38){$1$};
\node at (1,6)(39){$9$};
\node at (1,4)(40){$63$};

\node at (0,7)(41){$1$};
\node at (0,5)(42){$10$};
\node at (0,3)(43){$75$};
\node at (-1,6)(44){$1$};
\node at (-1,4)(45){$11$};
\node at (-2,5)(46){$1$};
\node at (-2,3)(47){$12$};
\node at (-3,4)(48){$1$};
\node at (-4,3)(49){$1$};

\node at (-5,2)(50){$$};
\node at (-3,2)(51){$$};
\node at (-1,2)(52){$$};
\node at (1,2)(53){$$};
\node at (3,2)(54){$$};
\node at (5,2)(55){$$};
\node at (7,2)(56){$$};

\node at (-6,1)(57){};
\node at (-4,1)(58){};
\node at (-2,1)(59){};
\node at (0,1)(60){};
\node at (2,1)(61){};
\node at (4,1)(62){};
\node at (6,1)(63){};


\node at (-4,14.5)(a){$P_1$};
\node at (-4,13.5)(a){$P_2$};
\node at (-4,12.5)(a){$P_3$};
\node at (-4,11.5)(a){$P_4$};
\node at (-4,10.5)(a){$P_5$};
\node at (-4,9.5)(a){$P_6$};
\node at (-4,8.5)(a){$P_7$};
\node at (-4,7.5)(a){$P_8$};
\node at (-4,6.5)(a){$P_9$};
\node at (-4,5.5)(a){$P_{10}$};

\path[->] (1) edge (8);
\path[->] (8) edge (2);
\path[->] (2) edge (9);
\path[->] (9) edge (3);
\path[->] (3) edge (10);
\path[->] (10) edge (4);
\path[->] (4) edge (11);
\path[->] (11) edge (5);
\path[->] (5) edge (12);
\path[->] (12) edge (6);
\path[->] (6) edge (13);
\path[->] (13) edge (7);

\path[->] (8) edge (14);
\path[->] (14) edge (9);
\path[->] (9) edge (15);
\path[->] (15) edge (10);
\path[->] (10) edge (16);
\path[->] (16) edge (11);
\path[->] (11) edge (17);
\path[->] (17) edge (12);
\path[->] (12) edge (18);
\path[->] (18) edge (13);
\path[->] (13) edge (19);

\path[->] (14) edge (20);
\path[->] (20) edge (15);
\path[->] (15) edge (21);
\path[->] (21) edge (16);
\path[->] (16) edge (22);
\path[->] (22) edge (17);
\path[->] (17) edge (23);
\path[->] (23) edge (18);
\path[->] (18) edge (24);
\path[->] (24) edge (19);

\path[->] (20) edge (25);
\path[->] (25) edge (21);
\path[->] (21) edge (26);
\path[->] (26) edge (22);
\path[->] (22) edge (27);
\path[->] (27) edge (23);
\path[->] (23) edge (28);
\path[->] (28) edge (24);
\path[->] (24) edge (29);

\path[->] (25) edge (30);
\path[->] (30) edge (26);
\path[->] (26) edge (31);
\path[->] (31) edge (27);
\path[->] (27) edge (32);
\path[->] (32) edge (28);
\path[->] (28) edge (33);
\path[->] (33) edge (29);

\path[->] (30) edge (34);
\path[->] (34) edge (31);
\path[->] (31) edge (35);
\path[->] (35) edge (32);
\path[->] (32) edge (36);
\path[->] (36) edge (33);
\path[->] (33) edge (37);

\path[->] (34) edge (38);
\path[->] (38) edge (35);
\path[->] (35) edge (39);
\path[->] (39) edge (36);
\path[->] (36) edge (40);
\path[->] (40) edge (37);

\path[->] (38) edge (41);
\path[->] (41) edge (39);
\path[->] (39) edge (42);
\path[->] (42) edge (40);
\path[->] (40) edge (43);

\path[->] (41) edge (44);
\path[->] (44) edge (42);
\path[->] (42) edge (45);
\path[->] (45) edge (43);

\path[->] (44) edge (46);
\path[->] (46) edge (45);
\path[->] (45) edge (47);

\path[->] (46) edge (48);
\path[->] (48) edge (47);

\path[->] (48) edge (49);

\draw[line width=2pt] (1) -- (2);
\draw[line width=2pt] (2) -- (3);
\draw[line width=2pt] (3) -- (4);
\draw[line width=2pt] (4) -- (5);
\draw[line width=2pt] (5) -- (6);
\draw[line width=2pt] (6) -- (7);


\path[->] (50) edge (57);
\path[->] (51) edge (58);
\path[->] (52) edge (59);
\path[->] (53) edge (60);
\path[->] (54) edge (61);
\path[->] (55) edge (62);
\path[->] (56) edge (63);


\path[-] (7) edge (56);
\path[-] (19) edge (56);
\path[-] (19) edge (55);
\path[-] (29) edge (54);
\path[-] (29) edge (55);
\path[-] (37) edge (53);
\path[-] (37) edge (54);
\path[-] (43) edge (52);
\path[-] (43) edge (53);
\path[-] (47) edge (51);
\path[-] (47) edge (52);
\path[-] (49) edge (50);
\path[-] (49) edge (51);

\draw[dotted] (-4.5,15) -- (7,15);
\draw[dotted] (-4.5,14) -- (6,14);
\draw[dotted] (-4.5,13) -- (5,13);
\draw[dotted] (-4.5,12) -- (4,12);
\draw[dotted] (-4.5,11) -- (3,11);
\draw[dotted] (-4.5,10) -- (2,10);
\draw[dotted] (-4.5,9) -- (1,9);
\draw[dotted] (-4.5,8) -- (0,8);
\draw[dotted] (-4.5,7) -- (-1,7);
\draw[dotted] (-4.5,6) -- (-2,6);
\draw[dotted] (-4.5,5) -- (-3,5);

\node at (8,16)(a){$a_t[0]$};
\node at (7,16)(a){$a_t[1]$};
\node at (6,16)(a){$a_t[2]$};
\node at (5,16)(a){$a_t[3]$};
\node at (4,16)(a){$a_t[4]$};
\node at (3,16)(a){$a_t[5]$};

\node at (2,16)(a){$\cdots$};
\label{triangle}
\end{tikzpicture}
\end{equation*}

For example for $t=3$ and $t=4$, we compute $f_{8}$ and $f_{10}$ as follows:\par\bigskip
\addtocounter{equation}{1}
\begin{equation}
\begin{split}
21&=f_{8}=0+3(3.2^0)+0+1(3.2^2),\\
55&=f_{10}=1.7+0+4(3.2^1)+0+1(3.2^{3}).
\end{split}
\label{3}
\end{equation}

Sequences $a_{t}[0]=d_{2i,0}$ and $a_{t}[1]=d_{2i+1,1}$ are encoded respectively as A132262 and A110122 in the OEIS (On-Line Encyclopedia of Integer Sequences). Actually, sequence $a_{t}[0]$ had not been registered in the OEIS before the publication of Ringel and Fahr. 

\par\bigskip

In a second stop of the trip to a categorification of Fibonacci numbers, Ringel and Fahr generalized the results obtained in \cite{Fahr1}, and proved that the following exact sequences (\ref{categorification1}), (\ref{categorification2}) and filtration (\ref{categorification3}) are categorifications of identities (\ref{categorification4}). Where, for $t\geq1$, $P_{t}(x)$ and $R_{t}(x, y)$ are indecomposable representations of the quiver $Q=(T,E,\Omega^{x}_{t})$ ($\Omega^{x}_{t}$ is a bipartite orientation such that $x$ is a sink in case $t$ is even and a source in case $t$ is odd) for which $s_{t}(x)$ and $r_{t}(x, y)$ denote respectively their corresponding dimension vectors, assuming that for even $t$ the vertex $x$ is a sink, and that for $t$ odd the vertex $x$ is a source. In this setting, the sequences $x_{0}, x_{1},\dots, x_{t}$ and $x_{-1},x_{0},\dots, x_{t},x_{t+1}$ denote suitable paths with $x_{0}$ as a sink and $z_{i}$ being a neighbor of $x_{i}$ different from $x_{i-1}$ and $x_{i+1}$. 

\begin{equation}\label{categorification1}
\begin{split}
0&\rightarrow P_{t-1}(y)\rightarrow P_{t}(x)\rightarrow R_{t}(x, y)\rightarrow 0\\
0&\rightarrow P_{t-1}(y')\rightarrow R_{t}(x, y)\rightarrow R_{t-1}(y'',x)\rightarrow 0.
\end{split}
\end{equation}

\begin{equation}\label{categorification2}
\begin{split}
0&\rightarrow P_{0}(z_{0})\oplus \dots\oplus P_{t}(z_{t})\rightarrow R_{t+1}(x_{t}, x_{t+1})\rightarrow R_{0}(x_{-1},x_{0}).
\end{split}
\end{equation}

\begin{equation}\label{categorification3}
\begin{split}
P_{0}(x_{0})&\subset P_{1}(x_{1})\subset\dots\subset P_{t}(x_{t})\hspace{0.1cm}\text{with factors}\\
P_{i}(x_{i})/P_{i-1}(x_{i-1})&=R_{i}(x_{i},x_{i-1}),\quad 1\leq i\leq t.
\end{split}
\end{equation}

\begin{equation}\label{categorification4}
\begin{split}
f_{t+1}&=f_{t-1}+f_{t},\\
f_{2t+1}&=1+\underset{i=1}{\overset{t}{\sum}}f_{2i} \hspace{0.2cm}\text{and}\\
f_{2t}&=\underset{i=1}{\overset{t}{\sum}}f_{2i-1}.
\end{split}
\end{equation}

Note that the Auslander-Reiten sequences
\begin{equation}\label{categorification5}
\begin{split}
0&\rightarrow P_{n-1}\rightarrow P^{3}_{n}\rightarrow P_{n+1}\rightarrow 0\hspace{0.2cm}\text{and}\\
0&\rightarrow R_{n-1,\lambda}\rightarrow E(n,\lambda)\rightarrow R_{n+1,\lambda}\rightarrow 0
\end{split}
\end{equation}

where $E(n,\lambda)$ is an indecomposable module having dimension vector $3(\mathrm{dim} R(n,\lambda))$ are categorifications of the identity
\begin{equation}\label{categorification6}
\begin{split}
f_{t-2}+f_{t+2}&=3f_{t}.
\end{split}
\end{equation}

In a third stop of the road to a categorification of Fibonacci numbers Ringel and Fahr \cite{Fahr3} named the  array (\ref{farray}) a Fibonacci triangle and stated that its entries (nonzero entries) are categorified by the modules $P_{n}=P_{n}(x)$ (called Fibonacci modules) provided that such entries give the Jordan-H\"older multiplicities of these modules.\par\bigskip

Finally, we point out that Ringel in \cite{Ringel} exhibits combinatorial data which can be derived from a category $\mathrm{Ind}\hspace{0.1cm}\Lambda$, where $\Lambda$ is a hereditary artin algebra of Dynkin type $\Delta$ and $\mathrm{Ind}\hspace{0.1cm}\Lambda$ is a set of indecomposable $\Lambda$-modules. He comments that many enumeration problems give rise to categorification of different integer sequences. For instance, the number of some tilting modules and the number of antichains in $\mathrm{mod}\hspace{0.1cm}\Lambda$ categorify the Catalan numbers if $\Lambda$ is an algebra of Dynkin type $\mathbb{A}_{n}$. Whereas, if $\Lambda$ is of Dynkin type $\mathbb{B}_{n}$ then such number of modules and antichains categorify the sequence $\binom{2n}{n}$.\par\bigskip 
According to Ringel \cite{Ringel}, the number of antichains in $\mathrm{mod}\hspace{0.1cm}\Lambda$ as well as the number of tilting modules and the number of indecomposable $\Lambda$-modules are examples of Dynkin functions which do not depend on the orientation. Regarding, this particular kind of functions he proposes to build an On-Line Encyclopedia of Dynkin Functions (OEDF) with the same purposes as the OEIS.\par\bigskip

In this work, in order to categorify integer sequences, we identify combinatorial information arising from the preprojective components of the 2-Kronecker algebra (or simply the Kronecker algebra) and the tetrad (or four incomparable points) with combinatorial information arising from indecomposable projective modules over some Brauer configuration algebras introduced recently by Green and Schroll \cite{Green}. In particular, we use these settings to define categorifications of the sequences encoded in the OEIS as A052558, A052591 and A100705. Such Brauer configuration algebras are defined by configurations of some multisets called polygons.
\par\bigskip
We recall here  that the Kronecker problem is equivalent to the problem of determining the indecomposable representations over a field $k$ of the following quiver $Q$ (2-Kronecker quiver):

\begin{equation}\label{QKronecker}
Q=\xymatrix@=30pt{
\overset{0}{\circ}&\ar@/_/[l]_{\alpha}\ar@/^/[l]^{\beta}\overset{1}{\circ}
}
\end{equation}

whereas to determine the indecomposable representations of four incomparable points (a tetrad) is  a very well known matrix problem named the four subspace problem (FSP). The solution of this problem is essentially equivalent to determine all of the indecomposable representations of the four subspace quiver $F$ with the following form \cites{Gelfand, Nazarova2, Simson, Zavadskij1}:

\begin{equation}\label{F}
 \xymatrix@=30pt{
&&\overset{1}{\circ}\ar@[black][d]&&\\
F=&\underset{2}{\circ}\ar@[black][r]&\underset{5}{\circ}&\underset{4}{\circ}\ar@[black][l]&\\
&&\underset{3}{\circ}\ar@[black][u]&&&
}
\end{equation}

We will see that some invariants associated to indecomposable projective modules over some suitable Brauer configuration algebras allow categorifications of integer numbers. In fact, since polygons in Brauer configurations are multisets, we will often assume that such polygons consists of words of the form

\begin{equation}\label{word}
\begin{split}
w&=x_{1}^{s_1}x_{2}^{s_2}\dots x^{s_{t-1}}_{t-1}x_{t}^{s_{t}}
\end{split}
\end{equation}

where for each $i$, $1\leq i\leq t$, $x_{i}$ is an element of the polygon called vertex and $s_{i}$ is the number of times that the vertex $x_{i}$ occurs in the polygon. In particular, if vertices $x_{i}$ in a polygon $V$ of a Brauer configuration are integer numbers then the corresponding word $w$ will be interpreted as a partition of an integer number $n_{V}$ associated to the polygon $V$ where it is assumed that each vertex $x_{i}$ is a part of the partition and $s_{i}$ is the number of times that the part $x_{i}$ occurs in the partition and $n_{V}=\underset{i=1}{\overset{t}{\sum}}s_{i}x_{i}$.\par\bigskip

In this paper, we prove the following results:

\begin{teorema}\label{a}
\textit{Let $(P,A,B,n)$, $(P',A',B',n)$, $H_{P}$ and $H_{P'}$ be two matrices of type $\mathscr{H}_{n}$ with corresponding sets of helices $H_{P}$ and $H_{P'}$ defined by systems of the form $(i_{P},j_{P}, P_{A}, P_{B})$ and $(f_{P'}, g_{P'}, P'_{A'}, P'_{B'})$, respectively. Then $|H_{P}|=h^{P}_{n}=|H_{P'}|=h^{P'}_{n}$}.
\end{teorema}
Helices associated to suitable $(n+1)\times 2n$-matrices of type $\mathscr{H}_{n}$ are introduced in section \ref{Kr2}.

\addtocounter{corolario}{1}
\begin{corolario}
\textit{If for $n\geq 1$, $P$ and $P'$ are equivalent preprojective Kronecker modules with dimension vector of the form $[n+1\hspace{0.2cm}n]$ and corresponding sets of helices $H_{P}$ and $H_{P'}$ then $|H_{P}|=|H_{P'}|$.}
\end{corolario}

\addtocounter{teorema}{1}
\begin{teorema}\label{TKr}
\textit{If for $n\geq1$, $P$ denotes a preprojective Kronecker module then the number of helices associated to $P$ is $h^{P}_{n}=n!\lceil\frac{n}{2}\rceil$ where $\lceil x \rceil$ denotes the smallest integer greatest than $x$}.

\end{teorema}

For a suitable Brauer configuration algebra $\Lambda_{K^{n}}$, the following result categorifies in the sense of Ringel and Fahr the number of helices associated to some preprojective Kronecker modules.

\addtocounter{corolario}{1}

\begin{corolario}\label{CKr}
\textit{For $n\geq3$ fixed and $1\leq t\leq n$, the number of summands in the heart of the indecomposable projective representation $V_{t}$ over the Brauer configuration algebra $\Lambda_{K^n}$ equals the number of helices associated to the preprojective Kronecker module $(2t+3,2t+2)$, $1\leq t\leq n$.}
\end{corolario}
\addtocounter{proposicion}{4}
\begin{proposicion}\label{word partition01}
\textit{If $\mathscr{W}_{P}$ is the set of matrix words associated to a matrix $P$ of type $\mathscr{H}_{n}$ then $|\mathscr{W}_{P}|$ equals $\underset{m=0}{\overset{n^{2}}{\sum}}P(n, n, m)=(n+1)C_{n}$, where $P(n, n, m)$ denotes the number of partitions of $m$ into $n$ parts, each $\leq n$, $P(n, n, 0)=1$, and $C_{n}$ denotes the $n$th Catalan number.}
\end{proposicion}

\addtocounter{teorema}{2}

The next theorem regards the number of cycles associated to preprojective representations of the tetrad and its relationship with the number of summands in the heart of an indecomposable module over a suitable Brauer configuration algebra. The notion of cycle associated to suitable $(2n+2)\times (4n+3)$-matrices of type $\mathscr{C}_{n}$ is introduced in section \ref{FSP01}.

\begin{teorema}\label{Tprojective quadruple}

\textit{For $j\geq2$ fixed and $1\leq i\leq j$, the number of summands in the heart of the indecomposable projective module $T_{P_{i}}$ over the Brauer configuration algebra $\Lambda_{E^{j}}$ equals the number of cycles associated to a preprojective representation of the tetrad of type $\mathrm{IV}$ and order $i+1$.}

\end{teorema}

The following result gives a way to categorify any counting function $u_{n}$, i.e., sequences whose elements count a given class of objects. Catalan numbers, Fibonacci numbers, Delannoy numbers and Dedekind numbers are some of the most well known counting functions. $D_{n}$ is a suitable Brauer configuration.

\begin{teorema}\label{Tgral1}
\textit{For $1\leq i\leq n$ and $n>1$ fixed, the number of summands in the heart of the indecomposable projective module $P_{i}$ over the Brauer configuration algebra $\Lambda_{D^{n}}$ is $u_{i}$.}

\end{teorema}

Dimension of the Brauer configuration algebras involved in the different results and their corresponding centers are given in Corollaries \ref{Kronecker2}, \ref{Kronecker3}, \ref{tetrad3}, \ref{tetrad4} and Theorems \ref{tetradnew1}, \ref{tetradnew2}, \ref{gral2} and \ref{gral3}. 
\par\bigskip

This paper is distributed as follows: In section 2, we recall main definitions and notation used throughout the document, in particular, in this section we define Brauer configuration algebras, the Kronecker problem and the four subspace problem.\par\bigskip In section 3, in order to categorify numbers in sequences A052558 and A052591, we prove that numbers in these sequences give the number of some helices associated to preprojective Kronecker modules. Besides, it is defined a sequence $\Lambda_{K^{n}}$ of Brauer configuration algebras whose indecomposable projective modules are in bijective correspondence with preprojective Kronecker modules via the number of summands in the heart of such indecomposable modules. Formulas for the dimension of this type of algebras and corresponding centers are given as well.
\par\bigskip
In section 4, numbers in the sequence A100705  are categorified by determining the number of cycles associated to preprojective representations of type IV of the tetrad and it is introduced a sequence $\Lambda_{E^{n}}$ of Brauer configuration algebras such that the number of summands in the heart of their indecomposable projective modules coincides with the number of cycles associated to the mentioned preprojective representations. The dimension of these algebras and corresponding quotients are also obtained in this section. \par\bigskip

In section 5, we describe how it is possible to use integer sequences in order to build Brauer configuration algebras, the process is applied to any counting function and in particular to the sequence A100705. Examples of helices are given in section \ref{Appendix}.

\section{Preliminaries}

In this section, we recall main definitions and notation to be used throughout the paper~\cites{Andrews, Canadas3, Green, Zavadskij1}.

\subsection{Brauer Configuration Algebras}

Brauer configuration algebras were introduced by Green and Schroll in \cite{Green} as a generalization of Brauer graph algebras which are biserial algebras of tame representation type and whose representation theory is encoded by some combinatorial data based on graphs. According to them, underlying every Brauer graph algebra is a finite graph with a cyclic orientation of the edges at every vertex and a multiplicity function \cite{Schroll}.  The construction of a Brauer graph algebra is a special case of the construction of a Brauer configuration algebra in the sense that every  Brauer graph is a Brauer configuration with the restriction that every polygon is a set with two vertices. In the sequel, we give precise definitions of a Brauer configuration and a Brauer configuration algebra.\par\bigskip

A \textit{Brauer configuration} $\Gamma$ is a quadruple of the form $\Gamma=(\Gamma_{0},\Gamma_{1},\mu,\mathcal{O})$ where:

\begin{enumerate}[$(B1)$]
\item $\Gamma_{0}$ is a finite set whose elements are called \textit{vertices},
\item  $\Gamma_{1}$ is a finite collection of multisets called \textit{polygons}. In this case, if $V\in \Gamma_{1}$ then the elements of $V$ are vertices possibly with repetitions,  $\mathrm{occ}(\alpha,V)$ denotes the frequency of the vertex $\alpha$ in the polygon $V$ and the \textit{valency} of $\alpha$ denoted  $val(\alpha)$ is defined  in such a way that:

\begin{equation}
\begin{split}
val(\alpha)&=\underset{V\in\Gamma_{1}}{\sum}\mathrm{occ}(\alpha,V),
\end{split}
\end{equation}

\item $\mu$ is an integer valued function such that $\mu:\Gamma_{0}\rightarrow \mathbb{N}$ where $\mathbb{N}$ denotes the set of positive integers, it is called the \textit{multiplicity function},
\item $\mathcal{O}$ denotes an orientation defined on $\Gamma_{1}$ which is a choice, for each vertex $\alpha \in \Gamma_0$, of a cyclic ordering of the polygons in which $\alpha$ occurs as a vertex, including repetitions, we denote $S_{\alpha}$ such collection of polygons.  More specifically, if $S_{\alpha}=\{V^{(\alpha_{1})}_{1},V^{(\alpha_{2})}_{2},\dots, V^{(\alpha_t)}_{t}\}$ is the collection of polygons where the vertex $\alpha$ occurs with $\alpha_{i}=\mathrm{occ}(\alpha,V_{i})$ and $V^{(\alpha_{i})}_{i}$ meaning that $S_{\alpha}$ has $\alpha_{i}$ copies of $V_{i}$ then an orientation $\mathcal{O}$ is obtained by endowing a linear order $<$ to $S_{\alpha}$  and adding a relation $V_{t}< V_{1}$, if $V_{1}=\mathrm{min}\hspace{0.1cm}S_{\alpha}$ and $V_{t}=\mathrm{max}\hspace{0.1cm}S_{\alpha}$. According to this order the $\alpha_{i}$ copies of $V_{i}$ can be ordered as $V_{1, i}<V_{2,i}<\dots < V_{(\alpha_{i}-1),i}<V_{\alpha_{i},i}$ and $S_{\alpha}$ can be ordered in the form $V^{(\alpha_{1})}_{1}<V^{(\alpha_{2})}_{2}<\dots< V^{(\alpha_{(t-1)})}_{(t-1)}<V^{\alpha_{t}}_{t}$,
\item Every vertex in $\Gamma_{0}$ is a vertex in at least one polygon in $\Gamma_{1}$,
\item Every polygon has at least two vertices,
\item Every polygon in $\Gamma_{1}$ has at least one vertex $\alpha$ such that $\mu(\alpha)val(\alpha)>1$.
\end{enumerate}

The set $(S_{\alpha},<)$ is called the \textit{successor sequence} at the vertex $\alpha$.\par\bigskip

A vertex $\alpha\in\Gamma_{0}$ is said to be \textit{truncated} if $val(\alpha)\mu(\alpha)=1$, that is, $\alpha$ is truncated if it occurs exactly once in exactly one $V\in\Gamma_{1}$ and $\mu(\alpha)=1$. A vertex is \textit{nontruncated} if it is not truncated.
\par\bigskip

\addtocounter{Nota}{7}

\begin{center}
\textbf{The Quiver of a Brauer Configuration Algebra}
\end{center}

The quiver $Q_{\Gamma}=((Q_{\Gamma})_{0},(Q_{\Gamma})_{1})$ of a Brauer configuration algebra is defined in such a way that the vertex set $(Q_{\Gamma})_{0}=\{v_{1},v_{2},\dots, v_{m}\}$ of $Q_{\Gamma}$ is in correspondence with the set of polygons $\{V_{1},V_{2},\dots,V_{m}\}$ in $\Gamma_{1}$, noting that there is one vertex in $(Q_{\Gamma})_{0}$ for every polygon in $\Gamma_{1}$.\par\bigskip

Arrows in $Q_{\Gamma}$ are defined by the successor sequences. That is, there is an arrow $v_{i}\stackrel{s_{i}}{\longrightarrow}v_{i+1}\in (Q_{\Gamma})_{1}$ provided that $V_{i}<V_{i+1}$ in $(S_{\alpha},<)\cup\{V_{t}<V_{1}\}$ for some nontruncated vertex $\alpha\in\Gamma_{0}$. In other words, for each nontruncated vertex $\alpha\in\Gamma_{0}$ and each successor $V'$ of $V$ at $\alpha$, there is an arrow from $v$ to $v'$ in $Q_{\Gamma}$ where $v$ and $v'$ are the vertices in $Q_{\Gamma}$ associated to the polygons $V$ and $V'$ in $\Gamma_{1}$, respectively.\par\bigskip

\begin{centering}

\textbf{The Ideal of Relations and Definition of a Brauer Configuration Algebra}\par\bigskip

\end{centering}

Fix a polygon $V\in\Gamma_{1}$ and suppose that $\mathrm{occ}(\alpha,V)=t\geq1$ then there are $t$ indices
$i_{1},\dots, i_{t}$ such that $V=V_{i_{j}}$. Then the \textit{special $\alpha$-cycles} at $v$ are the cycles $C_{i_{1}}, C_{i_{2}},\dots, C_{i_{t}}$ where $v$ is the vertex in the quiver of $Q_{\Gamma}$ associated to the polygon $V$.
If $\alpha$ occurs only once in $V$ and $\mu(\alpha)=1$ then there is only one special $\alpha$-cycle at $v$.

\par\bigskip

Let $k$ be a field and $\Gamma$ a Brauer configuration. The \textit{Brauer configuration algebra associated to $\Gamma$} is defined to be the bounded path algebra $\Lambda_{\Gamma}=kQ_{\Gamma}/I_{\Gamma}$, where $Q_{\Gamma}$ is the quiver associated to $\Gamma$ and $I_{\Gamma}$ is the \textit{ideal} in $kQ_{\Gamma}$ generated by the following set of relations $\rho_{\Gamma}$ of type I, II and III.\par\bigskip
\begin{enumerate}
\item \textbf{Relations of type I}. For each polygon $V=\{\alpha_{1},\dots, \alpha_{m}\}\in \Gamma_{1}$ and each pair of nontruncated vertices $\alpha_{i}$ and $\alpha_{j}$ in $V$, the set of relations $\rho_{\Gamma}$ contains all relations of the form $C^{\mu(\alpha_{i})}-C'^{\mu(\alpha_{j})}$ where $C$ is a special $\alpha_{i}$-cycle
and $C'$ is a special $\alpha_{j}$-cycle.

\item \textbf{Relations of type II}. Relations of type II are all paths of the form $C^{\mu(\alpha)}a$ where $C$ is a special $\alpha$-cycle and $a$ is the first arrow in $C$.

\item \textbf{Relations of type III}. These relations are quadratic monomial relations of the form $ab$ in $kQ_{\Gamma}$ where $ab$ is not a subpath of any special cycle unless $a=b$ and $a$ is a loop associated to a vertex of valency 1 and $\mu(\alpha)>1$.

\end{enumerate}

Henceforth, if there is no confusion, we will assume notations, $\Lambda$, $I$ and $\rho$ instead of $\Lambda_{\Gamma}$, $I_{\Gamma}$ and $\rho_{\Gamma}$ for a Brauer configuration algebra, the ideal and set of relations, respectively defined by a given Brauer configuration $\Gamma$.\par\bigskip

As a toy example consider a configuration $\Gamma=(\Gamma_{0},\Gamma_{1},\mu,\mathcal{O})$ such that:
\begin{enumerate}\label{example}
\item $\Gamma_{0}=\{1,2,3,4\}$,
\item $\Gamma_{1}=\{U=\{1,1,2,3,3,4\}, V=\{1,2,3,4,4\}\}$,
\item At vertex 1, it holds that; $S_{1}=\{U^{(2)}V^{(1)}\}$,\quad $U< U< V$,\quad $val(1)=3$,
\item At vertex 2, it holds that;  $S_{2}=\{U^{(1)}V^{(1)}\}$,\quad $U<V$,\quad $val(2)=2$,
\item At vertex 3, it holds that;  $S_{3}=\{U^{(2)}V^{(1)}\}$,\quad $U< U<V$,\quad$val(3)=3$,
\item At vertex 4, it holds that;  $S_{4}=\{U^{(1)}V^{(2)}\}$,\quad$U< V< V$,\quad$val(4)=3$,
\item $\mu(\alpha)=1$ for any vertex $\alpha$.

\end{enumerate}

The ideal $\mathrm{I}$ of the corresponding Brauer configuration algebra $\Lambda_{\Gamma}$ is generated by the following relations (see Figure (\ref{fexample})), for which it is assumed the following notation for the special cycles:
\begin{equation}\label{special}
\begin{split}
P^{u,1}_{1}&=c^{1}_{1}c^{1}_{2}c^{1}_{3},\\
P^{u,1}_{2}&=c^{1}_{2}c^{1}_{3}c^{1}_{1},\\
P^{u,2}_{3}&=c^{2}_{1}c^{2}_{2},\\
P^{u,3}_{4}&=c^{3}_{2}c^{3}_{3}c^{3}_{1},\\
P^{u,3}_{5}&=c^{3}_{1}c^{3}_{2}c^{3}_{3},\\
P^{u,4}_{6}&=c^{4}_{1}c^{4}_{2}c^{4}_{3},\\
P^{v,1}_{1}&=c^{1}_{3}c^{1}_{1}c^{1}_{2},\\
P^{v,2}_{2}&=c^{2}_{2}c^{2}_{1},\\
P^{v,3}_{3}&=c^{3}_{3}c^{3}_{1}c^{3}_{2},\\
P^{v,4}_{4}&=c^{4}_{2}c^{4}_{3}c^{4}_{1},\\
P^{v,4}_{5}&=c^{4}_{3}c^{4}_{1}c^{4}_{2},\\
\end{split}
\end{equation}

\begin{enumerate}
\item $c^{h}_{i}c^{s}_{r}$, if $h\neq s$, for all possible values of $i$ and $r$,
\item $(c^{1}_{1})^{2}$;\quad$(c^{3}_{1})^{2}$;\quad $(c^{4}_{2})^{2}$, $c^{4}_1c^4_3$,\quad$c^{1}_{3}c^1_{2}$,\quad$c^3_{3}c^3_2$,
\item $P^{u, i}_{j}-P^{u, t}_{l}$,\quad for all possible values of $i,j,t$ and $l$,
\item $P^{v, i}_{j}-P^{v, t}_{l}$,\quad for all possible values of $i,j,t$ and $l$,
\item $P^{u, j}_{i}a$ ($P^{v, j}_{i}a'$) , with $a$ ($a'$) being the first arrow of $P^{u,j}_{i}$ ($P^{v,j}_{i}$) for all $i,j$.
\end{enumerate}

The following diagrams (\ref{fexample}-\ref{Pv3}) show the quiver $Q_{\Gamma}$ associated to this configuration, the indecomposable projective modules $P_{U}$ and $P_{V}$, corresponding heart and radical square of these modules.

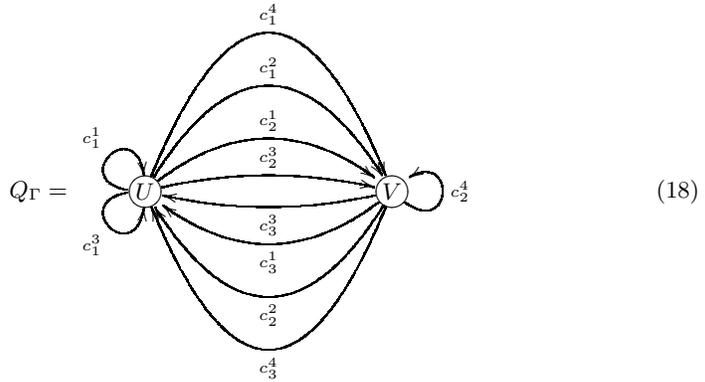
\begin{equation}\label{fexample}
Q_{\Gamma}=\xymatrix@=80pt{*+[o][F-]{U}\ar@(l,d)_{c_1^3}\ar@(l,u)@[blue]^{c_1^1}\ar@/^/[r]^{c^3_2}\ar@/^20pt/@[blue][r]^{c^1_2}\ar@/^40pt/@[red][r]^{c^2_1}\ar@/^60pt/@[green][r]^{c^4_1}&*+[o][F-]{V}\ar@(dr,ru)@[green]_{c_2^4}\ar@/^/[l]^{c^3_3}\ar@/^20pt/@[blue][l]^{c^1_3}\ar@/^40pt/@[red][l]^{c^2_2}\ar@/^60pt/@[green][l]^{c^4_3}}
\end{equation}

\begin{equation}\label{Pu1}
P(U): 
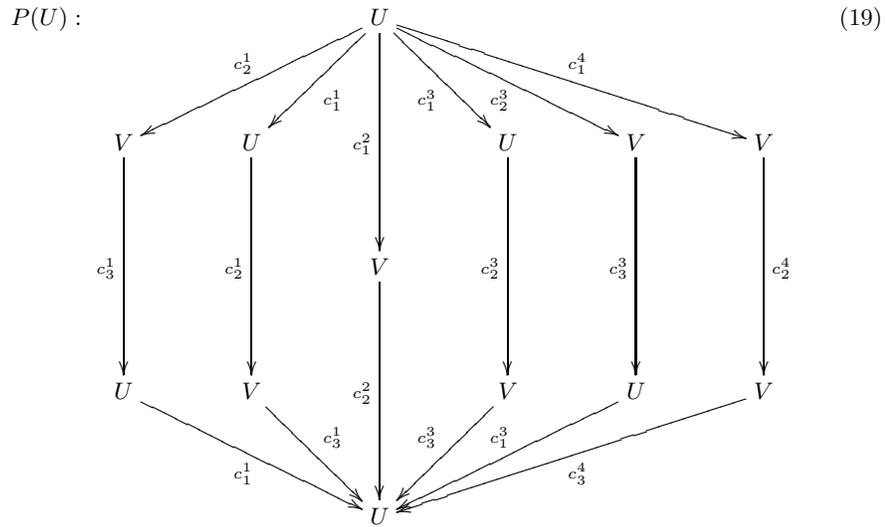
\xymatrix@=35pt{
&&\ar@[blue][lld]_{c_2^1}\ar@[blue][ld]^{{c_1^1}}\ar@[red][dd]_{c_1^2}\ar@[black][dr]_{c_1^3}\ar@[black][drr]_{c_2^3}\ar@[green][drrr]^{c_1^4}U&&&\\
\ar@[blue][dd]_{c_3^1}V&\ar@[blue][dd]_{c_2^1}U&&\ar@[black][dd]_{c_2^3}U&\ar@[black][dd]_{c_3^3}V&\ar@[green][dd]^{c_2^4}V\\
&&\ar@[red][dd]_{c_2^2}V&&&\\
\ar@[blue][drr]_{c_1^1}U&\ar@[blue][dr]^{{c_3^1}}V&&\ar@[black][dl]_{c_3^3}V&\ar@[black][dll]_{c_1^3}U&\ar@[green][dlll]^{c_3^4}V\\
&&U&&&
}
\end{equation}

\begin{equation}\label{Pu2}
\mathrm{Heart}(P(U)): 
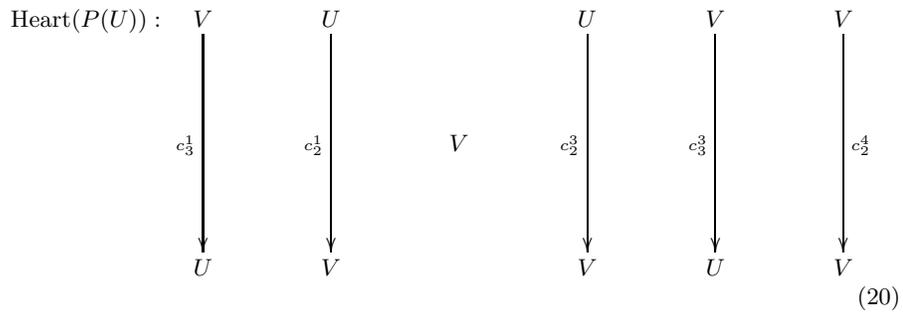
\xymatrix@=35pt{
\ar@[blue][dd]_{c_3^1}V&\ar@[blue][dd]_{c_2^1}U&&\ar@[black][dd]_{c_2^3}U&\ar@[black][dd]_{c_3^3}V&\ar@[green][dd]^{c_2^4}V\\
&&V&&&\\
U&V&&V&U&V\\
}
\end{equation}

\begin{equation}\label{Pu3}
\mathrm{rad}^2P(U): 
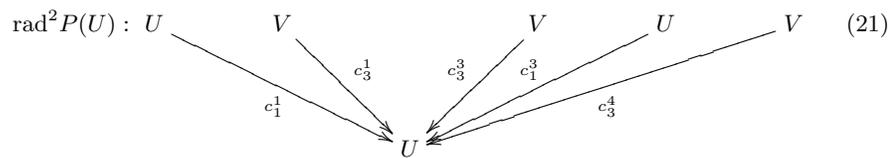
\xymatrix@=35pt{
\ar@[blue][drr]_{c_1^1}U&\ar@[blue][dr]^{{c_3^1}}V&&\ar@[black][dl]_{c_3^3}V&\ar@[black][dll]_{c_1^3}U&\ar@[green][dlll]^{c_3^4}V\\
&&U&&&
}
\end{equation}

\begin{equation}\label{Pv1}
P(V): 
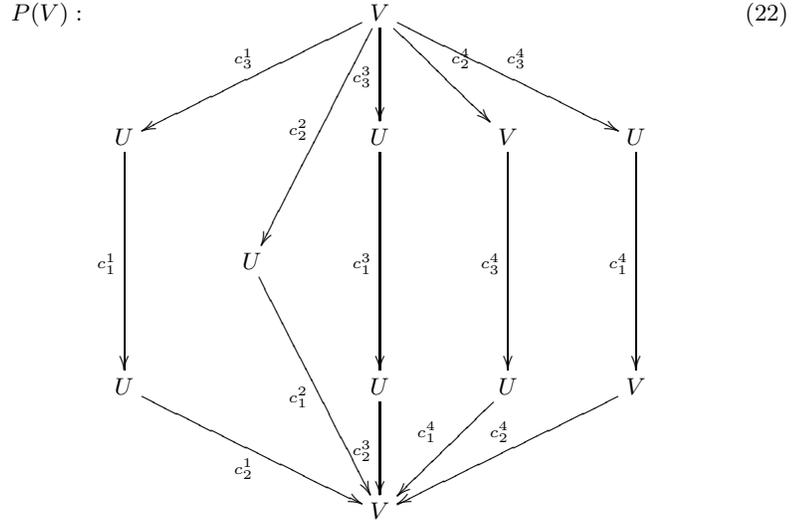
\xymatrix@=35pt{
&&\ar@[blue][lld]_{c_3^1}\ar@[red][ldd]_{c_2^2}\ar@[black][d]_{c_3^3}\ar@[green][dr]^{c_2^4}\ar@[green][drr]^{c_3^4}V&&\\
\ar@[blue][dd]_{c_1^1}U&&\ar@[black][dd]_{c_1^3}U&\ar@[green][dd]_{c_3^4}V&\ar@[green][dd]_{c_1^4}U\\
&\ar@[red][rdd]_{c_1^2}U&&&\\
\ar@[blue][drr]_{c_2^1}U&&\ar@[black][d]_{c_2^3}U&\ar@[green][dl]_{c_1^4}U&\ar@[green][dll]_{c_2^4}V\\
&&V&&\\
}
\end{equation}

\begin{equation}\label{Pv2}
\mathrm{Heart}\,P(V): 
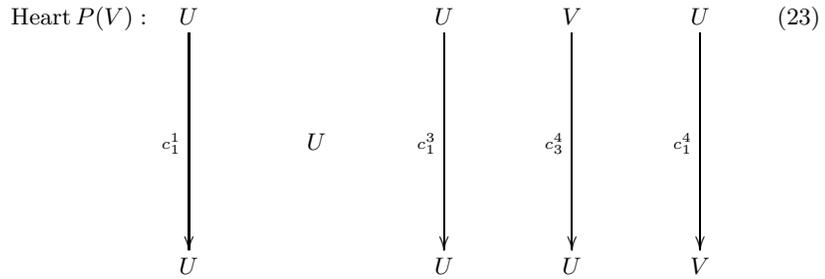
\xymatrix@=35pt{
\ar@[blue][dd]_{c_1^1}U&&\ar@[black][dd]_{c_1^3}U&\ar@[green][dd]_{c_3^4}V&\ar@[green][dd]_{c_1^4}U\\
&U&&&\\
U&&U&U&V\\
}
\end{equation}

\begin{equation}\label{Pv3}
\mathrm{rad}^2\, P(V): 
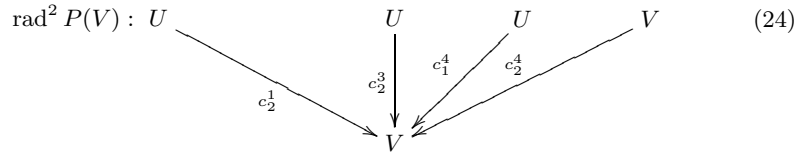
\xymatrix@=35pt{
\ar@[blue][drr]_{c_2^1}U&&\ar@[black][d]_{c_2^3}U&\ar@[green][dl]_{c_1^4}U&\ar@[green][dll]_{c_2^4}V\\
&&V&&\\
}
\end{equation}

\par\bigskip
The following results give some description of the structure of Brauer configuration algebras \cites{Green, Sierra}.

\begin{teor}[\cite{Green}, Theorem B, Proposition 2.7, Theorem 3.10, Corollary 3.12]\label{multiserial}
	\textit{Let $\Lambda$ be a Brauer configuration algebra with Brauer configuration $\Gamma$}.
	\begin{enumerate}
		\item \textit{There is a bijective correspondence between the set of indecomposable projective $\Lambda$-modules and the polygons in $\Gamma$}.
		\item \textit{If $P$ is an indecomposable projective $\Lambda$-module corresponding to a polygon $V$ in $\Gamma$. Then $\mathrm{rad}\hspace{0.1cm}P$ is a sum of $r$ indecomposable uniserial modules, where $r$ is the number
			of (nontruncated) vertices of $V$ and where the intersection of any two of the uniserial modules is a simple $\Lambda$-module}.
		\item \textit{A Brauer configuration algebra is a multiserial algebra}.
		\item \textit{The number of summands in the heart $ht(P)=\mathrm{rad}\hspace{0.1cm}P/\mathrm{soc}\hspace{0.1cm}P$ of an indecomposable projective $\Lambda$-module $P$ such that $\mathrm{rad}^{2}\hspace{0.1cm}P\neq 0$ equals the number of nontruncated vertices of the polygons in $\Gamma$ corresponding to $P$ counting repetitions}.
		\item \textit{If $\Lambda'$ is a Brauer configuration algebra obtained from $\Lambda$ by removing a truncated vertex of a polygon in $\Gamma_{1}$ with $d\geq 3$ vertices then $\Lambda$ is isomorphic to $\Lambda'$.}
		
	\end{enumerate}
	
\end{teor}

\addtocounter{prop}{1}

\begin{prop}[\cite{Green}, Proposition 3.3]\label{basis}
	\textit{Let $\Lambda$ be the Brauer configuration algebra associated to the Brauer configuration $\Gamma$. For each $V\in\Gamma_{1}$ choose a nontruncated vertex $\alpha$ and exactly one special $\alpha$-cycle $C_{V}$ at $V$,}	
\begin{equation*}
\begin{split}
A&=\{\overline{p}\mid p\hspace{0.1cm} is\hspace{0.1cm} a\hspace{0.1cm} proper\hspace{0.1cm} prefix\hspace{0.1cm} of\hspace{0.1cm} some\hspace{0.1cm}C^{\mu(\alpha)}\hspace{0.1cm} where\hspace{0.1cm} C\hspace{0.1cm} is\hspace{0.1cm} a\hspace{0.1cm} special\hspace{0.1cm}\alpha-cycle\},\\
B&= \{\overline{C^{\mu(\alpha)}_{V}}\mid V\in \Gamma_{1}\}.
\end{split}
\end{equation*}
	 \textit{Then $A\cup B$ is a $k$-basis of $\Lambda$}.
\end{prop}

\begin{prop}[\cite{Green}, Proposition 3.13]\label{dimension}
	\textit{Let $\Lambda$ be a Brauer configuration algebra associated to the Brauer configuration $\Gamma$ and let $\mathcal{C}=\{C_{1},\dots, C_{t}\}$ be a full set of equivalence class representatives of special cycles. Assume that for $i=1,\dots,t$, $C_{i}$ is a special $\alpha_{i}$-cycle where $\alpha_{i}$ is a nontruncated vertex in $\Gamma$}. \textit{Then}
	\begin{center}
		$\mathrm{dim}_{k}\hspace{0.1cm}\Lambda=2|Q_{0}|+\underset{C_{i}\in\mathcal{C}}{\sum}|C_{i}|(n_{i}|C_{i}|-1)$,
		
	\end{center}
	\textit{where $|Q_{0}|$ denotes the number of vertices of $Q$, $|C_{i}|$ denotes the number of arrows in the $\alpha_{i}$-cycle $C_{i}$ and $n_{i}=\mu(\alpha_{i})$}.
\end{prop}

The following result regards the center of a Brauer configuration algebra.

\addtocounter{teor}{2}
\begin{teor}[\cite{Sierra}, Theorem 4.9]\label{Serra}
	\textit{Let $\Gamma$ be a reduced and connected Brauer configuration and let $Q$ be its induced quiver and let $\Lambda$ be the induced Brauer configuration algebra such that $\mathrm{rad}^{2}\hspace{0.1cm}\Lambda \neq 0$ then the dimension of the center of $\Lambda$ denoted $\mathrm{dim}_{k}\hspace{0.1cm}Z(\Lambda)$ is given by the formula}:
	
	\begin{equation*}\label{Sierra}
	\begin{split}
	\mathrm{dim}_{k}\hspace{0.1cm}Z(\Lambda)&=1+\underset{\alpha\in\Gamma_{0}}{\sum}\mu(\alpha)+|\Gamma_{1}|-|\Gamma_{0}|+\#(Loops\hspace{0.1cm} Q)-|\mathscr{C}_{\Gamma}|.
	\end{split}
	\end{equation*}
	
	\textit{where} $|\mathscr{C}_{\Gamma}|=\{\alpha\in\Gamma_{0}\mid val(\alpha)=1, \hspace{0.1cm}and\hspace{0.1cm} \mu(\alpha)>1\}$.
	
\end{teor}

As an example the following are the numerical data associated to the algebra $\Lambda_{\Gamma}=kQ_{\Gamma}/I$ with $Q_{\Gamma}$ as shown in Figure (\ref{fexample}) and special cycles given in (\ref{special}), ($|r(Q_{\Gamma})|$ is the number of indecomposable projective modules, $r_{U}$ and $r_{V}$ denote the number of summands in the heart of the indecomposable projective modules $P(U)$ and $P(V)$). Note that, $|C_{i}|=val(i)$:

\begin{equation}\label{theexample}
\begin{split}
|r(Q_{\Gamma})|&=2,\\
r_{U}&=6,\quad r_{V}=5,\\
|C_{1}|&=3, \quad|C_{2}|=2,\quad|C_{3}|=3,\quad|C_{4}|=3,\\
\underset{\alpha\in\Gamma_{0}}{\sum}\underset{X\in\Gamma_{1}}{\sum}\mathrm{occ}(\alpha, X)&= 11,\quad\text{the number of special cycles},\\
\mathrm{dim}_{k}\hspace{0.1cm}\Lambda_{\Gamma}&=4+3(3-1)+3(3-1)+2(2-1)+3(3-1)=24,\\
\mathrm{dim}_{k}\hspace{0.1cm}Z(\Lambda_\Gamma)&=1+4+(2-4)+3=6. 
\end{split}
\end{equation}

\subsection{The Kronecker Problem}\label{Kr2}

The classification of indecomposable Kronecker modules was solved by Weierstrass in 1867 for some particular cases and by Kronecker in 1890 for the complex number field case. This problem is equivalent to the problem of finding canonical Jordan form of pairs of matrices $(A,B)$ (with the same size) with respect to the following elementary transformations over a field $k$ (for the sake of brevity, it is assumed that $k$ is an algebraically closed field):
\begin{enumerate}
    \item [(i)] All elementary transformations on rows of the block matrix $(A,B)$.
    \item[(ii)] All elementary transformations made simultaneously on columns of $A$ and $B$ having the same index number.  
\end{enumerate}

If the matrix blocks $P=(A, B)$ and $P'=(A', B')$ can be transformed one into the other by means of elementary transformations, then they are said to be \textit{equivalent} or isomorphic as Kronecker modules. The following is the matrix form (up to isomorphism) of the non-regular Kronecker modules \cites{Simson, Zavadskij2}:

\begin{centering}
 \setlength{\unitlength}{1pt}
\begin{picture}(100,100)

      \multiput(0,50)(30,0){3}{\line(0,1){30}}
        \multiput(0,50)(0,30){2}{\line(1,0){60}}
        \put(-45,65){$\mathrm{II}=\mathrm{III^{*}}$: }
     \put(10,65){$\overset{\rightarrow}{\mathrm{I}_{n}}$}
      \put(35,65){$\overset{\leftarrow}{\mathrm{I}_{n}}$}

       \multiput(0,0)(30,0){3}{\line(0,1){30}}
        \multiput(0,0)(0,30){2}{\line(1,0){60}}
        \put(-45,15){$\mathrm{III}=\mathrm{II^{*}}$: }
     \put(10,15){$\mathrm{I}^{\uparrow}_{n}$}
      \put(40,15){$\mathrm{I}^{\downarrow}_{n}$}

\end{picture}\par\bigskip

\end{centering}

In this case, $\overset{\rightarrow}{\mathrm{I}_{n}}$ ($\overset{\leftarrow}{\mathrm{I}_{n}}$, respectively)  denotes an $n\times (n+1)$ matrix obtained from the identity $\mathrm{I}_{n}$ by adding a column of zeroes in fact the last column (the first column, respectively) in these matrices consists only of zeroes. In the same way,  $I^{\uparrow}_{n}$ ($I^{\downarrow}_{n}$) denotes an $n+1\times n$ matrix obtained from an $n\times n$ identity matrix by adding at the top (at the bottom) a row of zeroes.\par\bigskip

We recall that the solution of the Kronecker matrix problem allows to classify the indecomposable representations of the path algebra $kQ$ with $Q$ a quiver of the form (\ref{QKronecker}).\par\bigskip

Figure (26) shows the preprojective component of the Auslander-Reiten quiver of the 2-Kronecker quiver which has as vertices isomorphism classes of indecomposable representations of type III ($[i+1\hspace{0.2cm} i]$ is a notation for the dimension vector of a preprojective representation (equivalently, preprojective module), whereas $[m\hspace{0.2cm} m+1]$ is the dimension vector of a preinjective module). The preinjective component has isomorphism classes of indecomposable representations of type $\mathrm{III}^{*}$ as vertices.
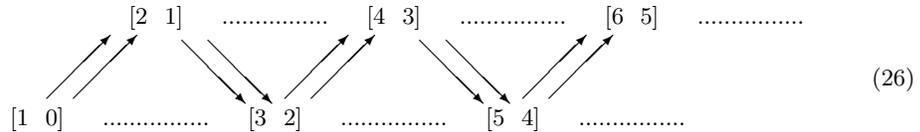
\begin{center}
\setlength{\unitlength}{1pt}
\begin{picture}(340,30)
	\put(0,-20){$[1 \hspace{0.22cm} 0]$}
	\put(90,-20){$[3 \hspace{0.22cm} 2]$}
	\put(180,-20){$[5 \hspace{0.22cm} 4]$}
	
	\put(45,18){$[2 \hspace{0.22cm} 1]$}
	\put(135,18){$[4 \hspace{0.22cm} 3]$}
	\put(225,18){$[6 \hspace{0.22cm} 5]$}
	
	\multiput(14,-10)(90,0){3}{\vector(1,1){24}}
	\multiput(24,-10)(90,0){3}{\vector(1,1){24}}
	
	\multiput(65,12)(90,0){2}{\vector(1,-1){24}}
	\multiput(75,12)(90,0){2}{\vector(1,-1){24}}
	
	\multiput(35,-20)(90,0){3}{$................$}
	\multiput(80,18)(90,0){3}{$................$}
	
	\put(326,-5){(26)}
\end{picture}
\end{center}
\par\bigskip\par\bigskip\par

\par\bigskip
\par\bigskip

Henceforth, we let $(n+1,n)$ ($(n, n+1)$) denote a representative of an isomorphism class of preprojective (preinjective) Kronecker modules obtained from a representation of type III (II) via elementary transformations. Actually, for the sake of simplicity, we will assume that such representatives have the form III (II).\par\bigskip

For $n\geq1$, let $P$ be an $(n+1)\times 2n$, $k$-matrix then $P$ can be partitioned into two $(n+1)\times n$ matrix blocks $A$ and $B$. In such a case we write $P=(P,A,B, n)$, where $A=(a_{i,j})=[C^{A}_{i_{1}},\dots, C^{A}_{i_{n}}]$, $B=(b_{i,j})=[C^{B}_{j_{1}},\dots, C^{B}_{j_{n}}]$, with $C^{A}_{i_{r}}$ ($C^{B}_{j_{s}}$) columns of $P$, if $I_{A}$ ($I_{B}$) is the set of indices $I_{A}=\{i_{r}\mid 1\leq r\leq n$\} ($I_{B}=\{j_{s}\mid 1\leq s\leq n\}$) then $I_{A}\cap I_{B}=\varnothing$, and $|I_A|=|I_B|=n$. In this case, each column of the matrix $P$ belongs either to the matrix $A$ or to the matrix $B$ and a word $W_{P}=l_{m_{1}}\dots l_{m_{n}}\dots l_{m_{2n}}$, \quad$l_{m_{h}}\in\{A, B\}$, $1\leq h\leq 2n$ is used to denote matrix $P$ by specifying the way that columns of $P$ have been assigned to the matrices $A$ and $B$.\par\bigskip

A row $r_{P}$ of $P$ has the form $(r_{A}, r_{B})$ with $r_{A}$ ($r_{B}$) being a row of the matrix block $A$ ($B$). We let $R_{A}$ ($R_{B}$) denote the set of rows of the matrix block $A$ ($B$), whereas $\mathscr{H}_{n}$ denotes the set of all matrices $P$ with the aforementioned properties.\par\bigskip

An helix associated to a matrix $P$ of type $\mathscr{H}_{n}$ is a connected directed graph $h$ whose construction goes as follows:

\begin{enumerate}[$(h_1)$]
\item  (\textit{Vertices}) Vertices of $h$ are entries of blocks $A$ and $B$. We let $h_{0}$ denote the set of vertices of $h$.
\item Fix two different rows $i_{P}=(i_{A}, i_{B})$ and $j_{P}=(j_{A}, j_{B})$ of $P$.
\item Choose sets $P_{A}$ and $P_{B}$ of \textit{pivoting entries} also called \textit{pivoting vertices}, $P_{A}\subset A$, $P_{B}\subset B$ such that $|P_{A}|=|P_{B}|=n$. Entries in $A\backslash P_{A}$ and $B\backslash P_{B}$ are said to be \textit{exterior entries} or \textit{exterior vertices}. In this case, if $x\in P_{A}$ ($x\in P_{B}$) then $x\notin i_{A}$ ($x\notin j_{B}$).\par\bigskip
$P_{A}$ and $P_{B}$ are sets of the form:
\addtocounter{equation}{1}
\begin{equation}\label{pivots}
\begin{split}
P_{A}&=\{a_{i_{1}, j_{1}}, a_{i_{2}, j_{2}},\dots, a_{i_{s}, j_{s}}\}, \quad j_{x}\neq j_{y}\hspace{0.1cm}\text{if and only if}\hspace{0.1cm}i_x\neq i_y,\\
P_{B}&=\{b_{t_{1}, h_{1}}, b_{t_{2}, h_{2}},\dots, b_{t_{s}, h_{s}}\}, \quad h_{x}\neq h_{y}\hspace{0.1cm}\text{if and only if}\hspace{0.1cm}t_x\neq t_y.
\end{split}
\end{equation}

Where, $a_{i_{r},j_{r}}\in R_{A}\backslash i_{A}$,\quad $b_{t_{m}, h_{m}}\in R_{B}\backslash j_{B}$, $1\leq r, m\leq s$. It is chosen just only one entry $a_{i_{r},j_{r}}$ ($b_{t_{m}, h_{m}}$) for each row in $R_{A}\backslash i_{A}$ ($R_{B}\backslash j_{B}$) and for each column $C_{A}$ ($C_{B}$) of $A$ ($B$).

\item (\textit{Arrows}) arrows in $h$ are defined in the following fashion;

\begin{enumerate}[$(a)$]
\item Arrows in $h$ are either horizontal or vertical. We let $h_{1}$ denote the set of arrows of $h$.
\item Horizontal arrows connect a vertex of the matrix block $A$ ($B$) with a vertex of the matrix block $B$ ($A$). Vertical arrows only connect vertices in the same matrix block. Starting and ending vertices of horizontal (vertical) arrows are entries of the same row (column) of $P$.
\item The starting vertex of a horizontal (vertical) arrow is an exterior (pivoting) vertex. The ending point of a horizontal (vertical) arrow is a pivoting (exterior) vertex. 
\item A pivoting vertex occurs as ending (starting) vertex just once. Thus, $h$ does not cross itself.

\item The first and last arrow of $h$ are horizontal and its starting vertex belongs to $i_{A}$. 

\item Each vertical arrow is preceded by a unique horizontal arrow, and unless the first arrow, any horizontal arrow is preceded by a vertical arrow. 

\item All the rows of $P$ are visited by $h$, and no row or column of $P$ is visited by arrows of $h$ more than once. 

\item There are not horizontal arrows connecting exterior vertices of $j_{A}$ with vertices of $j_{B}$.

\end{enumerate}

\end{enumerate}

\addtocounter{Nota}{-3}
\begin{Nota}\label{system}
We let $(i_{P}, j_{P}, P_{A}, P_{B})$ denote the set of all helices which can be built by fixing these data associated to a matrix $P$ of type $\mathscr{H}_{n}$, $h^{P}_{n}=|(i_{P}, j_{P}, P_{A}, P_{B})|$ denotes the corresponding cardinality.  See figures ((54), (56), (58), (60)) in section \ref{Appendix} where it is presented a set $(4_{P}, 2_{P}, P_{A}, P_{B})$ defined by the word $BAABAB$.
\end{Nota}

Matrix presentations of preprojective Kronecker modules $p$ of type III are of type $\mathscr{H}_{n}$, $n\geq1$ (In this case, $W_{p}=AA\dots ABB\dots B$). In \cite{Canadas3}, the first author et al studied sets of helices $(1_{p}, (n+1)_{p}, p_{A}=\{a_{i+1,i}\mid 1\leq i\leq n\}, p_{B}=\{b_{i, i}\mid 1\leq i\leq n\})$ associated to this kind of matrices.\par\bigskip

\addtocounter{prop}{2}
\begin{prop}\label{word partition}
\textit{If $\mathscr{W}_{P}$ is the set of matrix words associated to a matrix $P$ of type $\mathscr{H}_{n}$ then $|\mathscr{W}_{P}|$ equals $\underset{m=0}{\overset{n^{2}}{\sum}}P(n, n, m)=(n+1)C_{n}$, where $P(n, n, m)$ denotes the number of partitions of $m$ into $n$ parts, each $\leq n$, $P(n, n, 0)=1$, and $C_{n}$ denotes the $n$th Catalan number.}
\end{prop}

\textbf{Proof.} Each matrix word $W_{P}$ of the form $W_{P}=l_{m_{1}}\dots l_{m_{n}}\dots l_{m_{2n}}$, \quad$l_{m_{h}}\in\{A, B\}$, $1\leq h\leq 2n$, gives rise to an integer partition $\lambda=(\lambda_{1},\lambda_{2},\dots, \lambda_{t})$, $\lambda_{i}, t\leq n$ of a nonnegative integer number $m\leq n^{2}$ by defining $\lambda_{1}$ as the number of $A$'s after the first occurrence of the letter $B$, $\lambda_{2}$ is the number of $A'$s after the second occurrence of the letter $B$ and so on. Since there are $n$ letters $A'$s and $n$ letters $B'$s in $W_{P}$ then the number of words associated to $P$ is $\binom{2n}{n}$. The result holds.\hspace{0.5cm}$\square$ 
\par\bigskip

In Theorem \ref{Kronecker0} we prove that the number of helices $h^{p}_{n}$ associated to a preprojective Kronecker module, $p=(n+1,n)$ is  $h^{p}_{n}=n!\lceil \frac{n}{2}\rceil$. 
\par\bigskip
 If we associate to a set of $n$ equidistant points on a circle the rows of a representation $p=(n+1,n)$ then the number of helices containing the fixed arrow $a_{1,1}\rightarrow b_{1,1}$ equals the number $a(n)$ of ways of connecting $n+1$ equally spaced points on a circle with a path of $n$ line segments ignoring reflections. In this case, vertical edges in a helix are in bijective correspondence with the edges of the path in the circle (Figure (29) shows examples of helices and this kind of paths). Thus

 \par\bigskip
 
\begin{equation}\label{a(n)}
\begin{split}
a(n)&=\frac{h^{p}_{n}}{n}, \quad n\geq 1.
\end{split}
\end{equation}

\begin{center}
\begin{figure}[H]
\setlength{\unitlength}{1pt}
\begin{picture}(140,100)
\put(330,38){$(29)$}
\multiput(120,50)(30,0){3}{\line(0,1){30}}
\multiput(120,50)(0,30){2}{\line(1,0){60}}
\put(75,65){$(2,1)=$}

\put(130,70){$0$}
\put(135,72){$\vector(1,0){30}$}
\put(167,70){$\vector(0,-1){9}$}
\put(164,58){$\vector(-1,0){30}$}
\put(130,55){$1$}
\put(165,70){$1$}
\put(165,55){$0$}
\put(220,65){\circle{30}}
\put(220,20){\circle{30}}
\put(211,32){$\bullet$}
\put(213,33){\vector(0,-1){26}}
\put(211,4){$\bullet$}
\put(213,6){\vector(1,1){21}}
\put(233,25){$\bullet$}
\put(237,25){$2$}
\put(211,37){$1$}
\put(211,-3){$3$}
\put(204,65){\vector(1,0){32}}
\put(198,63){1}
\put(202,63){$\bullet$}
\put(238,63){2}
\put(234,63){$\bullet$}

\multiput(120,0)(30,0){3}{\line(0,1){40}}
\multiput(120,0)(0,40){2}{\line(1,0){60}}          
\put(75,15){$(3,2)=$}

\put(125,30){$0$}
\put(125,17.5){$1$}
\put(125,5){$0$}

\put(140,17.5){$0$}
\put(140,30){$0$}
\put(140,5){$1$}
\put(155,17.5){$0$}
\put(155,30){$1$}
\put(155,5){$0$}

\put(170,17.5){$1$}
\put(170,30){$0$}
\put(170,5){$0$}

\put(144,34){$\vector(1,0){12}$}
\put(157,35){$\vector(0,-1){23}$}
\put(155,9){$\vector(-1,0){25}$}
\put(127,12){$\vector(0,1){7}$}
\put(128,23){$\vector(1,0){42}$}

\end{picture}
\end{figure}

\end{center}

Sequence $a(n)$ is recorded as A052558 in the OEIS.

\subsection{The Four Subspace Problem (FSP)}\label{FSP01}
The four subspace problem is another example of a matrix problem, in this case if $k$ is an arbitrary field then a \textit{quadruple} of finite-dimensional $k$ vector spaces is a system of the form 
\par\smallskip
\begin{centering}
$U=(U_{0}, U_{1}, U_{2}, U_{3}, U_{4})$ \par\smallskip
\end{centering}

where $U_{0}$ is a finite- dimensional $k$ vector space and $U_{1},\dots, U_{4}$ is an ordered collection of four subspaces of $U_{0}$. Two quadruples are said to be isomorphic if there exists a $k$-space isomorphism $\varphi: U_{0}\rightarrow V_{0}$ such that $\varphi(U_{i})=V_{i}$ for all $i$. And a quadruple $U$ is \textit{decomposable} ($U=U'\oplus U''$) if some non-trivial direct sum decomposition $U_{0}=U'_{0}\oplus U''_{0}$ satisfies the identity $U_{i}=(U_{i}\cap U'_{0})\oplus (U_{i}\cap U''_{0})$ for each $i\in\{1,\dots, 4\}$.  \par\bigskip

The four subspace problem consists of classifying all indecomposable quadruples up to isomorphism, it is equivalent to determine indecomposable representations of four incomparable points or a tetrad. Actually, it is essentially equivalent to determine the indecomposable representations of the four subspace quiver (see Figure (14)). \par\bigskip

Given two matrix  representations $M=(M_{x_{i}},\mid 1\leq i\leq 4)$ and $M'=(M'_{x_{i}},\mid 1\leq i\leq 4)$ of the tetrad then $M$ and $M'$ are said to be equivalent or isomorphic, if one can be turned into the other by means of the following \textit{admissible transformations}:

\begin{enumerate}
\item $k$-elementary transformations of rows of the whole matrix.
\item $k$-elementary transformations of columns of matrices $M_{x_{i}}$.
\end{enumerate}

FSP was solved by Gelfand and Ponomarev in 1970 for $k$ algebraically closed and by Nazarova (1967-1973) for the arbitrary case. An advance to this problem was given by Brenner who described the indecomposable quadruples with non-zero defect $\partial(U)=\underset{i=1}{\overset{4}{\sum}}\mathrm{dim}\hspace{0.1cm}U_{i}-2\mathrm{dim}\hspace{0.1cm}U_{0}$ (called non-regular) in particular she extended the results of Gelfand and Ponomarev to the case of a skew field $k$. Afterwards, in 2004 Zavadskij and Medina gave an elementary solution of this problem \cites{Brenner1, Brenner2, Gelfand, Nazarova2, Zavadskij1}.
\par\bigskip

For $n\geq 1$, we consider $(2n+2)\times (4n+3)$-matrices of type $\mathscr{C}_{n}$, whose rows and columns are partitioned by four adjacent matrix blocks (from the left to the right), $U_{1}, U_{2}, U_{3}$ and $U_{4}$ denoted $(A_{1}, A'_{1});(A_{2},A'_{2});(A_{3},A'_{3}),(A_{4},A'_{4})$, respectively. $A_{i}$ and $A'_{i}$ are matrices vertically adjacent (blocks of type $A_{i}$ as well as those of type $A'_{i}$ are horizontally adjacent) of the same size, three of the four blocks $U_{i}$ consists of $(n+1)\times (n+1)$-matrices, and the remaining block consists of two $(n+1)\times n$-matrices.  \par\bigskip

If $U$ is a matrix of type $\mathscr{C}_{n}$ then $i$th row ($j$th column) $i^{R}_{U}$ ($j^{C}_{U}$)  of $U$ is given by the union

\[i^{R}_{U}=
\begin{cases}
\underset{m=1}{\overset{4}{\cup}}i^{R}_{A_{m}},  & \hspace{0.2cm}\text{if}\hspace{0.2cm}i^{R}_{U}\subset\underset{m=1}{\overset{4}{\cup}}A_{m}, \\
\underset{m=1}{\overset{4}{\cup}}i^{R}_{A'_{m}},  & \hspace{0.2cm}\text{if}\hspace{0.2cm}i^{R}_{U}\subset\underset{m=1}{\overset{4}{\cup}}A'_{m},
\end{cases}\]

where $i^{R}_{A_{m}}$ and $i^{R}_{A'_{m}}$ are corresponding rows in the matrices $A_{m}$ and $A'_{m}$, $m=1,\dots, 4$.

The $j$th column $j^{C}_{U}$  of $U$ is given by the union

\begin{equation}
j^{C}_{U}=j^{C}_{A_{s}}\cup j^{C}_{A'_{s}}\notag\quad\text{if}\hspace{0.2cm}j^{C}_{U}\subset A_{s}\cup A'_{s},
\end{equation}

The following is a typical shape of a matrix $U$ of type $\mathscr{C}_{n}$:

\begin{center}

$U=\begin{array}{|c|c|c|c|}\hline A_{1} & A_{2} & A_{3} & A_{4}  \\\hline A'_{1} & A'_{2} & A'_{3} & A'_{4}\\\hline\end{array}$

\end{center}

Matrices of type $\mathscr{C}_{n}$ have associated cycles which are connected oriented graph whose construction goes as follows:

\begin{enumerate}[$(cl_1)$]
\item (\textit{Vertices}) By definition any matrix $U$ of type $\mathscr{C}_{n}$ is defined by a word $w$ whose letters are the symbols $U_{1}, U_{2}, U_{3}$ and $U_{4}$, i.e., $w=U_{\sigma(1)}U_{\sigma(2)}U_{\sigma(3)}U_{\sigma(4)}$ where $\sigma$ is a permutation of $4$ elements. For the sake of clarity, later on, we assume that the matrix blocks of a matrix $U$ of type $\mathscr{C}_{n}$ are organized according to words of the form $w=U_{\sigma(1)}U_{\sigma(2)}U_{\sigma(3)}U_{\sigma(4)}$ with $\sigma(i)=i$, $1\leq i\leq 4$. Blocks $U_{1}$ and $U_{4}$ are said to be \textit{external blocks}, whereas blocks $U_{2}$ and $U_{3}$ are \textit{internal blocks}.

\item Fix three rows $i^{R}_{A_{i_{1}}}$, $j^{R}_{A_{i_{2}}}$ and $j'^{R}_{A'_{i_{2}}}$. In this case, $U_{i_{1}}$ is a $(2n+2)\times (n+1)$-matrix block, whereas $U_{i_{2}}$ is a $(2n+2)\times n$-matrix block. Thus, according to our choice, $i_{1}\in\{1,2,3\}$ and $i_{2}=4$.
\item Entries of matrix $U$ of type $\mathscr{C}_{n}$ are either pivoting or exterior vertices (or entries). They are defined as in the case for matrices of type $\mathscr{H}_{n}$, taking into account, that in this case there is not a row $i_{A}$ to apply the restriction for the possible choices of pivoting entries, and instead of $j_{B}$, we consider rows $j^{R}_{A_{4}}$ and $j'^{R}_{A'_{4}}$, which do not contain pivoting entries. If $P_{A^{r}_{m}}$ denotes the set of pivoting entries of the matrix $A^{r}_{m}$, $r=0,1$, $m=1,\dots,4$ with $A^{0}_{i}=A_{i}$, $A^{1}_{i}=A'_{i}$. Then $|P_{A^r_{i}}|=|P_{A^r_{j}}|=n+1$, if $i,j\in \{1,2,3\}$ and  $|P_{A_{4}}|=|P_{A'_{4}}|=n$.

\item Without loss of generality, we assume that the starting and ending vertex are the same and belong to $i^{R}_{A_{1}}\subset U_{1}$, it is a pivoting vertex. Furthermore,  $j^{R}_{A_{i_{2}}}\subset A_{4}$ and $j'^{R}_{A'_{i_{2}}}\subset A'_{4}$. 
\item According to the assumption introduced in $(cl_{1})$, the sequence of vertices connected by a cycle $\mathcal{C}$ is of the form,

 \begin{equation}
 \begin{split}
\mathcal{C}_{0}&=\{a_{i_{1}, j_{1}}, a'_{i'_{1}, j_{1}}, b'_{i'_{1},j'_{1}}, b_{i_{2}, j'_{1}}, c_{i_{2}, j_{2}}, c'_{i'_{3}, j_{2}}, d'_{i'_{3}, j'_{2}} \}\cup \{d_{i_{4}, j'_{2}}, c_{i_{4}, j_{3}},c'_{i'_{4}, j_{3}}, b'_{i'_{4},j'_{4}}\}\\\notag
&\cup\{b_{i_{1}, j'_{4}}, a_{i_{1},j_{1}}\}.
\end{split}
\end{equation}

Where entries $(a, a')$, $(b,b')$, $(c,c')$ and $(d,d')$ correspond respectively to the blocks $(A_{1}, A'_{1})$, $(A_{2}, A'_{2})$, $(A_{3}, A'_{3})$ and $(A_{4}, A'_{4})$.
\item Any cycle contains the fixed vertices $a_{i_{1}, j_{1}}$, $a'_{i'_{1}, j_{1}}$, $b'_{i'_{1}, j'_{1}}$, $b_{i_{2}, j_{1}}$ and $c_{i_{2}, j_{2}}$, which is a pivoting entry. Entries of the form $d'_{j'^{R}_{A_{i_{2}}}, j_{2}}\notin\mathcal{C}_{0}$. No horizontal arrow has an entry $d_{j^{R}_{A_{i_{2}}}, j}\in A_{4}$ as its ending vertex.

\item Vertices $a_{i_{1}, j_{1}}\in A_{1}$,  $b'_{i'_{1}, j'_{1}}\in A'_{2}$, $c_{i_{2}, j_{2}}\in A_{3}$, $d'_{i'_{3}, j'_{2}}\in A'_{4}$,  $c_{i_{4}, j_{3}}\in A_{3}$ and $b'_{i'_{4}, j'_{4}}\in A'_{2}$  are pivoting entries.

\item (\textit{Arrows}) arrows in $\mathcal{C}$ connect alternatively pivoting entries with exterior entries. They are defined in the following fashion;

\begin{enumerate}[$(a)$]
\item Arrows in $\mathcal{C}$ are either horizontal ($\rightarrow$, $\leftarrow$) or vertical ($\uparrow$, $\downarrow$). We let $\mathcal{C}_{1}$ denote the set of arrows of $\mathcal{C}$. We write, $X\rightarrow Y$, $X\leftarrow Y$, $X\uparrow Y$ and $X\downarrow Y$, the different ways of connecting matrices of the different blocks $U_{i}$. According to the sequence of vertices, the first vertical arrow of the cycle $\mathcal{C}$ is of the form $A_{1}\downarrow A'_{1}$.
\item Since we are assuming cycles associated to a word of the form $w=U_{1}U_{2}U_{3}U_{4}$ then horizontal arrows connect adjacent matrices $A_{i}$ ($A'_{i}$) with $A_{i+1}$ ($A'_{i+1}$) (conversely, $A_{i+1}$ ($A'_{i+1}$) with $A_{i}$ ($A'_{i}$)). The first horizontal arrow in our case is of the form $A'_{1}\rightarrow A'_{2}$. Vertical arrows connect matrices in the form $A_{i}\uparrow A'_{i}$ or $A'_{i}\downarrow A_{i}$. External blocks are connected by a unique vertical arrow (either $A_{1}\downarrow A'_{1}$ iff $A_{4}\uparrow A'_{4}$ or $A'_{1}\uparrow A_{1}$ iff $A_{4}\downarrow A'_{4}$). Two arrows, $U_{2}\stackrel{\rightarrow}{\leftarrow}U_{3}$, $U_{3}\stackrel{\rightarrow}{\leftarrow}U_{4}$ connect internal matrix blocks ($A_{2}\rightarrow A_{3}$ iff $A'_{2}\leftarrow A'_{3}$, $A_{3}\rightarrow A_{4}$ iff $A'_{3}\leftarrow A'_{4}$ ). Two vertical arrows connect internal matrix blocks $A_{2}\downarrow\downarrow A'_{2}$ or $A_{2}\uparrow\uparrow A'_{2}$, same conditions satisfy matrices $A_{3}$ and $A'_{3}$.
\item Each horizontal arrow is preceded by a unique vertical arrow, and unless the first arrow, any vertical arrow is preceded by an horizontal arrow. 

\item No row or column of a matrix $A^{r}_{m}$ in a matrix block $U_{m}$, $r=0,1$, $m=1,\dots, 4$ is visited by arrows of $\mathcal{C}$ more than once.

\end{enumerate}

\end{enumerate}

As for matrices of type $\mathscr{H}_{n}$, any matrix $U$ of type $\mathscr{C}_{n}$ can be described in the form $(U, U_{i},n\mid 1\leq i\leq 4)$. We let $(i^{R}_{A_{1}}, j^{R}_{A_{4}},j'^{R}_{A'_{4}}, P_{A^{r}_{i}}\mid 1\leq i\leq 4, r\in\{0,1\})$ denote the  system of all cycles which can be constructed by using these data associated to a matrix $U$ of type $\mathscr{C}_{n}$ and $\mathcal{C}_{U}=|(i^{R}_{A_{1}}, j^{R}_{A_{4}},j'^{R}_{A'_{4}}, P_{A^{r}_{i}}\mid 1\leq i\leq 4, r\in\{0,1\})|$ are notations for its corresponding cardinal.
\par\bigskip

In this work, we compute the number of cycles associated to preprojective representations of type IV  of the tetrad and establish a connection between these preprojective representations with indecomposable projective modules over some Brauer configuration algebras via such cycles. The following is the canonical matrix presentation of the mentioned preprojective representations, where $I_{n}$ is an $n\times n$ identity matrix. $n$ is said to be the \textit{order} of the representation. \par\smallskip

\begin{centering}

\begin{picture}(50,2)
\multiput(-10,-40)(30,0){5}{\line(0,1){36}}
\multiput(-10,-40)(0,18){3}{\line(1,0){120}}
\multiput(30,-35)(30,0){2}{$\mathrm{I}_{n+1}$}
\multiput(90,-35)(30,0){1}{$\mathrm{I}^{\downarrow}_{n}$}
\multiput(5,-35)(30,0){1}{$0$}
\multiput(0,-18)(30,0){1}{$\mathrm{I}_{n+1}$}
\multiput(30,-18)(30,0){1}{$0$}
\multiput(60,-18)(30,0){1}{$\mathrm{I}_{n+1}$}
\multiput(90,-18)(30,0){1}{$\mathrm{I}^{\uparrow}_{n}$}

\multiput(180,-35)(30,0){1}{(30)}

\end{picture}       
           
\par\bigskip
\end{centering}                    
\vspace{2cm}

The following is an example of a cycle of type \par\bigskip
\begin{centering}
$(1^{R}_{A_{1}}, 1^{R}_{A_{4}}, 4^{R}_{A'_{4}}, P_{A^{r}_{i}}=\underset{\underset{r=0,1}{1\leq i\leq 3}}{\cup}d^{r}_{i}\cup\underline{d_{4}}\cup d'_{4}\mid r\in\{0,1\})$ \par\bigskip
\end{centering}
associated to a preprojective representation of order $n=3$ of the tetrad, black arrows connect fixed vertices of the corresponding cycles, whereas $d^{r}_{i}$ ($\underline{d_{4}}$) denote the set of diagonal (subdiagonal) entries used to define the corresponding pivoting vertices.

\hspace{5mm}
\setlength{\unitlength}{0.97pt}
\begin{picture}(90,90)

 \multiput(0,0)(90,0){4}{\line(0,1){78}}
 \multiput(333,0)(90,0){1}{\line(0,1){78}}        
  \multiput(0,0)(0,39){3}{\line(1,0){333}}
   \multiput(5,72)(30,0){1}{\tiny 1}
    \multiput(30,72)(30,0){1}{\tiny 0}
  \multiput(55,72)(30,0){1}{\tiny 0}
    \multiput(80,72)(30,0){1}{\tiny 0}
      \multiput(5,62)(30,0){1}{\tiny 0}
    \multiput(30,62)(30,0){1}{\tiny 1}
  \multiput(55,62)(30,0){1}{\tiny 0}
    \multiput(80,62)(30,0){1}{\tiny 0}
 \multiput(5,52)(30,0){1}{\tiny 0}
    \multiput(30,52)(30,0){1}{\tiny 0}
  \multiput(55,52)(30,0){1}{\tiny 1}
    \multiput(80,52)(30,0){1}{\tiny 0}
       \multiput(5,42)(30,0){1}{\tiny 0}
    \multiput(30,42)(30,0){1}{\tiny 0}
  \multiput(55,42)(30,0){1}{\tiny 0}
    \multiput(80,42)(30,0){1}{\tiny 1}
 
   \multiput(95,72)(30,0){1}{\tiny 0}
    \multiput(120,72)(30,0){1}{\tiny 0}
  \multiput(145,72)(30,0){1}{\tiny 0}
    \multiput(170,72)(30,0){1}{\tiny 0}
      \multiput(95,62)(30,0){1}{\tiny 0}
    \multiput(120,62)(30,0){1}{\tiny 0}
  \multiput(145,62)(30,0){1}{\tiny 0}
    \multiput(170,62)(30,0){1}{\tiny 0}
 \multiput(95,52)(30,0){1}{\tiny 0}
    \multiput(120,52)(30,0){1}{\tiny 0}
  \multiput(145,52)(30,0){1}{\tiny 0}
    \multiput(170,52)(30,0){1}{\tiny 0}
       \multiput(95,42)(30,0){1}{\tiny 0}
    \multiput(120,42)(30,0){1}{\tiny 0}
  \multiput(145,42)(30,0){1}{\tiny 0}
    \multiput(170,42)(30,0){1}{\tiny 0}
   
    \multiput(185,72)(30,0){1}{\tiny 1}
    \multiput(210,72)(30,0){1}{\tiny 0}
  \multiput(235,72)(30,0){1}{\tiny 0}
    \multiput(260,72)(30,0){1}{\tiny 0}
      \multiput(185,62)(30,0){1}{\tiny 0}
    \multiput(210,62)(30,0){1}{\tiny 1}
  \multiput(235,62)(30,0){1}{\tiny 0}
    \multiput(260,62)(30,0){1}{\tiny 0}
 \multiput(185,52)(30,0){1}{\tiny 0}
    \multiput(210,52)(30,0){1}{\tiny 0}
  \multiput(235,52)(30,0){1}{\tiny 1}
    \multiput(260,52)(30,0){1}{\tiny 0}
       \multiput(185,42)(30,0){1}{\tiny 0}
    \multiput(210,42)(30,0){1}{\tiny 0}
  \multiput(235,42)(30,0){1}{\tiny 0}
    \multiput(260,42)(30,0){1}{\tiny 1}
     \multiput(275,72)(30,0){1}{\tiny 0}
    \multiput(300,72)(30,0){1}{\tiny 0}
  \multiput(325,72)(30,0){1}{\tiny 0}
       \multiput(275,62)(30,0){1}{\tiny 1}
    \multiput(300,62)(30,0){1}{\tiny 0}
  \multiput(325,62)(30,0){1}{\tiny 0}
 
 \multiput(275,52)(30,0){1}{\tiny 0}
    \multiput(300,52)(30,0){1}{\tiny 1}
  \multiput(325,52)(30,0){1}{\tiny 0}
    
       \multiput(275,42)(30,0){1}{\tiny 0}
    \multiput(300,42)(30,0){1}{\tiny 0}
  \multiput(325,42)(30,0){1}{\tiny 1}

    \multiput(5,2)(30,0){1}{\tiny 0}
    \multiput(30,2)(30,0){1}{\tiny 0}
  \multiput(55,2)(30,0){1}{\tiny 0}
    \multiput(80,2)(30,0){1}{\tiny 0}
      \multiput(5,12)(30,0){1}{\tiny 0}
    \multiput(30,12)(30,0){1}{\tiny 0}
  \multiput(55,12)(30,0){1}{\tiny 0}
    \multiput(80,12)(30,0){1}{\tiny 0}
 \multiput(5,22)(30,0){1}{\tiny 0}
    \multiput(30,22)(30,0){1}{\tiny 0}
  \multiput(55,22)(30,0){1}{\tiny 0}
    \multiput(80,22)(30,0){1}{\tiny 0}
       \multiput(5,32)(30,0){1}{\tiny 0}
    \multiput(30,32)(30,0){1}{\tiny 0}
  \multiput(55,32)(30,0){1}{\tiny 0}
    \multiput(80,32)(30,0){1}{\tiny 0}
   
       \multiput(95,2)(30,0){1}{\tiny 0}
    \multiput(120,2)(30,0){1}{\tiny 0}
  \multiput(145,2)(30,0){1}{\tiny 0}
    \multiput(170,2)(30,0){1}{\tiny 1}
      \multiput(95,12)(30,0){1}{\tiny 0}
    \multiput(120,12)(30,0){1}{\tiny 0}
  \multiput(145,12)(30,0){1}{\tiny 1}
    \multiput(170,12)(30,0){1}{\tiny 0}
 \multiput(95,22)(30,0){1}{\tiny 0}
    \multiput(120,22)(30,0){1}{\tiny 1}
  \multiput(145,22)(30,0){1}{\tiny 0}
    \multiput(170,22)(30,0){1}{\tiny 0}
       \multiput(95,32)(30,0){1}{\tiny 1}
    \multiput(120,32)(30,0){1}{\tiny 0}
  \multiput(145,32)(30,0){1}{\tiny 0}
    \multiput(170,32)(30,0){1}{\tiny 0}

      \multiput(185,2)(30,0){1}{\tiny 0}
    \multiput(210,2)(30,0){1}{\tiny 0}
  \multiput(235,2)(30,0){1}{\tiny 0}
    \multiput(260,2)(30,0){1}{\tiny 1}
      \multiput(185,12)(30,0){1}{\tiny 0}
    \multiput(210,12)(30,0){1}{\tiny 0}
  \multiput(235,12)(30,0){1}{\tiny 1}
    \multiput(260,12)(30,0){1}{\tiny 0}
 \multiput(185,22)(30,0){1}{\tiny 0}
    \multiput(210,22)(30,0){1}{\tiny 1}
  \multiput(235,22)(30,0){1}{\tiny 0}
    \multiput(260,22)(30,0){1}{\tiny 0}
       \multiput(185,32)(30,0){1}{\tiny 1}
    \multiput(210,32)(30,0){1}{\tiny 0}
  \multiput(235,32)(30,0){1}{\tiny 0}
    \multiput(260,32)(30,0){1}{\tiny 0}

     \multiput(275,2)(30,0){1}{\tiny 0}
    \multiput(300,2)(30,0){1}{\tiny 0}
  \multiput(325,2)(30,0){1}{\tiny 0}
  
      \multiput(275,12)(30,0){1}{\tiny 0}
    \multiput(300,12)(30,0){1}{\tiny 0}
  \multiput(325,12)(30,0){1}{\tiny 1}
   \multiput(275,22)(30,0){1}{\tiny 0}
    \multiput(300,22)(30,0){1}{\tiny 1}
  \multiput(325,22)(30,0){1}{\tiny 0}
    
       \multiput(275,32)(30,0){1}{\tiny 1}
    \multiput(300,32)(30,0){1}{\tiny 0}
  \multiput(325,32)(30,0){1}{\tiny 0}

  \multiput(7,72)(0,90){1}{\vector(0,-1){50}}    \multiput(7,22)(0,90){1}{\vector(1,0){115}} 
  \multiput(122,22)(0,90){1}{\vector(0,1){20}}  
 
     \multiput(122,42)(0,90){1}{\vector(1,0){138}}  
   
 \color{red}  
   
    \multiput(262,42)(0,90){1}{\vector(0,-1){30}} 
     \multiput(262,12)(0,90){1}{\vector(1,0){65}} 
      \multiput(326,12)(0,90){1}{\vector(0,1){61}} 
       \multiput(326,72)(90,0){1}{\vector(-1,0){140}} 
       \multiput(186,72)(90,0){1}{\vector(0,-1){60}} 
        \multiput(186,12)(90,0){1}{\vector(-1,0){40}} 
         \multiput(146,12)(90,0){1}{\vector(0,1){60}} 
          \multiput(146,72)(90,0){1}{\vector(-1,0){140}}

\end{picture}       

\par\bigskip

\subsection{Binomial trees and integer partitions}

Binomial trees appear in many fields of the mathematics, they are binary trees with the shape~\cite{Knuth}:

\addtocounter{equation}{2}
\begin{equation}
\begindc{\commdiag}[250]
\obj(0,0)[p0]{$\bullet$}
\obj(-2,-2)[p1]{$T_0$}
\obj(2,-2)[p2]{$T_1$}
\obj(4,-2)[p3]{$\cdots$}
\obj(6,-2)[p4]{$T_{n-1}$}
\mor{p0}{p1}{$0$}[0,0]
\mor{p0}{p2}{$1$}
\mor{p0}{p4}{$n-1$}
\enddc
\end{equation}

\par\bigskip

As an example $T_{4}$ has the following form:

\begin{equation}
\begindc{\commdiag}[250]
\obj(0,0)[p0]{$\bullet$}
\obj(-1,-2)[p1]{$\bullet$}
\obj(1,-2)[p2]{$\bullet$}
\obj(3,-2)[p3]{$\bullet$}
\obj(8,-2)[p4]{$\bullet$}
\obj(1,-4)[p5]{$\bullet$}
\obj(2,-4)[p6]{$\bullet$}
\obj(4,-4)[p7]{$\bullet$}
\obj(6,-4)[p8]{$\bullet$}
\obj(7,-4)[p9]{$\bullet$}
\obj(9,-4)[p10]{$\bullet$}
\obj(7,-6)[p11]{$\bullet$}
\obj(8,-6)[p12]{$\bullet$}
\obj(10,-6)[p13]{$\bullet$}
\obj(10,-8)[p14]{$\bullet$}
\mor{p0}{p1}{$0$}[0,0]
\mor{p0}{p2}{$1$}
\mor{p0}{p3}{$2$}
\mor{p0}{p4}{$3$}
\mor{p2}{p5}{$0$}
\mor{p3}{p6}{$0$}[0,0]
\mor{p3}{p7}{$1$}
\mor{p4}{p8}{$0$}[0,0]
\mor{p4}{p9}{$1$}
\mor{p4}{p10}{$2$}
\mor{p9}{p11}{$0$}
\mor{p10}{p12}{$0$}[0,0]
\mor{p10}{p13}{$1$}
\mor{p13}{p14}{$0$}
\enddc
\end{equation}

Note that, at each level $T_{4}$ gives integer partitions of numbers 1, 2 and 3 (without taking into account 0 as a part). Often, these types of trees are said to be \textit{partition trees} which can be used to store partitions of a given positive integer $n$ or of all positive integers $\leq n$~\cite{Lin}. In \cite{Luschny}, Luschny describes partition trees for different integer numbers and uses them to define some orders on the set of integer partitions.

\section{Categorification of the Sequences A052558 and A052591}
Results in this section can be interpreted as categorifications (see Remark \ref{A052558}) of the sequences $(n-1)!\lceil \frac{n}{2}\rceil$ and $n!\lceil \frac{n}{2}\rceil$  (A052558 and A052591 in the OEIS, respectively) via Kronecker modules and Brauer configuration algebras. Theorems \ref{invariant11} and \ref{Kronecker0} and Corollaries \ref{invariant2} and \ref{invariant3} prove that  the number of helices associated to preprojective Kronecker modules is invariant with respect to admissible transformations. 

\par\bigskip
The following results regard the number of helices associated to matrices of type $\mathscr{H}_{n}$ and in particular to preprojective Kronecker modules (see Figure (26)), $k$ is an algebraically closed field.
\addtocounter{teor}{2}

\begin{teor}\label{invariant11}
\textit{Let $(P,A,B,n)$, $(P',A',B',n)$, $H_{P}$ and $H_{P'}$ be two matrices of type $\mathscr{H}_{n}$ with corresponding sets of helices $H_{P}$ and $H_{P'}$ defined by systems of the form $(i_{P},j_{P}, P_{A}, P_{B})$ and $(f_{P'}, g_{P'}, P'_{A'}, P'_{B'})$, respectively. Then $|H_{P}|=h^{P}_{n}=|H_{P'}|=h^{P'}_{n}$}.
\end{teor}

 \textbf{Proof.} Firstly, we suppose without loss of generality that, $i_{P}\neq f_{P'}$ and $j_{P}\neq g_{P'}$. Then, we note that each helix $h\in (i_{P},j_{P}, P_{A}, P_{B})$ gives rise to a unique helix $h'\in (i_{P}, g_{P}, P_{A''}, P_{B''})$, where $P_{A''}$ and $P_{B''}$ are suitable sets of pivoting entries in $P$. The process consists of copying helix $h$, in such a way that each occurrence of entries of $j_{P}$ is substituted by a corresponding occurrence of $g_{P}$ (taking into account the new sets of pivoting vertices, $P_{A''}$ and $P_{B''}$), conversely, each occurrence of $g_{P}$ is substituted by a corresponding occurrence of $j_{P}$, keeping without changes the remaining rows visited by the helix $h$. For example, 
  if a vertical arrow $v\in h_{1}$ connects entries $p_{i, j}$ (starting vertex) and $p_{i', j}$ (ending vertex) in $P$ then the corresponding vertical arrow $v'\in h'_{1}$ connects entries of the rows $i$ and $i'$ if $i\in\{ j_{P}, g_{P}\}$ and $i'\notin \{j_{P}, g_{P}\}$ or if $i$ and $i'$ are such that $i, i'\notin\{j_{P}, g_{P}\}$, $v'$ connects rows $i$ and $j$ ($g$) if $i'=g$ ($i'=j$). We let $\sigma$ denote the bijection, $\sigma:(i_{P}, j_{P}, P_{A}, P_{B})\rightarrow (i_{P}, g_{P}, P_{A''}, P_{B''})$ defined by these substitutions. Thus, if a bijection $\delta:(i_{P}, g_{P}, P_{A''}, P_{B''})\rightarrow (f_{P},g_{P}, P_{A'},P_{B'})$ is defined as $\sigma$ where $P_{A'}$ and $P_{B'}$ are sets of pivoting entries of $P$ given by $P'_{A'}$ and $P'_{B'}$ respectively,  then the maps composition $\delta\sigma$ is also a bijection from $(i_{P},j_{P}, P_{A}, P_{B})$ to $(f_{P}, g_{P}, P_{A'}, P_{B'})$. Any helix $h'\in (f_{P},g_{P}, P_{A'},P_{B'})$ corresponds uniquely to an helix $h''\in(f_{P'},g_{P'}, P'_{A'},P'_{B'})$ via the identification $\tau: P\rightarrow P'$ such that $\tau(p_{i,j})=p'_{i,j}$, in this case $p_{i,j}\in P_{A'}$ ($p_{i,j}\in P_{B'}$) if and only if $p'_{i, j}\in P'_{A'}$ ($p'_{i, j}\in P'_{B'}$). In general, a copy $h''\in (f_{P'},g_{P'}, P'_{A'},P'_{B'})$ of an helix $h'\in (f_{P},g_{P}, P_{A},P_{B})$ 
  can be built taking into account that an initial exterior vertex $e_{f, j}\in h'$ has a vertex $e'_{f, j'}\in f_{A'}$ as its corresponding initial exterior copy and $h''$ visits the same rows in the same order as those visited previously by $h'$. We are done. \hspace{0.5cm}$\square$\par\bigskip
  
 Examples of copies of elements of the set  
 
 \par\bigskip
 \begin{centering}
 $(4_{P}, 2_{P}, P_{A}=\{p_{3,2},p_{1,3},p_{2,5}\}, P_{B}=\{p_{4,1},p_{1,4}, p_{3,6}\})$\par\bigskip
 \end{centering}
 associated to a matrix $P$ of type $\mathscr{H}_{3}$ and defined by the word $W_{P}=BAABAB$ are given in section \ref{Appendix} (see figures ((55), (57), (59), (61)). In such a case, the corresponding copies belong to the set \par\bigskip
 \begin{centering}
 $(4_{P}, 1_{P}, P_{A}=\{p_{2, 1}, p_{1, 3}, p_{3, 6}\}, P_{B}=\{p_{3,2}, p_{2, 4}, p_{4, 5}\})$\par\bigskip
 \end{centering}
  and the matrix $P$ is partitioned according to the word $ABABBA$.

\addtocounter{corol}{7}

\begin{corol}\label{invariant2}

\textit{Let $W_{P}= l_{m_{1}}\dots l_{m_{2n}}$ and $W'_{P}= l'_{m_{1}}\dots l'_{m_{2n}}$;\quad $l_{m_{n}}, l'_{m_{n}}\in \{A,B\}$ be two words associated to a matrix $P$ of type $\mathscr{H}_{n}$ with corresponding sets of helices $H_{P}=(i_{P}, j_{P}, P_{A}, P_{B})$, and $H'_{P}=(f'_{P}, g'_{P}, P'_{A}, P'_{B})$. Then $|H_{P}|=|H'_{P}|$.}
\end{corol}
\textbf{Proof.} The result follows from Theorem \ref{invariant11} by replacing, $P', A',B', P'_{A'}, P'_{B'}, f_{P'}$ and $g_{P'}$ for $P, A, B, P'_{A}, P'_{B}, f'_{P}$ and $g'_{P}$, respectively.  \hspace{0.5cm}$\square$

\begin{corol}\label{invariant3}
\textit{If for $n\geq 1$, $P$ and $P'$ are equivalent preprojective Kronecker modules with dimension vector of the form $[n+1\hspace{0.2cm}n]$ and corresponding sets of helices $H_{P}$ and $H_{P'}$ then $|H_{P}|=|H_{P'}|$.}

\end{corol}

\textbf{Proof.} Matrix presentations of preprojective Kronecker modules $P$ and $P'$  are both of type $\mathscr{H}_{n}$ defined by words of the form $AA\dots ABB\dots BB$. \hspace{0.5cm}$\square$

\addtocounter{teor}{2}

\begin{teor}\label{Kronecker0}
\textit{If for $n\geq1$, $P$ denotes a preprojective Kronecker module then the number of helices associated to $P$ is $h^{P}_{n}=n!\lceil\frac{n}{2}\rceil$ where $\lceil x \rceil$ denotes the smallest integer greatest than $x$}.

\end{teor}

\textbf{Proof.}  According to Theorem \ref{invariant11} and Corollary \ref{invariant3}, it suffices to determine the number of helices $|(1_{p}, (n+1)_{p}, p_{A}=\{a_{i+1,i}\mid 1\leq i\leq n\}, p_{B}=\{b_{i, i}\mid 1\leq i\leq n\})|$ associated to preprojective Kronecker modules $p=(n+1,n)$ of type III and words of the form $W_{p}=AA\dots ABB\dots B$. \par\bigskip

Firstly, we note that there is only one helix associated to the indecomposable preprojective modules $(2,1)$ and $(3,2)$. And the vertices sequence of helices associated to the indecomposable $(4,3)$ with $a_{1,j}$ fixed are:

\begin{equation}
\begin{split}
hl_{1} &=\{a_{1,j}, b_{1,1},b_{2,1},a_{2,1},a_{3,1},b_{3,3},b_{4,3},a_{4,3}\},\\
hl_{2}&= \{a_{1,j}, b_{1,1},b_{3,1},a_{3,2},a_{2,2},b_{2,2},b_{4,2},a_{4,3}\},\\
hl_{3}&=\{a_{1,j}, b_{1,1},b_{4,1},a_{4,3},a_{3,3},b_{3,3},b_{2,3},a_{2,1}\},\\
hl_{4}&=\{a_{1,j}, b_{1,1},b_{4,1},a_{4,3},a_{2,3},b_{2,2},b_{3,2},a_{3,2}\}.\\
\end{split}  
\end{equation}

 The number of helices is given by the number of vertices at the last level of the following associated tree:    
 
 
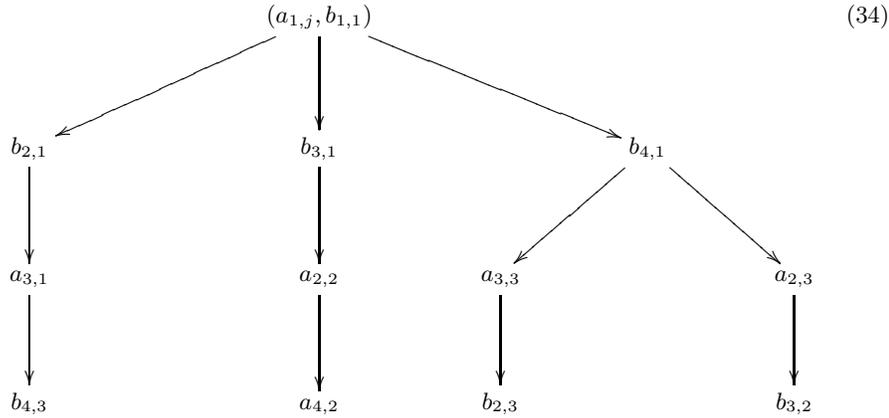
\begin{equation}\label{Ks}
 \xymatrix@=35pt{
&&\ar@[black][lld]\ar@[black][d]\ar@[black][drr](a_{1, j}, b_{1,1})&&&\\
\ar@[black][d]b_{2, 1}&&\ar@[black][d]b_{3, 1}&&\ar@[black][ld]\ar@[black][rd]b_{4,1}&\\
\ar@[black][d]a_{3, 1}&&\ar@[black][d]a_{2, 2}&\ar@[black][d]a_{3,3}&&\ar@[black][d]a_{2,3}\\
b_{4, 3}&&a_{4,2}&b_{2,3}&&b_{3,2}
}
\end{equation}

 Suppose now that the result is true for any indecomposable preprojective Kronecker module $(t+1,t)$, $1\leq t<n$ then we can see that in general the rooted tree $T_{n}$ associated to the indecomposable preprojective Kronecker module $(n+1,n)$ has the following characteristics bearing in mind that vertex $b_{1,1}$ gives the root node $a^{0}_{1}$:  
\begin{enumerate}[$(a)$]
\item $a^{0}_{1}$ has $n$ children enumerated from the left to the right as $(a^{1}_{1}, a^{1}_{2},\dots, a^{1}_{n})$,
\item For $1\leq i\leq n-1$ each vertex $a^{1}_{i}$ has $n-2$ children enumerated from the left to the right as $(a^{1}_{i,1}, a^{1}_{i,2},\dots, a^{1}_{i, n-2})$ whereas vertex $a^{1}_{n}$ has $n-1$ children of the form $(a^{1}_{n,1},a^{1}_{n,2},\dots, a^{1}_{n, n-1})$, each children of a vertex $a^{1}_{n, l_{1}}$, $1\leq l_{1}\leq n-1$ has $n-2$ children $a^{1}_{n, l_{1},l_{2}}$ with $1\leq l_{2}\leq n-2$, in general for this particular tree a vertex $a^{1}_{n, l_{1},l_{2},l_{3},\dots, l_{t}}$ has $n-(t+1)$ children, $1\leq t\leq n-2$. Note that the number of vertices at the last level of the rooted tree $T'_{n}$ with $a^{1}_{n}$ as root node is $(n-1)!$,
\item For each $h$, $1\leq h\leq n-2$, vertex $a^{1}_{i, h}$ is a root node of the tree $T_{n-2}$.

\end{enumerate}
The following picture shows the general structure of the rooted tree $T_{n}$
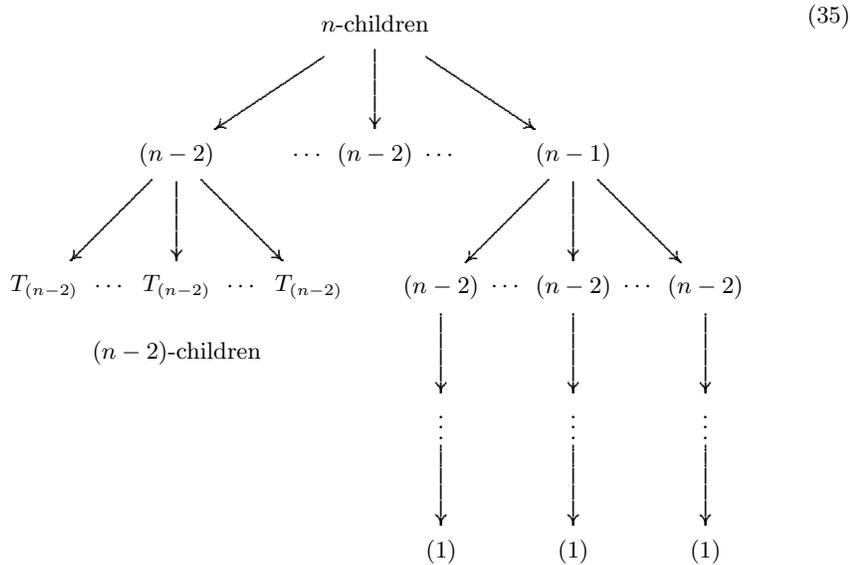
\begin{equation}
\begindc{\commdiag}[250]
\obj(0,0)[n00]{$n$-children}
\obj(-3,-2)[n11]{$(n-2)$}
\obj(-1,-2)[n21]{$\cdots$}
\obj(0,-2)[n31]{$(n-2)$}
\obj(1,-2)[n41]{$\cdots$}
\obj(3,-2)[n51]{$(n-1)$}
\obj(-5,-4)[n12]{$T_{(n-2)}$}
\obj(-4,-4)[n22]{$\cdots$}
\obj(-3,-4)[n32]{$T_{(n-2)}$}
\obj(-2,-4)[n42]{$\cdots$}
\obj(-1,-4)[n52]{$T_{(n-2)}$}
\obj(1,-4)[n62]{$(n-2)$}
\obj(1,-6)[n84]{$\vdots$}
\obj(1,-8)[n86]{$(1)$}
\obj(2,-4)[n72]{$\cdots$}
\obj(3,-4)[n82]{$(n-2)$}
\obj(3,-6)[n83]{$\vdots$}
\obj(3,-8)[n87]{$(1)$}
\obj(4,-4)[n92]{$\cdots$}
\obj(5,-4)[n102]{$(n-2)$}
\obj(5,-6)[n85]{$\vdots$}
\obj(5,-8)[n88]{$(1)$}
\obj(-3,-5)[p0]{$(n-2)$-children}
\mor{n00}{n11}{}
\mor{n00}{n31}{}
\mor{n00}{n51}{}
\mor{n11}{n12}{}
\mor{n11}{n32}{}
\mor{n11}{n52}{}
\mor{n51}{n62}{}
\mor{n51}{n82}{}
\mor{n51}{n102}{}
\mor{n85}{n88}{}
\mor{n83}{n87}{}
\mor{n84}{n86}{}
\mor{n82}{n83}{}
\mor{n102}{n85}{}
\mor{n62}{n84}{}
\enddc
\end{equation}

  According to the rules $(a)-(c)$ the number of vertices $L_{T_{n}}$ at the last level of the tree $T_{n}$ is given by the formula

  \begin{equation}
  \begin{split}
 L_{T_{n}}&=(n-1)(n-2)L_{T_{n-2}}+L(T'_{n})=(n-1)(n-2)\frac{h^{p}_{n-2}}{n-2}+(n-1)!\\
& =(n-1)!\lceil\frac{n}{2}\rceil= \frac{h^{p}_{n}}{n}.  
  \end{split}
  \end{equation}
We are done.\hspace{0.5cm}$\square$

\addtocounter{Nota}{5}
\begin{Nota}\label{A052558}
Sequence A052558 is categorified via the number of helices associated to preprojective Kronecker modules, if in the condition $(e)$ of its definition, it is assumed that the starting vertex is fixed. Without such fixing condition the number of helices associated to a preprojective Kronecker module is given in the sequence encoded as A052591 in the OEIS.   

\end{Nota}

\subsection{Sequence A052591 Via Brauer Configuration Algebras}

In this section, categorification of elements of the integer sequence A052591 are given via the number of summands in the heart of indecomposable projective modules over the Brauer configuration algebra $\Lambda_{K^n}$ defined by the Brauer configuration $K^{n}=(K^{n}_{0},K^{n}_{1}, \mu, \mathcal{O})$ with the following properties for $n\geq3$ fixed:

\begin{enumerate}
\item \begin{equation}\label{sequence1}
\begin{split}
K^n_{0}&=\{x_{1},x_{2}\},\\
K^n_{1}&=\{V_{t}=x^{(2t+2)!}_{1}x_{2}^{({(t)(2t+2)!})}\}_{1\leq t\leq n}. 
\end{split}
\end{equation}

\item The orientation $\mathcal{O}$ is defined in such a way that for $t\geq1$

 \begin{equation}
\begin{split}
\mathrm{at \hspace{0.1cm}vertex}\hspace{0.1cm}x_1;&\hspace{0.1cm}V^{(4!)}_{1}< V^{(6!)}_{2}< V^{(8!)}_{3}<\dots< V^{((2n+2)!)}_{n},\\
\mathrm{at \hspace{0.1cm}vertex}\hspace{0.1cm}x_2;&\hspace{0.1cm}V^{2(12)}_{1}< V^{2(720)}_{2}<V^{2(60480)}_{3}< \dots<V^{(({(n)(2n+2)!}))}_{n}.                
\end{split}
\end{equation}

\item The multiplicity function $\mu$ is such that $\mu(x_{1})=\mu(x_{2})=1$.
\end{enumerate}
Where the symbol $x^{(j)}_{i}$ in a given polygon $V_{t}$ means that $\mathrm{occ}(x_{i}, V_{t})=j$. Note that, the specializations $x_{1}=1$ and $x_{2}=2$ allows to describe polygons $V_{t}$ as integer partitions of numbers in the sequence A052591 (see identity (\ref{word})).\par\bigskip

The following is the Brauer quiver $Q_{K^{n}}$ associated to this configuration (numbers $n_{t_{1}} (n_{t_{2}})$ attached to the loops denote the occurrence of the vertex ($x_{1}$ above, $x_{2}$ below) in the corresponding polygon $V_{t}$, $1\leq t\leq n$), $c^{1}_{i_{t}} (c^{2}_{j_{t}})$ denotes a set of loops associated to the vertex $V_{t}$, $|c^{1}_{i_{t}}|=l^{1}_{t}= n_{t_{1}}-1 $,\hspace{0.2cm} $ |c^{2}_{j_{t}}|=l^{2}_{t}= n_{t_{2}}-1$.

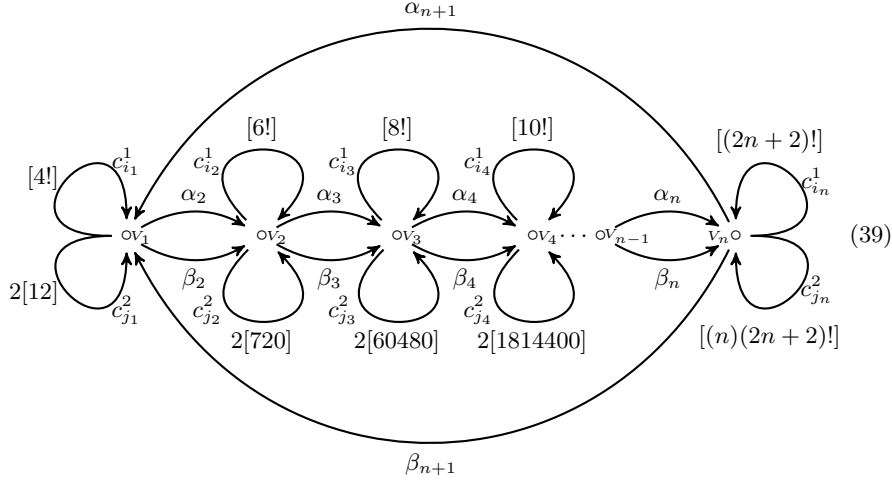
\begin{center}
\begin{tikzpicture}[node distance=2mm,->,>=stealth',shorten >=1pt,thick,scale=0.9]

\node at (-1.79,2.98) {\tiny{$V_1$}};
\node at (0.22,2.98) {\tiny{$V_2$}};
\node at (2.22,2.98){\tiny{$V_3$}};
\node at (4.22,2.98){\tiny{$V_4$}};

\node at (5.4,2.98){\tiny{$V_{n-1}$}};
\node at (6.75,2.98){\tiny{$V_n$}};

\node at (-2,4.1)(0){$c_{i_{{1}}}^1$};
\node at (-0.8,4.1)(0){$c_{i_{{2}}}^1$};
\node at (1.2,4.1)(0){$c_{i_{{3}}}^1$};
\node at (3.2,4.1)(0){$c_{i_{{4}}}^1$};
\node at (8.2,3.8)(0){$c_{i_{{n}}}^1$};
\node at (7.5,4.4)(0){$[(2n+2)!]$};

\node at (-2,1.9)(0){$c^2_{j_1}$};
\node at (-0.8,1.9)(0){$c_{j_{{2}}}^2$};
\node at (1.2,1.9)(0){$c_{j_{{3}}}^2$};
\node at (3.2,1.9)(0){$c_{j_{{4}}}^2$};
\node at (8.2,2.2)(0){$c_{j_{{n}}}^2$};

\node at (7.5,1.5)(0){$\left[ {(n)(2n+2)!}\right]$};

\node[inner sep=2.8pt,label=below:] at (-2,3)(1){$\circ$};
\node[inner sep=2.8pt,label=below:] at (0,3)(2){$\circ$};
\node[inner sep=2.8pt,label=below:] at (2,3)(3){$\circ$};
\node[inner sep=2.8pt,label=below:] at (4,3)(4){$\circ$};
\node[inner sep=2.8pt] at (4.6,3)(7){$\ldots$};
\node[inner sep=2.8pt] at (5,3)(8){$\circ$};
\node[inner sep=2.8pt,label=below:] at (7,3)(5){$\circ$};


\path[->](1) [bend left] edge node[above] {$\alpha_2$} (2);
\path[->](2) [bend left] edge node[above] {$\alpha_3$}(3);
\path[->](3) [bend left] edge node[above] {$\alpha_4$}(4);
\path[->](8) [bend left] edge node[above] {$\alpha_{n}$}(5);
\draw[->] (5).. controls ($ (5) -(2,4)$) and ($ (1) +(2,-4)$) .. (1) node[pos=0.5, inner sep=-1pt, label=below:{$\beta_{n+1}$}] {};

\draw[->] (5).. controls ($ (5) +(-2,4)$) and ($ (1) -(-2,-4)$) .. (1) node[pos=0.5, inner sep=-1pt, label=above:{$\alpha_{n+1}$}] {};


\path[->](1) [bend right] edge node[below] {$\beta_2$} (2);
\path[->](2) [bend right] edge node[below] {$\beta_3$} (3);
\path[->](3) [bend right] edge node[below] {$\beta_4$} (4);
\path[->](8) [bend right] edge node[below] {$\beta_{n}$} (5);

\Loop[dist=2cm,dir=SOWE,style={thick},label= {$2[12]$},labelstyle=left](1)
\Loop[dist=2cm,dir=SO,style={thick},label= {$2[720]$},labelstyle=below](2)
\Loop[dist=2cm,dir=SO,style={thick},label= {$2[60480]$},labelstyle=below](3)
\Loop[dist=2cm,dir=SO,style={thick},label={$2[1814400]$},labelstyle=below](4)
\Loop[dist=2cm,dir=SOEA,label= {},labelstyle=right](5)


\Loop[dist=2cm,dir=NO,label= {$[6!]$},labelstyle=above](2)
\Loop[dist=2cm,dir=NO,label= {$[8!]$},labelstyle=above](3)
\Loop[dist=2cm,dir=NO,label= {$[10!]$},labelstyle=above](4)
\Loop[dist=2cm,dir=NOWE,label={$[4!]$},labelstyle=left](1)
\Loop[dist=2cm,dir=NOEA, style={thick},label={},labelstyle=right](5)

\node at (9,3)(0){$(39)$};
\end{tikzpicture}

\end{center}

The admissible ideal ${I}$ is generated by the following relations (in this case, if there are associated $l^{1}_{t}$ ($l^{2}_{t}$) loops at the vertex $V_{t}$ associated to $x_{1}$ (associated to $x_{2}$)  then we let $P^{j}_{t}$ denote the product of $j\leq l^{m}_{t}$ loops, $m\in\{1,2\}$), $c^m_{h_s}$ is a notation for a set of cycles $\{c^m_{h_{s,1}}, c^m_{h_{s,2}},\dots, c^m_{h_{s,l^{m}_{s}}}, m\in\{1,2\}, h\in\{i,j\}, s\in\{1,2,\dots,n\}\}$:
\begin{enumerate}
\item $c^{1}_{i_{s,x}}c^{1}_{i_{s,y}}-c^{1}_{i_{s,y}}c^{1}_{i_{s,x}}$, for all possible values of $i,s,x,y$,
\item $c^{2}_{j_{s,x}}c^{2}_{j_{s,y}}-c^{2}_{j_{s,y}}c^{2}_{j_{s,x}}$, for all possible values of $i,s,x,y$,
\item  $c^{1}_{i_{s,x}}c^{2}_{j_{s,y}}$   \quad and\quad  $c^{2}_{j_{s,x}}c^{1}_{i_{s,y}}$  , for all possible values of $i,s,x,y$,
\item $c^{1}_{i_{s,x}}\beta_{s+1}$;   $c^{2}_{j_{s,y}}\alpha_{s+1}$; $\beta_{s}c^{1}_{i_{s,x}}$; $\alpha_{s}c^{2}_{j_{s,x}}$,  for all possible values of $i,s,x,y$,
\item   $(c^{1}_{i_{s,x}})^2$;\quad $(c^{2}_{j_{s,y}})^2$, for all possible values of $i,s,x,y$,
\item $\alpha_{t}\alpha_{t+1}$;\quad$\alpha_{n+1}\alpha_{2}$;\quad$\beta_{t}\beta_{t+1}$;\quad$\beta_{n+1}\beta_{2}$$;\quad\alpha_{t}\beta_{t+1}$;\quad$\beta_{j}\alpha_{j+1}$;\quad $\alpha_{n+1}\beta_{2}$;\quad$\beta_{n+1}\alpha_{2}$, for all possible values of $j,t$,
\item $\alpha_{i}P^{j}_{i}\gamma_{i+1}$;\quad $\alpha_{n+1}P^{j}_{1}\gamma_{2}$;\quad$\beta_{t}P^{h}_{t}\gamma_{t+1}$;\quad $\beta_{n+1}P^{h}_{1}\gamma_{2}$;\quad$0<j<l^{1}_{i}$, $0<h<l^{2}_{t}$, $1\leq i,t\leq n$, $\gamma\in\{\alpha,\beta\}$,
\item For all the possible products (special cycles) of the form:
\addtocounter{equation}{1}
\begin{equation}
\begin{split}
\varepsilon^{1}_{1}&=\alpha_{t}P^{l^{1}_{t}}_{t}\alpha_{t+1}P^{l^{1}_{t+1}}_{t+1}\dots\alpha_{n}P^{l^{1}_{n}}_{n}\alpha_{n+1}P^{l^{1}_{1}}_{1}\dots\alpha_{t-1}P^{l^{1}_{t-1}}_{t-1},\\
\varepsilon^{2}_{1}&=P^{m}_{t-1}\alpha_{t}P^{l^{1}_{t}}_{t}\alpha_{t+1}P^{l^{1}_{t+1}}_{t+1}\dots\alpha_{n}P^{l^{1}_{n}}_{n}\alpha_{n+1}P^{l^{1}_{1}}_{1}\dots\alpha_{t-1}P^{l^{1}_{t-1}-j}_{t-1},\\
\varepsilon^{3}_{2}&=\beta_{t}P^{l^{2}_{t}}_{t}\beta_{t+1}P^{l^{2}_{t+1}}_{t+1}\dots\beta_{n}P^{l^{2}_{n}}_{n}\beta_{n+1}P^{l^{2}_{1}}_{1}\dots\beta_{t-1}P^{l^{2}_{t-1}}_{t-1},\\
\varepsilon^{4}_{2}&=P^{h}_{t-1}\beta_{t}P^{l^{2}_{t}}_{t}\beta_{t+1}P^{l^{2}_{t+1}}_{t+1}\dots\beta_{n}P^{l^{2}_{n}}_{n}\beta_{n+1}P^{l^{2}_{1}}_{1}\dots\beta_{t-1}P^{l^{2}_{t-1}-h}_{t-1},
\end{split}
\end{equation}

relations of the form $\varepsilon^{r}_{i}-\varepsilon^{s}_{j}$, $r, s\in\{1,2,3,4\}$, $i,j\in\{1,2\}$ take place. Note that, products of the form $P^0_{t-1}$ correspond to suitable orthogonal primitive idempotents $e_{t}$, $1\leq t\leq n$, 

\item  $\varepsilon^{1}_{1}\alpha_{t}$,  $\varepsilon^{3}_{2}\beta_{t}$.

\end{enumerate}

\par\bigskip
The following result holds for indecomposable projective modules over the algebra $\Lambda_{K^n}$.

\addtocounter{corol}{2}

\begin{corol}\label{Kronecker1}
\textit{For $n\geq3$ fixed and $1\leq t\leq n$, the number of summands in the heart of the indecomposable projective representation $V_{t}$ over the Brauer configuration algebra $\Lambda_{K^n}$ equals the number of helices associated to the preprojective Kronecker module $(2t+3,2t+2)$, $1\leq t\leq n$.}

\end{corol}

\textbf{Proof.} Firstly we note that for any $t$, $\mathrm{rad}^2\hspace{0.1cm}V_{t}\neq0$. Thus according to the Theorem \ref{multiserial} the number of summands in the heart of any of the indecomposable projective modules $V_{t}$ equals $\mathrm{occ}(x_{1},V_{t})+\mathrm{occ}(x_{2},V_{t})=(2t+2)!+t(2t+2)!=h^{p}_{2t+2}=h^{(2t+3, 2t+2)}_{2t+2}$ which is the number of helices associated in a unique form to the indecomposable preprojective Kronecker module $(2t+3, 2t+2)$. We are done.\hspace{0.5cm}$\square$
\par\bigskip

The following results regard the dimension of algebras of type $\Lambda_{K^n}$.

\begin{corol}\label{Kronecker2}
\textit{For $n\geq3$ fixed, it holds that $\frac{1}{2}(\mathrm{dim}_{k}\hspace{0.1cm}\Lambda_{K^{n}})=n+t_{\gamma_{n}-1}+t_{\delta_{n}-1}$, where $\gamma_{n}=\underset{m=1}{\overset{n}{\sum}}m(2m+2)!$, $\delta_{n}=\underset{m=1}{\overset{n}{\sum}}(2m+2)!$, and $t_{h}$ denotes the $h$th triangular number.}

\end{corol}

\textbf{Proof.} Proposition \ref{dimension} allows to conclude that  $\mathrm{dim}_{k}\hspace{0.1cm}\Lambda_{K^n}/{I}=2n+\underset{i=1}{\overset{2}{\sum}}|C_{i}|(|C_{i}|-1)$ where for each $i=1,2$, $|C_{i}|=val({x_{i}})$. The theorem holds taking into account that  for any $j\geq 2$, $j(j-1)=2t_{j-1}$.\hspace{0.5cm}$\square$

\begin{corol}\label{Kronecker3}
\textit{For $n\geq3$ fixed, it holds that $\mathrm{dim}_{k}\hspace{0.1cm}Z(\Lambda_{K^{n}})=n-1+\underset{t=1}{\overset{n}{\sum}}h^{p}_{2t+2}$.}

\end{corol}

\textbf{Proof.} Since $\mathrm{rad}^2\hspace{0.1cm}\Lambda_{K^n}\neq0$, the result is a consequence of Theorem \ref{Serra} with $\mu(x_{1})=\mu(x_{2})=1$, $|K^{n}_{0}|=2$,
$|K^{n}_{1}|=n$ and $\mathrm{occ}(x_{1},V_{t})+\mathrm{occ}(x_{2},V_{t})=h^{p}_{2t+2}$.\hspace{0.5cm}$\square$

\addtocounter{Nota}{3}

\begin{Nota}
Similar results as in Corollaries \ref{Kronecker1}-\ref{Kronecker3} can be obtained for preprojective Kronecker modules of the form $(4t+2,4t+1), t\geq1$ by considering in the original Brauer configuration that
\item \begin{equation}\label{sequence2}
\begin{split}
K^n_{0}&=\{x_{1},x_{2}\},\\
K^n_{1}&=\{V_{t}=x^{(4t+1)!}_{1}x_{2}^{{2t(4t+1)!})}\}_{1\leq t\leq n},
\end{split}
\end{equation}
and keeping the relations in the quiver without changes (bearing in mind of course the new occurrences of the vertices for the different products). In particular, it holds that
$\mathrm{dim}_{k}\hspace{0.1cm}Z(\Lambda_{K^{n}})=n-1+\underset{t=1}{\overset{n}{\sum}}h^{p}_{4t+1}.$
\end{Nota}

\section{Categorification of the Sequence A100705 and Some Related Integer Sequences}
In this section elements of the integer sequence $h_{n}=n^{3}+(n+1)^{2}$, $n\geq 1$ are interpreted as the number of cycles associated to preprojective representations of the tetrad, such interpretation allows to categorify this integer sequence (encoded in the OEIS as A100705). Besides, some new Brauer configuration algebras are defined in order to get alternative categorifications of this sequence and some additional integer sequences. We note that the Brauer configuration (\ref{sequence}) allows to see each polygon $V_{n}$ as a partition of the number $h_{n}$ into two parts of the form $\{n, n+1\}$ where $n$ occurs $(n)^{2}$ times and $n+1$ occurs $n+1$ times. Assuming the classical notation for partitions \cite{Andrews} each number $h_{n}$ can be expressed as follows:

\begin{equation}\label{equation}
\begin{split}
&(n)^{(n^{2})} (n+1)^{(n+1)},\quad n\geq1,
\end{split}
\end{equation}

we let $\mathcal{P}_{n}$ denote such a partition. The partition tree $T_{\mathcal{P}_{n}}$ associated to each partition of the form $\mathcal{P}_{n}$ is obtained by assuming the notation:

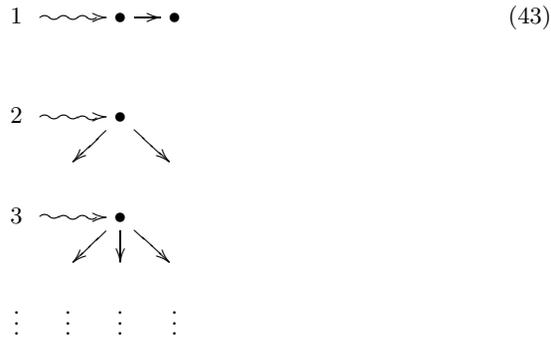
\begin{equation}
\xymatrix@=10pt{
1\ar@{~>}[rr]&&\bullet\ar[r]&\bullet\\
&&&\\
2\ar@{~>}[rr]&&\ar[ld]\bullet\ar[rd]&\\
&&&\\
3\ar@{~>}[rr]&&\ar[ld]\ar[d]\bullet\ar[rd]&\\
&&&\\
\vdots&\vdots&\vdots&\vdots
}
\end{equation}

In this case, $T_{\mathcal{P}_{n}}$ has a root node with $n+1$ children, $n$ of them have $n$ children and the last one has $n+1$ children in such a way that in the last level of $T_{\mathcal{P}_{n}}$, $n$ of these children represent a partition of the form $(n)^{(n-1)}(n+1)^{(1)}$ and the last one represents a partition of the form $(n)^{(n)}(n+1)^{(1)}$. Partition trees of the form $T_{\mathcal{P}_{n}}$ are used in the proof of Theorem \ref{projective quadruple}.\par\bigskip

The following results regard the number of cycles associated to matrices of type $\mathscr{C}_{n}$. In particular to preprojective representations of the tetrad of type IV.

\par\bigskip
 The next Theorem \ref{invariant12} can be proved by using similar arguments as those posed in the proof of Theorem \ref{invariant11}.
\addtocounter{teor}{5}
\begin{teor}\label{invariant12}
\textit{For $n\geq2$, let, $(U,U_{i},n\mid 1\leq i\leq 4)$, $(U',U'_{i},n\mid 1\leq i\leq 4)$, $\mathcal{C}_{U}$ and $\mathcal{C}_{U'}$ be two matrices of type $\mathscr{C}_{n}$ with corresponding sets of cycles $\mathcal{C}_{U}$ and $\mathcal{C}_{U'}$ defined by systems of the form $(i^{R}_{A_{1}}, j^{R}_{A_{4}},j'^{R}_{A'_{4}}, P_{A^{r}_{i}}\mid 1\leq i\leq 4, r\in\{0,1\})$ and $(f^{R}_{B_{1}}, g^{R}_{B_{4}},g'^{R}_{B'_{4}}, P_{B^{r}_{i}}\mid 1\leq i\leq 4, r\in\{0,1\})$, respectively. Then $|\mathcal{C}_{U}|=|\mathcal{C}_{U'}|$}.
\end{teor}

\textbf{Proof.} To each cycle $\mathcal{C}\in(i^{R}_{A_{1}}, j^{R}_{A_{4}},j'^{R}_{A'_{4}}, P_{A^{r}_{i}}\mid 1\leq i\leq 4, r\in\{0,1\})$, it is possible to build a unique copy $\mathcal{C}'\in (f^{R}_{B_{1}}, g^{R}_{B_{4}},g'^{R}_{B'_{4}}, P_{B^{r}_{i}}\mid 1\leq i\leq 4, r\in\{0,1\}))$  by using the same procedures described in the proof of Theorem \ref{invariant11}.\hspace{0.5cm}$\square$\par\bigskip

The next corollary shows that the number of cycles associated to an indecomposable representation of type IV of the tetrad is invariant under admissible matrix transformations.
\addtocounter{corol}{2}
\begin{corol}
\textit{If for $n\geq 2$, $U$ and $U'$ are equivalent preprojective representations of the tetrad of type IV with corresponding sets of cycles $\mathcal{C}_{U}$ and $\mathcal{C}_{U'}$ then $|\mathcal{C}_{U}|=|\mathcal{C}_{U'}|$.}

\end{corol}

\textbf{Proof.} The matrix presentations of preprojective representations of type IV of the tetrad are of type $\mathscr{C}_{n}$.\hspace{0.5cm}$\square$\par\bigskip

\addtocounter{teor}{1}

\begin{teor}\label{projective quadruple}

\textit{For $n\geq 2$ fixed and $1\leq i\leq n$ the number of cycles associated to  an indecomposable preprojective  representation of type $\mathrm{IV}$ and order $i+1$ is $h_{i}=i^{3}+ (i+1)^{2}$.}

\end{teor}

\textbf{Proof.} Fix a preprojective representation $U^{n}$ of type IV and order $n\geq 2$, and denote its different blocks as follows:

\begin{center}

$U^{n}=\begin{array}{|c|c|c|c|}\hline A & B & C & D  \\\hline A' & B' & C' & D'\\\hline\end{array}$

\end{center}

We note that all the cycles associated to $U^{n}$ can be seen as trees $T_{c_{(n+1),(n+1)}}$ which have the entry $c_{(n+1),(n+1)}$ as root node with $n$ branches whose successors are given by entries

\begin{center}
 $c'_{1,(n+1)}, c'_{2,(n+1)},\dots, c'_{n,(n+1)}$.
\end{center}
 
Each entry $c'_{i,(n+1)}$ has $n-1$ branches if $i\neq n$, whereas $c'_{n,(n+1)}$ has $n$ branches. Besides, all of these entries give rise to an arrow\par\bigskip

\begin{center}

$c'_{i,(n+1)}\rightarrow d_{j,i}$,

\end{center}

for some entry $d_{j, i}\in D$. Actually,  $d_{j, i}$ is a successor root of $c'_{i,(n+1)}$ with $(n-1)$ branches in the tree whenever $j\in\{1,\dots, n\}$ and $i\neq 1$. If $i=1$ then  $d_{1,j}$ has by construction $n$ branches in $C'$. Therefore, the structure of
$T_{c_{(n+1),(n+1)}}$ has the following shape:

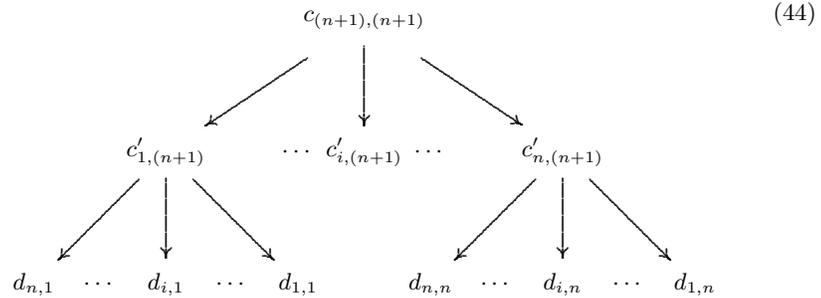
\begin{equation}
\begindc{\commdiag}[250]
\obj(0,0)[n00]{$c_{(n+1),(n+1)}$}
\obj(-3,-2)[n11]{$c'_{1,(n+1)}$}
\obj(-1,-2)[n21]{$\cdots$}
\obj(0,-2)[n31]{$c'_{i,(n+1)}$}
\obj(1,-2)[n41]{$\cdots$}
\obj(3,-2)[n51]{$c'_{n,(n+1)}$}
\obj(-5,-4)[n12]{$d_{n,1}$}
\obj(-4,-4)[n22]{$\cdots$}
\obj(-3,-4)[n32]{$d_{i,1}$}
\obj(-2,-4)[n42]{$\cdots$}
\obj(-1,-4)[n52]{$d_{1,1}$}
\obj(1,-4)[n62]{$d_{n,n}$}
\obj(2,-4)[n72]{$\cdots$}
\obj(3,-4)[n82]{$d_{i,n}$}
\obj(4,-4)[n92]{$\cdots$}
\obj(5,-4)[n102]{$d_{1,n}$}
\mor{n00}{n11}{}
\mor{n00}{n31}{}
\mor{n00}{n51}{}
\mor{n11}{n12}{}
\mor{n11}{n32}{}
\mor{n11}{n52}{}
\mor{n51}{n62}{}
\mor{n51}{n82}{}
\mor{n51}{n102}{}
\enddc
\end{equation}

Which corresponds to the partition tree $T_{\mathcal{P}_{(n-1)}}$ of $h_{(n-1)}=(n-1)^{3}+(n)^{2}$. \hspace{0.5cm}$\square$\par\bigskip

As an example the following is the diagram of $T_{c_{4,4}}$ such that the number of vertices in the last level gives the number of associated cycles (described in the proof of Theorem \ref{projective quadruple}) to the indecomposable representation of the tetrad $U^{3}$:

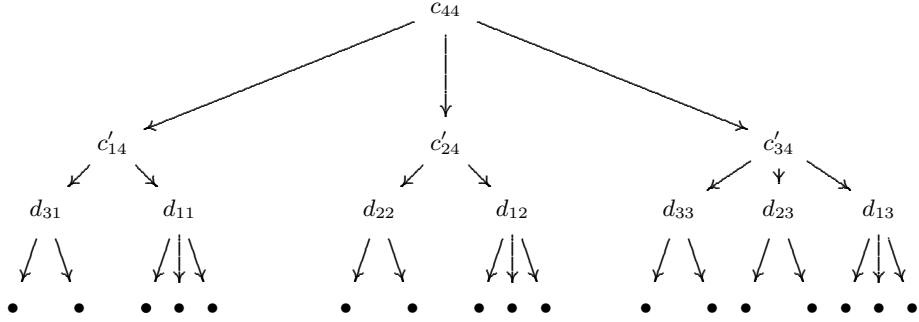
\begin{equation}
\begindc{\commdiag}[126]
\obj(0,0)[p0]{$c_{44}$}
\obj(-10,-4)[p1]{$c'_{14}$}
\obj(0,-4)[p2]{$c'_{24}$}
\obj(10,-4)[p3]{$c'_{34}$}
\obj(-12,-6)[p4]{$d_{31}$}
\obj(-8,-6)[p5]{$d_{11}$}
\obj(-2,-6)[p6]{$d_{22}$}\notag
\obj(2,-6)[p7]{$d_{12}$}
\obj(7,-6)[p8]{$d_{33}$}
\obj(10,-6)[p9]{$d_{23}$}
\obj(13,-6)[p10]{$d_{13}$}
\obj(-13,-9)[p11]{$\bullet$}
\obj(-11,-9)[p12]{$\bullet$}
\obj(-9,-9)[p13]{$\bullet$}
\obj(-8,-9)[p14]{$\bullet$}
\obj(-7,-9)[p15]{$\bullet$}
\obj(-3,-9)[p16]{$\bullet$}
\obj(-1,-9)[p17]{$\bullet$}
\obj(-9,-9)[p18]{$\bullet$}
\obj(1,-9)[p20]{$\bullet$}
\obj(2,-9)[p21]{$\bullet$}
\obj(3,-9)[p22]{$\bullet$}
\obj(6,-9)[p23]{$\bullet$}
\obj(8,-9)[p24]{$\bullet$}
\obj(9,-9)[p25]{$\bullet$}
\obj(11,-9)[p26]{$\bullet$}
\obj(12,-9)[p27]{$\bullet$}
\obj(13,-9)[p28]{$\bullet$}
\obj(14,-9)[p29]{$\bullet$}
\mor{p0}{p1}{}
\mor{p0}{p2}{}
\mor{p0}{p3}{}
\mor{p1}{p4}{}
\mor{p1}{p5}{}
\mor{p2}{p6}{}
\mor{p2}{p7}{}
\mor{p3}{p8}{}
\mor{p3}{p9}{}
\mor{p3}{p10}{}
\mor{p4}{p11}{}
\mor{p4}{p12}{}
\mor{p5}{p13}{}
\mor{p5}{p14}{}
\mor{p5}{p15}{}
\mor{p6}{p16}{}
\mor{p6}{p17}{}
\mor{p7}{p20}{}
\mor{p7}{p21}{}
\mor{p7}{p22}{}
\mor{p8}{p23}{}
\mor{p8}{p24}{}
\mor{p9}{p25}{}
\mor{p9}{p26}{}
\mor{p10}{p27}{}
\mor{p10}{p28}{}
\mor{p10}{p29}{}
\enddc
\end{equation}  
\par\bigskip
The number of cycles associated to the indecomposable preprojective representation of the tetrad $U^{3}$ (shown below) equals the second term of the integer sequence A100705. Actually, the number of cycles associated to $U^{n}$ is given by $h_{(n-1)}=(n-1)^{3}+(n)^{2}$, $n\geq2$, which is the $(n-1)$th term of this sequence. Black arrows denote the common part of all these cycles.

\par\bigskip

\hspace{5mm}
\setlength{\unitlength}{0.97pt}
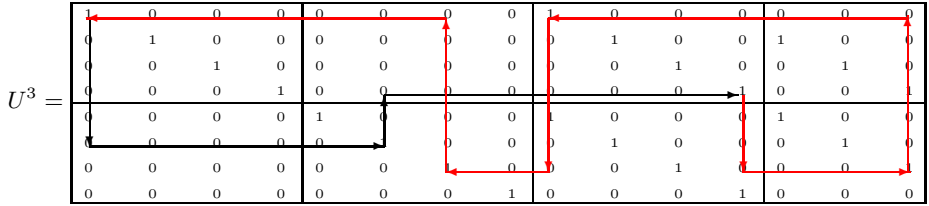
\begin{picture}(90,90)
 
 \multiput(-25,37)(0,30){1}{$U^{3}=$}
 \multiput(0,0)(90,0){4}{\line(0,1){78}}
 \multiput(333,0)(90,0){1}{\line(0,1){78}}        
  \multiput(0,0)(0,39){3}{\line(1,0){333}}
   \multiput(5,72)(30,0){1}{\tiny 1}
    \multiput(30,72)(30,0){1}{\tiny 0}
  \multiput(55,72)(30,0){1}{\tiny 0}
    \multiput(80,72)(30,0){1}{\tiny 0}
      \multiput(5,62)(30,0){1}{\tiny 0}
    \multiput(30,62)(30,0){1}{\tiny 1}
  \multiput(55,62)(30,0){1}{\tiny 0}
    \multiput(80,62)(30,0){1}{\tiny 0}
 \multiput(5,52)(30,0){1}{\tiny 0}
    \multiput(30,52)(30,0){1}{\tiny 0}
  \multiput(55,52)(30,0){1}{\tiny 1}
    \multiput(80,52)(30,0){1}{\tiny 0}
       \multiput(5,42)(30,0){1}{\tiny 0}
    \multiput(30,42)(30,0){1}{\tiny 0}
  \multiput(55,42)(30,0){1}{\tiny 0}
    \multiput(80,42)(30,0){1}{\tiny 1}
 
   \multiput(95,72)(30,0){1}{\tiny 0}
    \multiput(120,72)(30,0){1}{\tiny 0}
  \multiput(145,72)(30,0){1}{\tiny 0}
    \multiput(170,72)(30,0){1}{\tiny 0}
      \multiput(95,62)(30,0){1}{\tiny 0}
    \multiput(120,62)(30,0){1}{\tiny 0}
  \multiput(145,62)(30,0){1}{\tiny 0}
    \multiput(170,62)(30,0){1}{\tiny 0}
 \multiput(95,52)(30,0){1}{\tiny 0}
    \multiput(120,52)(30,0){1}{\tiny 0}
  \multiput(145,52)(30,0){1}{\tiny 0}
    \multiput(170,52)(30,0){1}{\tiny 0}
       \multiput(95,42)(30,0){1}{\tiny 0}
    \multiput(120,42)(30,0){1}{\tiny 0}
  \multiput(145,42)(30,0){1}{\tiny 0}
    \multiput(170,42)(30,0){1}{\tiny 0}
   
    \multiput(185,72)(30,0){1}{\tiny 1}
    \multiput(210,72)(30,0){1}{\tiny 0}
  \multiput(235,72)(30,0){1}{\tiny 0}
    \multiput(260,72)(30,0){1}{\tiny 0}
      \multiput(185,62)(30,0){1}{\tiny 0}
    \multiput(210,62)(30,0){1}{\tiny 1}
  \multiput(235,62)(30,0){1}{\tiny 0}
    \multiput(260,62)(30,0){1}{\tiny 0}
 \multiput(185,52)(30,0){1}{\tiny 0}
    \multiput(210,52)(30,0){1}{\tiny 0}
  \multiput(235,52)(30,0){1}{\tiny 1}
    \multiput(260,52)(30,0){1}{\tiny 0}
       \multiput(185,42)(30,0){1}{\tiny 0}
    \multiput(210,42)(30,0){1}{\tiny 0}
  \multiput(235,42)(30,0){1}{\tiny 0}
    \multiput(260,42)(30,0){1}{\tiny 1}
     \multiput(275,72)(30,0){1}{\tiny 0}
    \multiput(300,72)(30,0){1}{\tiny 0}
  \multiput(325,72)(30,0){1}{\tiny 0}
       \multiput(275,62)(30,0){1}{\tiny 1}
    \multiput(300,62)(30,0){1}{\tiny 0}
  \multiput(325,62)(30,0){1}{\tiny 0}
 
 \multiput(275,52)(30,0){1}{\tiny 0}
    \multiput(300,52)(30,0){1}{\tiny 1}
  \multiput(325,52)(30,0){1}{\tiny 0}
    
       \multiput(275,42)(30,0){1}{\tiny 0}
    \multiput(300,42)(30,0){1}{\tiny 0}
  \multiput(325,42)(30,0){1}{\tiny 1}

    \multiput(5,2)(30,0){1}{\tiny 0}
    \multiput(30,2)(30,0){1}{\tiny 0}
  \multiput(55,2)(30,0){1}{\tiny 0}
    \multiput(80,2)(30,0){1}{\tiny 0}
      \multiput(5,12)(30,0){1}{\tiny 0}
    \multiput(30,12)(30,0){1}{\tiny 0}
  \multiput(55,12)(30,0){1}{\tiny 0}
    \multiput(80,12)(30,0){1}{\tiny 0}
 \multiput(5,22)(30,0){1}{\tiny 0}
    \multiput(30,22)(30,0){1}{\tiny 0}
  \multiput(55,22)(30,0){1}{\tiny 0}
    \multiput(80,22)(30,0){1}{\tiny 0}
       \multiput(5,32)(30,0){1}{\tiny 0}
    \multiput(30,32)(30,0){1}{\tiny 0}
  \multiput(55,32)(30,0){1}{\tiny 0}
    \multiput(80,32)(30,0){1}{\tiny 0}
   
       \multiput(95,2)(30,0){1}{\tiny 0}
    \multiput(120,2)(30,0){1}{\tiny 0}
  \multiput(145,2)(30,0){1}{\tiny 0}
    \multiput(170,2)(30,0){1}{\tiny 1}
      \multiput(95,12)(30,0){1}{\tiny 0}
    \multiput(120,12)(30,0){1}{\tiny 0}
  \multiput(145,12)(30,0){1}{\tiny 1}
    \multiput(170,12)(30,0){1}{\tiny 0}
 \multiput(95,22)(30,0){1}{\tiny 0}
    \multiput(120,22)(30,0){1}{\tiny 1}
  \multiput(145,22)(30,0){1}{\tiny 0}
    \multiput(170,22)(30,0){1}{\tiny 0}
       \multiput(95,32)(30,0){1}{\tiny 1}
    \multiput(120,32)(30,0){1}{\tiny 0}
  \multiput(145,32)(30,0){1}{\tiny 0}
    \multiput(170,32)(30,0){1}{\tiny 0}

      \multiput(185,2)(30,0){1}{\tiny 0}
    \multiput(210,2)(30,0){1}{\tiny 0}
  \multiput(235,2)(30,0){1}{\tiny 0}
    \multiput(260,2)(30,0){1}{\tiny 1}
      \multiput(185,12)(30,0){1}{\tiny 0}
    \multiput(210,12)(30,0){1}{\tiny 0}
  \multiput(235,12)(30,0){1}{\tiny 1}
    \multiput(260,12)(30,0){1}{\tiny 0}
 \multiput(185,22)(30,0){1}{\tiny 0}
    \multiput(210,22)(30,0){1}{\tiny 1}
  \multiput(235,22)(30,0){1}{\tiny 0}
    \multiput(260,22)(30,0){1}{\tiny 0}
       \multiput(185,32)(30,0){1}{\tiny 1}
    \multiput(210,32)(30,0){1}{\tiny 0}
  \multiput(235,32)(30,0){1}{\tiny 0}
    \multiput(260,32)(30,0){1}{\tiny 0}

     \multiput(275,2)(30,0){1}{\tiny 0}
    \multiput(300,2)(30,0){1}{\tiny 0}
  \multiput(325,2)(30,0){1}{\tiny 0}
  
      \multiput(275,12)(30,0){1}{\tiny 0}
    \multiput(300,12)(30,0){1}{\tiny 0}
  \multiput(325,12)(30,0){1}{\tiny 1}
   \multiput(275,22)(30,0){1}{\tiny 0}
    \multiput(300,22)(30,0){1}{\tiny 1}
  \multiput(325,22)(30,0){1}{\tiny 0}
    
       \multiput(275,32)(30,0){1}{\tiny 1}
    \multiput(300,32)(30,0){1}{\tiny 0}
  \multiput(325,32)(30,0){1}{\tiny 0}

  \multiput(7,72)(0,90){1}{\vector(0,-1){50}}    \multiput(7,22)(0,90){1}{\vector(1,0){115}} 
  \multiput(122,22)(0,90){1}{\vector(0,1){20}}  
 
     \multiput(122,42)(0,90){1}{\vector(1,0){138}}  
   
 \color{red}  
   
    \multiput(262,42)(0,90){1}{\vector(0,-1){30}} 
     \multiput(262,12)(0,90){1}{\vector(1,0){65}} 
      \multiput(326,12)(0,90){1}{\vector(0,1){61}} 
       \multiput(326,72)(90,0){1}{\vector(-1,0){140}} 
       \multiput(186,72)(90,0){1}{\vector(0,-1){60}} 
        \multiput(186,12)(90,0){1}{\vector(-1,0){40}} 
         \multiput(146,12)(90,0){1}{\vector(0,1){60}} 
          \multiput(146,72)(90,0){1}{\vector(-1,0){140}}

\end{picture}

\par\bigskip
   
  \subsection{Sequence A100705 Via Brauer Configuration Algebras}

Partition trees $T_{P_{j}}$ associated to numbers $h_{j}=j^{3}+(j+1)^2$ in the proof of Theorem \ref{projective quadruple} define a sequence of Brauer configuration algebras $\Lambda_{E^{j}}$, $j\geq2$ induced by Brauer configurations $E^{j}$ whose set of vertices are 4-vertex paths contained in such trees. In order to describe Brauer configurations $E^{j}$, we assume the notation $L^{i_{1}}_{i_{2}, i_{3}+1}$ for the $i_{2}$th,  4-vertex path occurring in a third ramification of size $i_{1}$ in the partition tree $T_{P_{i_{3}}}$. In fact, vertices in $E^{j}$ is a labeling of 4-vertex paths in partition trees $T_{P_{i}}$, $1\leq i\leq j$. As an example, the five 4-vertex paths of $T_{P_{1}}$ belong to $P^{1}=P^{1}_{1}\cup P^{1}_{2}\cup P^{2}_{2}$, where $P^{1}_{1}=\{L^{1}_{1,1}\}$, $P^{2}_{2}=\{L^{2}_{1,2}, L^{2}_{2,2}, L^{2}_{3,2}, L^{2}_{4,2}\}$, and $P^{1}_{2}=\varnothing$. The seventeen 4-vertex paths $P^{2}$ of $T_{P_{2}}=T_{c_{4}, 4}$ are given by the following identities:

\begin{equation}
\begin{split}
P^{2}&=P^{2}_{2}\cup P^{2}_{3}\cup P^{3}_{3},\\
P^{2}_{3}&=\{L^{2}_{5,3}, L^{2}_{6,3}, L^{2}_{7,3}, L^{2}_{8,3}\},\\
P^{3}_{3}&=\{L^{3}_{1,3}, L^{3}_{2,3}, L^{3}_{3,3}, L^{3}_{4,3}, L^{3}_{5,3}, L^{3}_{6,3}, L^{3}_{7,3}, L^{3}_{8,3}, L^{3}_{9,3}\}.
\end{split}
\end{equation}

For $j\geq2$ fixed, $1\leq h_{i}\leq i^{2}$, $i^{2}+1\leq h'_{i}\leq i^{3}$, and $1\leq h_{i+1}\leq (i+1)^{2}$. The Brauer configuration $E^{j}=(E^{j}_{0},E^{j}_{1},\mu^{j},\mathcal{O}^{j})$ is defined in the following fashion:

\begin{enumerate}
\item 
\begin{equation}
\begin{split}
E^{j}_{0}&=\{L^{i}_{h_{i},i}, L^{i}_{h'_{i},i+1}, L^{i+1}_{h_{i+1}, i+1}\mid 2\leq i\leq j\}\cup\{L^{1}_{1, 1}\},\\
E^{j}_{1}&=\{T_{P_{i}}\mid 1\leq i\leq j\},\\
T_{P_{i}}&=P^{i}=P^{i}_{i}\cup P^{i}_{i+1}\cup P^{i+1}_{i+1},\\
P^{i}_{i}&=\{L^{i}_{1, i}, L^{i}_{2, i},\dots, L^{i}_{i^2, i} \},\quad{for}\hspace{0.1cm}i\geq2,\\
P^{i}_{i+1}&=\{L^{i}_{i^2+1, i+1},L^{i}_{i^2+2,i+1}, \dots, L^{i}_{i^3,i+1} \},\quad{for}\hspace{0.1cm}i\geq2,\\
P^{i+1}_{i+1}&=\{L^{i+1}_{1, i+1}, L^{i+1}_{2, i+1},\dots, L^{i+1}_{(i+1)^2, i+1} \},\quad{for}\hspace{0.1cm}i\geq2.
\end{split}
\end{equation}

\item The orientation $\mathcal{O}^{j}$ for successor sequences is defined by the usual order of natural numbers, i.e., any successor sequence has the form $T^{(s_{j_{1}})}_{P_{1}}<T^{(s_{j_{2}})}_{P_{2}}<\dots<T^{(s_{j_{(j-1)}})}_{P_{j-1}}<T^{(s_{j_{j}})}_{P_{j}}$, for some nonnegative integers $s_{j_{m}}$ ($s_{j_{m}}=0$, means that the vertex does not occur in the polygon $T_{P_{m}}$). 

\item $\mu^{j}$ is a multiplicity function such that, $\mu^{j}(L)=2$ for any $L\in\{L^{1}_{1,1}\}\cup P^{i}_{i+1}\cup P^{j+1}_{j+1}$, $2\leq i\leq j$, $\mu^{j}(L)=1$, otherwise. This multiplicity function is defined to avoid the presence of truncated vertices in the configuration.

\end{enumerate}

The Brauer quiver $Q_{E^{j}}$ has the following shape, where notation $c_{L^{i}_{h_{r},s}}$ means that the corresponding polygon has associated loops of type $c_{L^{i}_{h_{r},i+1}}$, $c_{L^{j+1}_{h_{m},j+1}}$,  $i^{2}+1\leq h_{r}\leq i^{3}$, $1\leq h_{m}\leq (i+1)^{2}$ defined by vertices $L^{i}_{h_{r},s}$ ($h_{1}=1$). In this case, $\alpha_{L^{i}_{i_{r},s}}$ ($\beta_{L^{i}_{i_{r},s}}$) denotes arrows determined by polygons $T_{P_{(i-1)}}$ and $T_{P_{i}}$, $1\leq i_{r}\leq  i^{2}$, $i\geq 2$. Loops $c_{L^{6}_{h_{6},6}}$ in the diagram appear if $j=5$. If $j>5$ then at the vertex $T_{P_{5}}$ there are associated only loops of type $c_{L^{5}_{h_{5},5}}$ and in the vertex $T_{P_{j}}$ there are attached loops of the form $c_{L^{j}_{h_{j},j}}$ and $c_{L^{j+1}_{h_{j+1}, j+1}}$ .
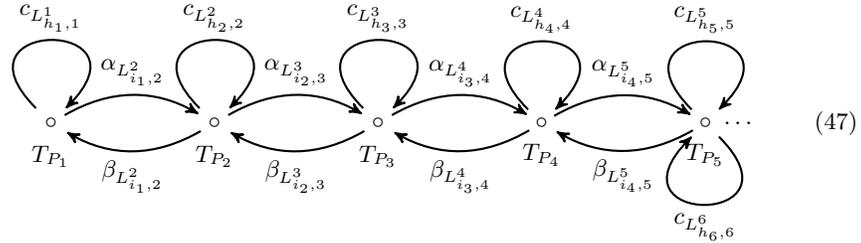
\begin{center}
\begin{tikzpicture}[->,>=stealth',shorten >=1pt,thick,scale=0.87]





\node[inner sep=2.5pt,label=below:$T_{P_{1}}$] at (1,3)(1){$\circ$};
\node[inner sep=2.5pt,label=below:$T_{P_{2}}$] at (3.5,3)(2){$\circ$};
\node[inner sep=2.5pt,label=below:$T_{P_{3}}$] at (6,3)(3){$\circ$};
\node[inner sep=2pt,label=below:$T_{P_{4}}$] at (8.5,3)(4){$\circ$};
\node[inner sep=2pt,label=below:$T_{P_{5}}$] at (11,3)(5){$\circ$};
\node[inner sep=2pt] at (11.5,3)(7){$\ldots$};


\path[->](1) [bend left] edge node[above] {$\alpha_{L^{2}_{i_{1},2}}$} (2);
\path[->](2) [bend left] edge node[above] {$\alpha_{L^{3}_{i_{2},3}}$}(3);
\path[->](3) [bend left] edge node[above] {$\alpha_{L^{4}_{i_{3},4}}$}(4);
\path[->](4) [bend left] edge node[above] {$\alpha_{L^{5}_{i_{4},5}}$}(5);


\path[->](2) [bend left] edge node[below] {$\beta_{L^{2}_{i_{1},2}}$} (1);
\path[->](3) [bend left] edge node[below] {$\beta_{L^{3}_{i_{2},3}}$} (2);
\path[->](4) [bend left] edge node[below] {$\beta_{L^{4}_{i_{3},4}}$} (3);
\path[->](5) [bend left] edge node[below] {$\beta_{L^{5}_{i_{4},5}}$} (4);

\Loop[dist=2cm,dir=SO,label= {$c_{L^{6}_{h_{6},6}}$},labelstyle=below](5)

\Loop[dist=2cm,dir=NO,label= {$c_{L^{2}_{h_{2},2}}$},labelstyle=above](2)
\Loop[dist=2cm,dir=NO,label= {$c_{L^{1}_{h_{1},1}}$},labelstyle=above](1)
\Loop[dist=2cm,dir=NO,label= {$c_{L^{3}_{h_{3},3}}$},labelstyle=above](3)
\Loop[dist=2cm,dir=NO,label= {$c_{L^{4}_{h_{4},4}}$},labelstyle=above](4)
\Loop[dist=2cm,dir=NO,label= {$c_{L^{5}_{h_{5},5}}$},labelstyle=above](5)
\node at (13,3)(1){$(47)$};
\end{tikzpicture}
\end{center}

The ideal $J$ in this case is generated by the following set of relations defined for all the possible values of $h, i, j, m, t, t'$ and $u$:

\begin{enumerate}
\item $(c_{L^{i}_{j,m}})^{2}$,
\item $c_{L^{j}_{h_{j},j}}c_{L^{j+1}_{h_{(j+1)}, j+1}}$,
\item $c_{L^{i}_{h_{j},m}}\alpha_{L^{i+1}_{i_{j},m+1}}$, 
\item $\alpha_{L^{i}_{i_{j},m}}c_{L^{i}_{h_{(j+1)},m}}$, 
\item $\alpha_{L^{j}_{i_{(j-1)},j}}c_{L^{j+1}_{h_{(j+1)},j+1}}$, 
\item $c_{L^{i}_{h_{j},m}}\beta_{L^{i}_{i_{(j-1)},m}}$, 
\item $c_{L^{j+1}_{h_{(j+1)},j+1}}\beta_{L^{j}_{i_{(j-1)},j}}$, 
\item $\beta_{L^{i}_{i_{j}, m}}c_{L^{i-1}_{h_{j}, m-1}}$,
\item $\alpha_{L^{i}_{i_{j},m}}\alpha_{L^{i+1}_{i_{(j+1)},m+1}}$, 
\item $\beta_{L^{i}_{i_{j},m}}\beta_{L^{i-1}_{i_{(j-1)},m-1}}$,
\item $s^{2}_{L^{i}_{t, u}}a$, where $a$ is the first arrow of the special cycle $s_{L^{i}_{t, u}}$ associated to the vertex $L^{i}_{t, u}$,
\item $s^{2}_{L^{i}_{t, i}}-s^{2}_{L^{i}_{t', i}}$,\quad $s^{2}_{L^{i}_{t,i}}-s^{2}_{L^{i}_{h,i+1}}$,\quad $s^{2}_{L^{i}_{t, i}}-s^{2}_{L^{i+1}_{m,i+1}}$.
\end{enumerate}
If we let $\Lambda_{E^{j}}$ denote the algebra $kQ_{E^{j}}/J$, then the following result categorifies sequence $h_{n}=n^{3}+(n+1)^{2}$, $n\geq1$, by considering its elements as invariants of objects of the category $\mathrm{mod}\hspace{0.1cm}\Lambda_{E^{j}}$.

 \begin{teor}\label{tetradnew}
\textit{For $j\geq2$ fixed and $1\leq i\leq j$, the number of summands in the heart of the indecomposable projective module $T_{P_{i}}$ over the Brauer configuration algebra $\Lambda_{E^{j}}$ equals the number of cycles associated to a preprojective representation of the tetrad of type $\mathrm{IV}$ and order $i+1$.}

\end{teor}
\textbf{Proof.} Since for any indecomposable projective module $T_{P_{i}}$, it holds that $\mathrm{rad}^2\hspace{0.1cm} T_{P_{i}}\neq0$ then the theorem follows from Theorem \ref{multiserial} and the definition of the polygon $T_{P_{i}}$ which has $i^{3}+(i+1)^{2}$ nontruncated vertices counting repetitions.\hspace{0.5cm}$\square$
\par\bigskip
Numbers in the sequence A100705 appear by computing the dimension of quotient spaces of the form $\Lambda_{E^{i+1}}/F_{i+1}$, where $F_{i+1}$ is a $k$-subspace of $\Lambda_{E^{i+1}}$ isomorphic to $\Lambda_{E^{i}}$. Actually the following results hold.

\begin{teor}\label{tetradnew1}
\textit{For $j\geq 2$ fixed, it holds that $\mathrm{dim}_{k}\hspace{0.1cm}\Lambda_{E^{j}}=\frac{(j+1)(j+2)(2j+3)}{6}+(\frac{j(j+1)}{2})^2+2j-1$. And $\mathrm{dim}_{k}\hspace{0.1cm}\Lambda_{E^{n+1}}-\mathrm{dim}_{k}\hspace{0.1cm}\Lambda_{E^{n}}=h_{n+1}+2$, where $h_{n}=n^{3}+(n+1)^{2}$, $n\geq2$.}

\end{teor}
\textbf{Proof.} Note that, $|E^{j}_{1}|=j$, and $E^{j}_{0}=H_{1}\cup H_{2}$, where $H_{1}=\{L^{i}_{h^{1}_{i, i}}\mid 2\leq i\leq j+1\}$ and $H_{2}=\{L^{m}_{h^{2}_{m, m+1}}\mid 1\leq m\leq j\}$, $L^{1}_{h^{2}_{1,2}}=L^{1}_{1,1}$, $1\leq h^{1}_{i}\leq i^{2}$, $1\leq h^{2}_{m}\leq m^{3}-m^{2}$ (terms $L^{m}_{h^{2}_{s, m+1}}$ correspond bijectively to vertices $L^{m}_{m^{2}+s, m+1}\in P^{m}_{m+1}$). If $x\in H_{1}$ then $val(x)=2$ and $\mu^{j}(x)=1$, whereas, if $y\in H_{2}$ then $val(y)=1$ and $\mu^{j}(y)=2$. The result follows bearing in mind that $val(L^{1}_{h^{2}_{1,2}})=1$, $\mu^{j}(L^{1}_{h^{2}_{1,2}})=2$ and that for $i$ and $m$ fixed, there are $i^{2}$ vertices of type $H_{1}$, $1\leq i\leq j+1$ and $m^{3}-m^{2}$ vertices of type $H_{2}$, $2\leq m\leq j$. We are done.\hspace{0.5cm}$\square$\par\bigskip

\begin{teor}\label{tetradnew2}
For $j\geq 3$, $\mathrm{dim}_{k}\hspace{0.1cm}Z(\Lambda_{E^{j}})-\mathrm{dim}_{k}\hspace{0.1cm}Z(\Lambda_{E^{(j-1)}})=j^{3}-j^{2}+1$.

\end{teor}

\textbf{Proof.} Since $\mathrm{rad}^{2}\hspace{0.1cm}\Lambda_{E^{j}}\neq0$, $|E^{j}_{0}|=\underset{i=1}{\overset{j}{\sum}}i^{3}+(j+1)^{2}$ and the number of loops in the quiver $Q_{E^{j}}$ equals $|\mathscr{C}_{E^{j}}|$ then $\mathrm{dim}_{k}\hspace{0.1cm}  Z(\Lambda_{E^{j}})=\underset{i=1}{\overset{j}{\sum}}(i^{3}-i^{2})+(j+2)$. We are done.\hspace{0.1cm}$\square$

\par\bigskip
The sequence $j^{3}-j^2+1$ appears encoded in the OEIS as A100104.

\section{Brauer Configuration Algebras Via Integer Sequences}

In this section, we describe the way that some Brauer configuration algebras can be defined by using integer sequences. 

\subsection{Brauer Configuration Algebras Arising From the Sequence A100705}

We consider Brauer configuration algebras of the form $\Lambda_{\Gamma_n}=kQ_{\Gamma_{n}}/{J}$ induced by the Brauer configuration $\Gamma_{n}$ whose polygons are defined by integer partitions of the elements in the sequence A100705 (see \ref{equation}). And such that
For $n\geq 2$ fixed, $\Gamma_{n}=(\Gamma_{0},\Gamma_{1}, \mu, \mathcal{O})$ with
\begin{enumerate}
\addtocounter{equation}{1}
\item \begin{equation}\label{sequence}
\begin{split}
\Gamma_{0}&=\{1,2,3,\dots,n,n+1\},\\
\Gamma_{1}&=\{V_{t}=t^{(t^2)}(t+1)^{(t+1)}\}_{1\leq t\leq n},\hspace{0.1cm}i.e., \mathrm{occ}(t,V_{t})=t^2, \mathrm{occ}(t+1,V_{t})=t+1.
\end{split}
\end{equation}
\item The orientation $\mathcal{O}$ is defined in such a way that for $2\leq i\leq n$ \textit{at vertex} $i$, $V^{(i,<)}_{i-1}< V^{(i^{2},<)}_{i}$, where $V^{(y,<)}_{x}$ means that the polygon $V_{x}$ occurs $y$ times in the successor sequence of the corresponding vertex, in particular, $V_{i-1}<V_{i}$. At the vertex $n+1$, the successor sequence has the form $V^{(n+1,<)}_{n}$, in this case, $V_{n,1}<V_{n,2}<\dots<V _{n, n}< V_{n, n+1}$ where $V_{n, i}$ denotes the $i$th occurrence of the polygon $V_{n}$ in the sequence. \par\bigskip

\item The multiplicity function $\mu$ is such that $\mu(j)=1$, for any $j\in \Gamma_{0}$.
\end{enumerate}

 The following is the quiver $Q_{\Gamma_{n}}$ associated to the Brauer configuration $\Gamma_{n}$, worth noting that there is no arrow connecting vertex 1 with any other vertex provided that it is truncated (see Theorem \ref{multiserial}, item 5), besides we use the symbol $[x_j; y_j]$ to denote that the vertex $x_j$ occurs $y_j$ times at the polygon $h_{j}=j^{3}+(j+1)^2$ (see identity (\ref{word})). And $c^{i}_{j}$ is a set of loops $\{c^{i}_{j_y}\mid 1\leq y\leq \mathrm{occ}(x_j, h_{j})-1, 2\leq i\leq n+1\}$. For instance, at $17$ there are associated the loops, $c^{2}_{17_1}, c^{2}_{17_2},c^{2}_{17_3}$ and $c^{3}_{17_1}, c^{3}_{17_2}$.

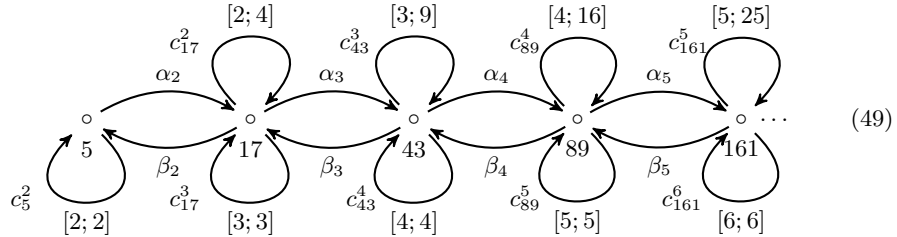
\begin{center}
\begin{tikzpicture}[->,>=stealth',shorten >=1pt,thick,scale=0.87]


\node at (2.5,4.2)(0){$c_{17}^{2}$};
\node at (5.1,4.2)(0){$c_{43}^3$};
\node at (7.7,4.2)(0){$c_{89}^4$};
\node at (10.2,4.2)(0){$c_{161}^5$};

\node at (0.0,1.8)(0){$c_{5}^{2}$};
\node at (2.5,1.8)(0){$c_{17}^{3}$};
\node at (5.2,1.8)(0){$c_{43}^4$};
\node at (7.7,1.8)(0){$c_{89}^5$};
\node at (10.1,1.8)(0){$c_{161}^6$};


\node[inner sep=2.5pt,label=below:$5$] at (1,3)(1){$\circ$};
\node[inner sep=2.5pt,label=below:$17$] at (3.5,3)(2){$\circ$};
\node[inner sep=2.5pt,label=below:$43$] at (6,3)(3){$\circ$};
\node[inner sep=2pt,label=below:$89$] at (8.5,3)(4){$\circ$};
\node[inner sep=2pt,label=below:$161$] at (11,3)(5){$\circ$};
\node[inner sep=2pt] at (11.5,3)(7){$\ldots$};


\path[->](1) [bend left] edge node[above] {$\alpha_2$} (2);
\path[->](2) [bend left] edge node[above] {$\alpha_3$}(3);
\path[->](3) [bend left] edge node[above] {$\alpha_4$}(4);
\path[->](4) [bend left] edge node[above] {$\alpha_5$}(5);


\path[->](2) [bend left] edge node[below] {$\beta_2$} (1);
\path[->](3) [bend left] edge node[below] {$\beta_3$} (2);
\path[->](4) [bend left] edge node[below] {$\beta_4$} (3);
\path[->](5) [bend left] edge node[below] {$\beta_5$} (4);

\Loop[dist=2cm,dir=SO,label= {$[2;2]$},labelstyle=below](1)
\Loop[dist=2cm,dir=SO,label= {$[3;3]$},labelstyle=below](2)
\Loop[dist=2cm,dir=SO,label= {$[4;4]$},labelstyle=below](3)
\Loop[dist=2cm,dir=SO,label= {$[5;5]$},labelstyle=below](4)
\Loop[dist=2cm,dir=SO,label= {$[6;6]$},labelstyle=below](5)

\Loop[dist=2cm,dir=NO,label= {$[2;4]$},labelstyle=above](2)
\Loop[dist=2cm,dir=NO,label= {$[3;9]$},labelstyle=above](3)
\Loop[dist=2cm,dir=NO,label= {$[4;16]$},labelstyle=above](4)
\Loop[dist=2cm,dir=NO,label= {$[5;25]$},labelstyle=above](5)
\node at (13,3)(1){$(49)$};
\end{tikzpicture}
\end{center}

The following are examples of polygons in a Brauer configuration $\Gamma_n$:
\addtocounter{equation}{1}
\begin{equation}
\begin{split}
{5}&={(1)+(2+2)}=(1)^{(1)}(2)^{(2)},\\
{17}&={(2+3)+(2+3)+(2+2+3)}=(2)^{(4)}(3)^{(3)},\\
43&=(3+3+4)+(3+3+4)+(3+3+4)+(3+3+3+4)=(3)^{(9)}4^{(4)},\\
{89}&={(4+4+4+5)+(4+4+4+5)+(4+4+4+5)+(4+4+4+5)+(4+4+4+4+5)}\\
\vdots&=\vdots
\end{split}
\end{equation}

The ideal ${J}$ is generated by the following relations where for a fixed $2\leq l\leq n+1$, $P^{i,l}_{h_{j}}$ is the product of $i$ loops of type $l$ ($1\leq i\leq \mathrm{occ}(l, h_j)-1$) attached to the polygon $h_{j}$ with $y_{j}-1$ being the total number of such loops ($y_{j}\in \{j^{2},j\}$):

\begin{enumerate}
\item $c^{u}_{j_x}c^{v}_{j_y}$, if $u\neq v$, for all the possible values of $u$, $v$, $x$, $y$ and $j$, 
\item $c^{t}_{j_x}c^{t}_{j_y}=c^{t}_{j_y}c^{t}_{j_x}$, for all the possible values of $x$, $y$, $t$, and $j$,
\item $(c^{t}_{j_x})^2$ for all the possible values of $j$, $t$ and $x$,
\item$c^{h}_{j_{x}}\alpha_{h+1}$;\quad $\alpha_{h}c^{h+1}_{(j+1)_{x}}$;\quad$c^{h}_{j_{x}}\beta_{h-1}$;\quad $\beta_{h}c^{h-1}_{(j-1)_{x}}$,\quad$\alpha_{j}\beta_{j}$ for all the possible values of $h, j$ and $x$,
\item $\alpha_{i}\alpha_{i+1}$;\quad $\beta_{j+1}\beta_{j}$, $2\leq i\leq n-1$, $2\leq j\leq n-1$,
\item If 
\begin{equation}
\begin{split}
\varepsilon^{1}_{j}&=P^{u, j}_{h_{j}}\alpha_{j}P^{y_{j+1}-1,j}_{h_{j+1}}\beta_{j}P^{y_{j}-(1+u),j}_{h_{j}},\\
\varepsilon^{2}_{j}&=\alpha_{j}P^{y_{j+1}-1,j}_{h_{j+1}}\beta_{j}P^{y_{j}-1,j}_{h_{j}},\\
\varepsilon^{3}_{j}&=P^{u, j}_{h_{j+1}}\beta_{j}P^{y_{j}-1,j}_{h_{j}}\alpha_{j}P^{y_{j+1}-(1+u),j}_{h_{j+1}},\\
\varepsilon^{4}_{j}&=\beta_{j}P^{y_{j}-1,j}_{h_{j}}\alpha_{j}P^{y_{j+1}-1,j}_{h_{j+1}},\\
\varepsilon^{5}_{j+1}&=P^{v ,j+1}_{h_{j+1}}\alpha_{j+1}P^{y_{j+2}-1,j+1}_{h_{j+2}}\beta_{j+1}P^{y_{j+1}-(1+v), j+1}_{h_{j+1}},\\
\varepsilon^{6}_{j+1}&=\alpha_{j+1}P^{y_{j+2}-1,j+1}_{h_{j+2}}\beta_{j+1}P^{y_{j+1}-1, j+1}_{h_{j+1}},\\
\varepsilon^{7}_{j+1}&=P^{v, +1}_{h_{j+2}}\beta_{j+1}P^{y_{j+1}-1, j+1}_{h_{j+1}}\alpha_{j+1}P^{y_{j+2}-(1+v), j+1}_{h_{j+2}},\\
\varepsilon^{8}_{j+1}&=\beta_{j+1}P^{y_{j+1}-1, j+1}_{h_{j+1}}\alpha_{j+1}P^{y_{j+2}-1, j+1}_{h_{j+2}},
\end{split}
\end{equation}
then there are relations of the form $\varepsilon^{r}_{s}-\varepsilon^{r'}_{s'}$ where $r, r'\in\{1,\dots, 8\}$, $r\neq r'$ and $s, s'\in\{j, j+1\}$, for all the possible values of $u$, $v$ and $j$,
\item $\varepsilon^{2}_{j}\alpha_{j}$,\quad $\varepsilon^{4}_{j}\beta_{j}$,\quad $\varepsilon^{6}_{j+1}\alpha_{j+1}$,\quad $\varepsilon^{8}_{j}\beta_{j+1}$.
\end{enumerate}

\par\bigskip
 The following results are consequences of Theorems \ref{multiserial}, \ref{Serra}, and Proposition \ref{dimension}. 
   \addtocounter{corol}{4}

  \begin{corol}\label{tetrad2}
\textit{For $n\geq2$ fixed and $2\leq i\leq n$, the number of summands in the heart of the indecomposable projective module $V_{i}$ over the algebra $\Lambda_{\Gamma_{n}}$ is $i^{2}+i+1$.}

\end{corol}
\textbf{Proof.} Since for any indecomposable projective module $V_{i}$, it holds that $\mathrm{rad}^2\hspace{0.1cm} V_{i}\neq0$ then the theorem follows from Theorem \ref{multiserial} and the definition of the polygon $V_{i}$ which has $i^{2}+i+1$ nontruncated vertices counting repetitions.\hspace{0.5cm}$\square$.

\begin{corol}\label{tetrad3}
\textit{For $n\geq2$ fixed, $\mathrm{dim}_{k}\hspace{0.1cm}\Lambda_{\Gamma_{n}}=\underset{m=2}{\overset{n}{\sum}}(m(m+1))^{2}-\frac{1}{3}(n-3)(n+1)(n+2)$. And $\mathrm{dim}_{k}\hspace{0.1cm}\Lambda_{\Gamma_{n+1} }/G_{n+1}=2(1-t_{n})+[(n+1)(n+2)]^2$, where for $i\geq1$, $t_{i}$ denotes the $i$th triangular number. And $G_{n+1}$ is a $k$-subspace of $\Lambda_{\Gamma_{n+1}}$ isomorphic to $\Lambda_{\Gamma_{n}}$}\par\bigskip

\textbf{Proof.} It is enough to observe that for $n\geq2$ and $2\leq j<n+1$, it holds that $val(j)=j^{2}+j$, whereas $val(1)=1$ and $val (n+1)=n+1$.  The theorem holds as a consequence of Proposition \ref{dimension}.\hspace{0.5cm} $\square$

\end{corol}

\begin{corol}\label{tetrad4}
\textit{For $n\geq2$ fixed, it holds that $\mathrm{dim}_{k}\hspace{0.1cm}Z(\Lambda_{\Gamma_{n}})=\frac{n(n+1)(n+2)}{3}+1$. And $\mathrm{dim}_{k}\hspace{0.1cm}Z(\Lambda_{\Gamma_{n+1}})/Z_{n+1}=2t_{n+1}$, where $Z_{n+1}$ is a $k$-subspace of $Z(\Lambda_{\Gamma_{n+1}})$ isomorphic to $Z(\Lambda_{\Gamma_{n}})$. }

\end{corol}

\textbf{Proof.} Since $\mathrm{rad}^2\hspace{0.1cm}\Lambda_{\Gamma_n}\neq 0$, the result is a consequence of Theorem \ref{Serra} with $\mu(i)=1$, for any $2\leq i \leq n+1$, $|\Gamma_{0}|=n$,
$|\Gamma_{1}|=n$, $\mathrm{occ}(i, h_i)+\mathrm{occ}(i+1, h_i)=i^2+i+1$, $2\leq i\leq n$, and $\mathrm{occ}(2, h_1)=2$.\hspace{0.5cm}$\square$

\addtocounter{Nota}{9}
\begin{Nota}
  Note that Corollaries \ref{tetrad2}-\ref{tetrad4} are categorifications of the integer sequences $n^{2}+n+1$ (encoded in the OEIS as A002061), $\underset{m=2}{\overset{n}{\sum}}(m(m+1))^{2}-\frac{1}{3}(n-3)(n+1)(n+2)$, and $\frac{n(n+1)(n+2)}{3}+1$ (which is the sequence A064999). Elements of the sequence A064999 appear as coefficients (in the case $t=3$) of the generating polynomial of a $n$-twist knot with the form $P_{n}(x)=\underset{t\geq0}{\sum} a_{n, t}x^t$. And the sequence $\underset{t=2}{\overset{n}{\sum}}(t(t+1))^{2}=\underset{1\leq i<j\leq n}{\sum}(j-i)^{3}$ is encoded A024166 in the OEIS. 
  \end{Nota}

\subsection{Brauer Configuration Algebras Arising From Counting Functions}

In this section,  we consider Brauer configurations arising from counting functions which are strictly increasing integer sequence  whose elements count a given class of objects $\mathscr{D}_{n}$. For instance, $\mathscr{D}_{n}$ can be the set of linear extensions of a poset $(\mathscr{P}_{n},\unlhd)=\{(i,j)\in\mathbb{N}^{2}\mid 0\leq i\leq j\leq n\}$, where $\unlhd$ is a partial order defined on $\mathscr{P}_{n}$ such that $(i, j)\unlhd (i', j')$ if and only if $i\leq i'$ and $j\leq j'$. According to Stanley \cite{Stanley}, the number of linear extensions $e(\mathscr{P}_{n})$ of $\mathscr{P}_{n}$ is equal to the number of lattice paths from $(0,0)$ to $(n, n)$ with steps $(1,0)$ and $(0,1)$, which never rise above the main diagonal $x=y$ of the plane $(x, y)$-plane.  It can be shown that $e(\mathscr{P}_{n})$ is given by the $n$th Catalan number $C_{n}=\frac{1}{n+1}\binom{2n}{n}$.\par\bigskip

We define now a family of Brauer configuration algebras $\Lambda_{D^{n}}$, $n>1$ arisen from Brauer configurations $D^{n}$ whose nontruncated vertices are in correspondence with objects of type $\mathscr{D}_{n}$, polygons are obtained by choosing objects of type $\mathscr{D}_{s}$, for $1\leq s\leq n$. We assume the notation $L^{s}_{j, n}\in\mathscr{D}_{s}$  for the $j$th object of type $s$ in a given polygon. Without loss of generality, we assume that for the first polygon $P_{1}$, it holds that $|P_{1}|=u_{1}>1$.\par\bigskip

For $n\geq 2$ fixed, the definition of the Brauer configuration 
\begin{center}
$D^{n}=(D_{0}^{n}, D_{1}^{n},  \mu^{n}, \mathcal{O}^{n})$ 
\end{center}
goes as follows: 
\begin{equation}\label{universal}
\begin{split}
D_{0}^{n}&=\{L^{s}_{i_{s}, n}, 2\leq s\leq n, 1\leq i_{s}\leq u_{s}-u_{(s-1)}\}\cup P_{1},\\
P_{1}&=\{L^{1}_{i_{{1}}, n}\mid 1\leq i_{{1}}\leq u_{1}\},\\
D_{1}^{n}&=\{P_{h}\mid 1\leq h\leq n\},\\
P_{h}&= P^{h}_{(h-1)}\cup P^{h}_{h},\quad | P^{h}_{(h-1)}|= | P_{(h-1)}|,\quad 2\leq h\leq n,\\
P^{h}_{(h-1)}&=\{L^{s}_{i_{s}, n}\mid 1\leq s\leq h-1\},\\
P^{h}_{h}&=\{L^{h}_{i_{h}, n}\mid 1\leq i_{h}\leq u_{h}-u_{(h-1)}, 2\leq h\leq n\},\\
\mu^{n}(L)&=1, \hspace{0.1cm}\text{for any vertex}\hspace{0.1cm}L\in (D_{0})^{n}\backslash P^{n}_{n},\\
\mu^{n}(L)&=2, \hspace{0.1cm}\text{for any vertex}\hspace{0.1cm}L\in P^{n}_{n},
\end{split}
\end{equation}

In $P^{h}_{(h-1)}$, it holds that, $1\leq i_{1}\leq  u_{1}$ if $s=1$, and $1\leq i_{s}\leq u_{s}-u_{(s-1)}$, if $s>2$.\par\bigskip

The orientation $\mathcal{O}^{n}$ is defined by the usual order of natural numbers. Thus, for a vertex $L^{i}_{j, n}\in D_{0}^{n}\backslash P^{n}_{n}$, the successor sequence has the form\par\bigskip 

\begin{centering}
$P_{i}<P_{(i+1)}<\dots< P_{(n-1)}<P_{n}$\par\bigskip
\end{centering}
For vertices $L^{(n-1)}_{r, n}$, the successor sequence has the form $P_{(n-1)}<P_{n}$, whereas for vertices of the form $L^{n}_{r, n}\in P^{n}_{n}$, the orientation is of the form  $P_{n}< P_{n}$.
\par\bigskip

The following is the shape of the Brauer quiver $Q_{D^{n}}$ defined by $D^{n}$.

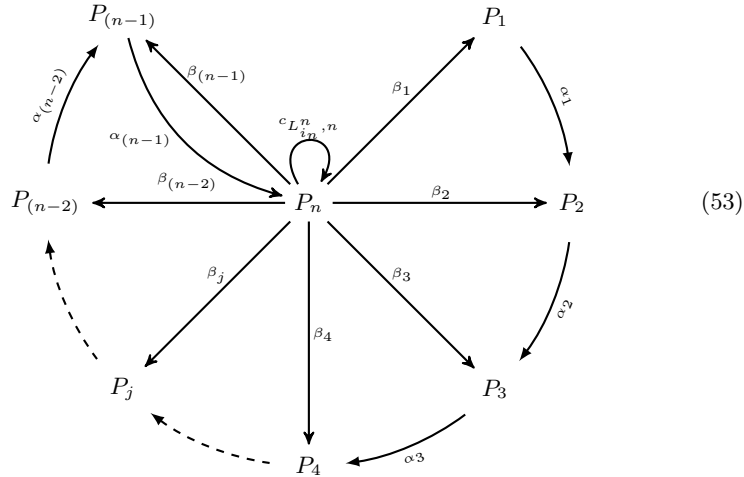
\begin{center}
	\begin{tikzpicture}[->,>=stealth',shorten >=1pt,thick,scale=0.5]
	\def \radius {7cm}
	\def \margin {8} 
	{
		\node at ({360/8 * (1 - 1)}:\radius)(U2) {{$P_2$}};
		
		\node at ({360/8 * (2 - 1)}:\radius)(U3) {{$P_1$}};
		
		\node at ({360/8 * (3 - 1)}:\radius) (U4){};
		
		\node at ({360/8 * (4 - 1)}:\radius)(U5) {{$P_{(n-1)}$}};

		\node at ({360/8 * (5 - 1)}:\radius)(U6) {{$P_{(n-2)}$}};
		
		\node at ({360/8 * (6 - 1)}:\radius)(U7) {{$P_j$}};
		
		\node at ({360/8 * (7 - 1)}:\radius)(Uk-1) {{$P_4$}};
		
		\node at ({360/8 * (8 - 1)}:\radius)(Uk) {{$P_3$}};

		\node at (0,0)(U1) {$P_{n}$};

		
		\draw[<-, >=latex,black] ({360/8* (1 - 1)+\margin}:\radius)
		arc ({360/8 * (1 - 1)+\margin}:{360/8 * (1)-\margin}:\radius)node[midway,sloped,above] {\tiny$\alpha_1$};
		
		\draw[<-, >=latex,black] ({360/8* (8 - 1)+\margin}:\radius)
		arc ({360/8 * (8 - 1)+\margin}:{360/8 * (8)-\margin}:\radius)node[midway,sloped,below] {\tiny$\alpha_{2}$};

		\draw[<-, >=latex,black] ({360/8* (7 - 1)+\margin}:\radius)
		arc ({360/8 * (7 - 1)+\margin}:{360/8 * (7)-\margin}:\radius)node[midway,sloped,below] {\tiny$\alpha_{3}$};

		\draw[dashed,<-, >=latex,black] ({360/8* (6 - 1)+\margin}:\radius)
		arc ({360/8 * (6 - 1)+\margin}:{360/8 * (6)-\margin}:\radius)node[midway,sloped,below] {};
		
		\draw[dashed,<-, >=latex,black] ({360/8* (5 - 1)+\margin}:\radius)
		arc ({360/8 * (5 - 1)+\margin}:{360/8 * (5)-\margin}:\radius)node[midway,sloped,below] {};

		\draw[<-, >=latex,black] ({360/8* (4 - 1)+\margin}:\radius)
		arc ({360/8 * (4 - 1)+\margin}:{360/8 * (4)-\margin}:\radius)node[midway,sloped,above] {\tiny$\alpha_{(n-2)}$};

		\path[->,black](U1) edge node[above=2,font=\tiny] {$\beta_1$} (U3);
		\path[->,black](U1) edge node[above=-2,font=\tiny] {$\beta_2$} (U2);
		
		\path[->,black](U1) edge node[above=2,font=\tiny] {$\beta_3$}  (Uk);
		
		\path[->,black](U1) edge node[right=-2,font=\tiny] {$\beta_{4}$} (Uk-1);
		\path[->,black](U1) edge node[above=1,font=\tiny] {$\beta_j$} (U7);
		
		\path[->,black](U1) edge node[above=1,font=\tiny] {$\beta_{(n-2)}$} (U6);

		\path[->,black](U5) [bend right] edge node[left=1,font=\tiny] {$\alpha_{(n-1)}$} (U1);
		\path[->,black](U1)  edge node[above=7,font=\tiny] {$\beta_{(n-1)}$} (U5);
	}
	
	
	\path[->,black] (U1) edge  [in=60,out=120,loop] node[above=-3,font=\tiny] {$c_{L^{n}_{i_{n}}, n}$}
	();

	\node at (-12,0)(U11) {\color{white}$(46)$};
	\node at (11,0)(U10) {$(53)$};
	\end{tikzpicture}
	\label{examplebca}
\end{center}

$\alpha_{j}$, $\beta_{j}$ and $c_{L^{n}_{i_{n}},n}$ denote  $j\times 1$-matrices of the form:
\par\bigskip
$\alpha_{j}=\left[\begin{array}{ccc}\alpha^{j}_{L^{1}_{i_{1}},n}  \\\alpha^{j}_{L^{2}_{i_{2}},n} \\\vdots  \\\alpha^{j}_{L^{s}_{{i_{s}}},n}  \\\vdots \\\alpha^{j}_{L^{j}_{i_{j}}, n}\end{array}\right]$, \hspace{0.2cm} $\beta_{j}=\left[\begin{array}{ccc}\beta^{j}_{L^{1}_{i_{1}},n}  \\\beta^{j}_{L^{2}_{i_{2}},n} \\\vdots  \\\beta^{j}_{L^{s}_{i_{s}},n}  \\\vdots \\\beta^{j}_{L^{j}_{i_{j}}, n}\end{array}\right]$, \hspace{0.2cm} $c_{L^{n}_{i_{n}},n}=\left[\begin{array}{ccc}c_{L^{n}_{{1}},n}  \\c_{L^{n}_{{2}},n} \\\vdots  \\c_{L^{n}_{{i}},n} \\\vdots \\c_{L^{n}_{(u_{n}-u_{(n-1)}),n}}\end{array}\right] $\par\bigskip

where $\alpha^{j}_{L^{s}_{i_{s}},n}$ ($\beta^{j}_{L^{s}_{i_{s}},n}$) is a set of arrows defined by the successor sequence at vertex $L^{s}_{i_{s}, n}$ connecting the corresponding polygons (polygon $P_{n}$ with the corresponding $P_{j}$). And $c_{L^{n}_{i_{n},n}}$ is a set of loops defined by vertices ${L^{n}_{i_{n},n}}$, $1\leq i_{n}\leq u_{n}-u_{(n-1)}$.

\par\bigskip

The following relations generate the admissible ideal $J$ of the Brauer configuration algebra $\Lambda_{D^{n}}=kQ_{D^{n}}/J$, for all possible values of $i, i', j, j', r, r'$ and $n$.

\begin{enumerate}
\item $(c_{L^{n}_{{r},n}})^{2}$,\quad $c_{L^{n}_{{r},n}}c_{L^{n}_{{r'},n}}$,\quad $r\neq r'$,
\item $\alpha^{j}_{L^{i}_{{r}},n} \alpha^{j}_{L^{i'}_{{r'}}, n} $,\quad$i\neq i'$,
\item $\alpha^{j}_{L^{i}_{{r}},n} \alpha^{(j+1)}_{L^{i'}_{{r'}}, n}$,
\item  $c_{L^{n}_{i_{n},n}}\beta_{i}$,
\item  $\alpha_{(n-1)}c_{L^{n}_{i_{n},n}}$,
\item For $1\leq j\leq n$, fixed and $1\leq i\leq j$, $s_{L^{i}_{j, n}}-s_{L^{j}_{j', n}}$ where $s_{x}$ is a special cycle associated to the vertex $x$,
\item $\beta_{i}\alpha_{i}$,
\end{enumerate}
products of the form; $\beta_{i}\alpha_{i}$, $\alpha_{(n-1)}c_{L^{n}_{i_{n},n}}$, $c_{L^{n}_{i_{n},n}}\beta_{i}$ means that relations of the form $xx'$, $y'y$ and $z'z$ have place where $x'$, $y'$ and $z'$ are entries of the corresponding matrices.\par\bigskip
The following result categorifies numbers of a counting function $u_{t}$, for $t\geq1$.
\addtocounter{teor}{4}

\begin{teor}\label{gral1}
\textit{For $1\leq i\leq n$ and $n>1$ fixed, the number of summands in the heart of the indecomposable projective module $P_{i}$ over the algebra $\Lambda_{D^{n}}$ is $u_{i}$.}

\end{teor}

\textbf{Proof.}  Since $\mathrm{rad}^{2}\hspace{0.1cm}P_{i}\neq0$, then the number of summands in the heart of the indecomposable projective module $P_{i}$ equals the number of its nontruncated vertices counting repetitions, which by definition is given by the sum $u_{1}+(u_{2}-u_{1})+\dots+(u_{(i-1)}-u_{(i-2)})+(u_{i}-u_{(i-1)})=u_{i}$. We are done.\hspace{0.5cm} $\square$

\begin{teor}\label{gral2}
\textit{For $n\geq 1$ fixed, $\mathrm{dim}_{k}\hspace{0.1cm}\Lambda_{D^{n}}=2n+n(n-1)u_{1}+2\underset{i=2}{\overset{n}{\sum}}t_{(n-i)}(u_{i}-u_{(i-1)})$.}

\end{teor}

\textbf{Proof.} It suffices to note that $val(L^{i}_{i_{s},n})=n-i+1$ and $\mu^{n}(L^{i}_{i_{s},n})=1$ for any $L^{i}_{i_{s},n}\in D_{0}^{n}\backslash P^{n}_{n}$, whereas for any $x\in P^{n}_{n}$, it holds that $val(x)=1$ and $\mu^{n}(x)=2$. We are done.\hspace{0.5cm}$\square$
\par\bigskip
Since for any $n>1$, $\mathrm{rad}^{2}\hspace{0.1cm}\Lambda_{D^{n}}\neq0$, then we have the following result regarding the center of these algebras.
\begin{teor}\label{gral3}
\textit{For $n\geq 2$ fixed, $\mathrm{dim}_{k}\hspace{0.1cm}Z(\Lambda_{D^{n}})=(u_{n}-u_{(n-1)})+(n+1)$.}
\end{teor}

\textbf{Proof.} Note that $|D_{0}^{n}|=u_{n}$,\quad $|D_{1}^{n}|=n$,\quad $\underset{\alpha\in D_{0}^{n} }{\sum}\mu^{n}(\alpha)=2u_{n}-u_{(n-1)}$. Since $\#(Loops\hspace{0.1cm} Q_{D^{(n-1)}})=|\mathscr{C}_{\mathscr{D}^{n}}|$, the theorem holds.\hspace{0.5cm}$\square$

\addtocounter{Nota}{3}

\begin{Nota}\label{general}
Perhaps, the sequence $C_{n}$ of Catalan numbers is one of the most interesting counting functions, they count the number of plane binary trees with $n+1$ endpoints (or $2n+1$ vertices), the number of triangulations of an $(n+3)$ polygon, or the number of paths $L$ in the $(x, y)$-plane from $(0,0)$ to $(2n, 0)$ with steps $(1, 1)$ and $(1, -1)$ that never pass below the $x$-axis, such paths are called \textit{Dyck paths} \cite{Stanley}. Thus, if $u_{n}=C_{n+1}$, $n\geq 1$ then Theorem \ref{gral1} categorifies numbers in this sequence $C_{n}$ via these enumeration problems.\par\bigskip

Since the number of compositions (partitions in which the order of the summands is considered) of a positive integer $n$ in which no 1's appear is the Fibonacci number $f_{(n-1)}$ \cite{Andrews}. Then Theorem \ref{gral1} categorifies these numbers by assuming that $u_{n}=f_{(n-1)}$ with $n\geq 4$. If $j>4$, then $\mathrm{dim}_{k}\hspace{0.1cm}(Z(\Lambda_{D^{j}})/C_{j})-1=f_{(j-4)}$, where $C_{j}$ is a $k$-subspace of $Z(\Lambda_{D^{j}})$ isomorphic to $Z(\Lambda_{D^{j-1}})$.\par\bigskip

Theorem \ref{gral1} categorifies the sequence $p(n)$ which gives the number of partitions of a positive integer $n$, recall the Hardy-Ramanujan theorem which states that for large $n$, $p(n)\sim \frac{1}{4n\sqrt{3}}e^{\pi\sqrt{\frac{2n}{3}}}$ (see also the sequence A002865 whose numbers give the differences $p(n)-p(n-1)$). In this case, we use $u_{n}=p(n+1)$, $n\geq1$.\par\bigskip

Another interesting sequence categorified by Theorem \ref{gral1} is the sequence $M(n)$ encoded in the OEIS as A000372, which consists of Dedekind numbers, these numbers count the number of antichains in the powerset $2^{\textbf{n}}$ (i.e., the set consisting of all the subsets of $\textbf{n}=\{1, 2,3,\dots, n\}$) ordered by inclusion or the number of elements in a free distributive lattice on $n$ generators. In this case $u_{n}=M(n)$, for $n\geq1$, worth noting that up to date only 8 numbers of this sequence are known.

\end{Nota}

Agust\'{\i}n Moreno Ca\~{n}adas\hspace{30mm}Pedro Fernando Fern\'andez Espinosa\\
amorenoca@unal.edu.co\hspace{32mm}pffernandeze@unal.edu.co\\
Department of Mathematics\hspace{26mm}Department of Mathematics\\
Universidad Nacional de Colombia\hspace{17mm}Universidad Nacional de Colombia.\\
Kra 30 No 45-03\\
ZIP Code 11001000\\
Bogot\'a-Colombia

\par\bigskip

Isa\'ias David Mar\'in Gaviria\hspace{27mm}Gabriel Bravo R{i}os\\
imaringa@unal.edu.co\hspace{35mm}gbravor@unal.edu.co\\
Department of Mathematics\hspace{26mm}Department of Mathematics\\
Universidad Nacional de Colombia\hspace{17mm}Universidad Nacional de Colombia.

\section{Appendix}\label{Appendix}

In this section, we present the set of helices\par\smallskip
\begin{centering}
 $(4_{P}, 2_{P}, P_{A}=\{p_{3,2},p_{1,3},p_{2,5}\}, P_{B}=\{p_{4,1},p_{1,4}, p_{3,6}\})$\par\smallskip
 \end{centering}
 
 associated to a matrix $P$ of type $\mathscr{H}_{3}$ and defined by the word $W_{P}=BAABAB$. The corresponding copies (see Theorem \ref{invariant11}) in \par\smallskip
 
 \begin{centering}$(4_{P}, 1_{P}, P_{A}=\{p_{2, 1}, p_{1, 3}, p_{3, 6}\}, P_{B}=\{p_{3,2}, p_{2, 4}, p_{4, 5}\})$\par\smallskip
 \end{centering}
 
 are also shown according with the associated word $ABABBA$.
\addtocounter{equation}{1}
\begin{equation}\label{hx1}
 \xymatrix@=30pt{
p_{1,1}&p_{1,2}\ar@[blue][rr]&\underline{p_{1,3}}&\ar@[blue][d]\underline{p_{1,4}}&p_{1,5}&p_{1,6}&\\
p_{2,1}&p_{2,2}&p_{2,3}&p_{2,4}\ar@[blue][r]&\underline{p_{2,5}}&p_{2,6}&j_{B}\\
p_{3,1}\ar@[blue][r]&\underline{p_{3,2}}\ar@[blue][uu]&p_{3,3}&p_{3,4}&p_{3,5}&\underline{p_{3,6}}&\\
\underline{p_{4,1}}\ar@[blue][u]&\ar@[blue][l]p_{4,2}&p_{4,3}&p_{4,4}&p_{4,5}&p_{4,6}&i_{A}\\
}
\end{equation}

\begin{equation}\label{hx2}
 \xymatrix@=30pt{
p_{1,1}&p_{1,2}&\underline{p_{1,3}}&{p_{1,4}}\ar@[blue][l]&p_{1,5}&p_{1,6}&j_{B}\\
\underline{p_{2,1}}&p_{2,2}&p_{2,3}&\underline{p_{2,4}}\ar@[blue][u]&{p_{2,5}}&p_{2,6}\ar@[blue][ll]&\\
p_{3,1}&\underline{p_{3,2}}&p_{3,3}&p_{3,4}&p_{3,5}\ar@[blue][r]&\underline{p_{3,6}}\ar@[blue][u]&\\
{p_{4,1}}&p_{4,2}&p_{4,3}\ar@[blue][rr]&p_{4,4}&\underline{p_{4,5}}\ar@[blue][u]&p_{4,6}&i_{A}\\
}
\end{equation}

\begin{equation}\label{hx3}
 \xymatrix@=30pt{
p_{1,1}&p_{1,2}&\underline{p_{1,3}}&\underline{p_{1,4}}&p_{1,5}&\ar@[red][lll]p_{1,6}&\\
p_{2,1}\ar@[red][rrrr]&p_{2,2}&p_{2,3}&p_{2,4}&\underline{p_{2,5}}\ar@[red][d]&p_{2,6}&j_{B}\\
p_{3,1}&\underline{p_{3,2}}&p_{3,3}&p_{3,4}&p_{3,5}\ar@[red][r]&\underline{p_{3,6}}\ar@[red][uu]&\\
\underline{p_{4,1}}\ar@[red][uu]&\ar@[red][l]p_{4,2}&p_{4,3}&p_{4,4}&p_{4,5}&p_{4,6}&i_{A}\\
}
\end{equation}

\begin{equation}\label{hx4}
 \xymatrix@=30pt{
p_{1,1}&p_{1,2}&\underline{p_{1,3}}\ar@[red][dd]&{p_{1,4}}&p_{1,5}\ar@[red][ll]&p_{1,6}&j_{B}\\
\underline{p_{2,1}}&p_{2,2}\ar@[red][rr]&p_{2,3}&\underline{p_{2,4}}&{p_{2,5}}&p_{2,6}&\\
p_{3,1}&\underline{p_{3,2}}\ar@[red][u]&p_{3,3}\ar@[red][l]&p_{3,4}&p_{3,5}&\underline{p_{3,6}}&\\
{p_{4,1}}&p_{4,2}&p_{4,3}\ar@[red][rr]&p_{4,4}&\underline{p_{4,5}}\ar@[red][uuu]&p_{4,6}&i_{A}\\
}
\end{equation}

\begin{equation}\label{hx5}
 \xymatrix@=30pt{
p_{1,1}&p_{1,2}&\underline{p_{1,3}}&\underline{p_{1,4}}\ar@[green][dd]&\ar@[green][l]p_{1,5}&p_{1,6}&\\
p_{2,1}\ar@[green][rrrr]&p_{2,2}&p_{2,3}&p_{2,4}&\underline{p_{2,5}}\ar@[green][u]&p_{2,6}&j_{B}\\
p_{3,1}&\underline{p_{3,2}}&p_{3,3}&\ar@[green][ll]p_{3,4}&p_{3,5}&\underline{p_{3,6}}&\\
\underline{p_{4,1}}\ar@[green][uu]&\ar@[green][l]p_{4,2}&p_{4,3}&p_{4,4}&p_{4,5}&p_{4,6}&i_{A}\\
}
\end{equation}

\begin{equation}\label{hx6}
 \xymatrix@=30pt{
p_{1,1}&p_{1,2}&\underline{p_{1,3}}\ar@[green][d]&{p_{1,4}}&p_{1,5}\ar@[green][ll]&p_{1,6}&j_{B}\\
\underline{p_{2,1}}&p_{2,2}&p_{2,3}\ar@[green][r]&\underline{p_{2,4}}\ar@[green][d]&{p_{2,5}}&p_{2,6}&\\
p_{3,1}&\underline{p_{3,2}}&p_{3,3}&p_{3,4}\ar@[green][rr]&p_{3,5}&\underline{p_{3,6}}&\\
{p_{4,1}}&p_{4,2}&p_{4,3}\ar@[green][rr]&p_{4,4}&\underline{p_{4,5}}\ar@[green][uuu]&p_{4,6}&i_{A}\\
}
\end{equation}

\begin{equation}\label{hx7}
 \xymatrix@=30pt{
p_{1,1}\ar@[black][rr]&p_{1,2}&\underline{p_{1,3}}\ar@[black][dd]&\underline{p_{1,4}}&p_{1,5}&p_{1,6}&\\
p_{2,1}&p_{2,2}&p_{2,3}&p_{2,4}&\underline{p_{2,5}}&\ar@[black][l]p_{2,6}&j_{B}\\
p_{3,1}&\underline{p_{3,2}}&p_{3,3}\ar@[black][rrr]&p_{3,4}&p_{3,5}&\ar@[black][u]\underline{p_{3,6}}&\\
\underline{p_{4,1}}\ar@[black][uuu]&\ar@[black][l]p_{4,2}&p_{4,3}&p_{4,4}&p_{4,5}&p_{4,6}&i_{A}\\
}
\end{equation}

\begin{equation}\label{hx8}
 \xymatrix@=30pt{
p_{1,1}&p_{1,2}\ar@[black][r]&\underline{p_{1,3}}&{p_{1,4}}&p_{1,5}&p_{1,6}&j_{B}\\
\underline{p_{2,1}}\ar@[black][d]&p_{2,2}&p_{2,3}&\underline{p_{2,4}}&{p_{2,5}}\ar@[black][llll]&p_{2,6}&\\
p_{3,1}\ar@[black][r]&\underline{p_{3,2}}\ar@[black][uu]&p_{3,3}&p_{3,4}&p_{3,5}&\underline{p_{3,6}}&\\
{p_{4,1}}&p_{4,2}&p_{4,3}\ar@[black][rr]&p_{4,4}&\underline{p_{4,5}}\ar@[black][uu]&p_{4,6}&i_{A}\\
}
\end{equation}

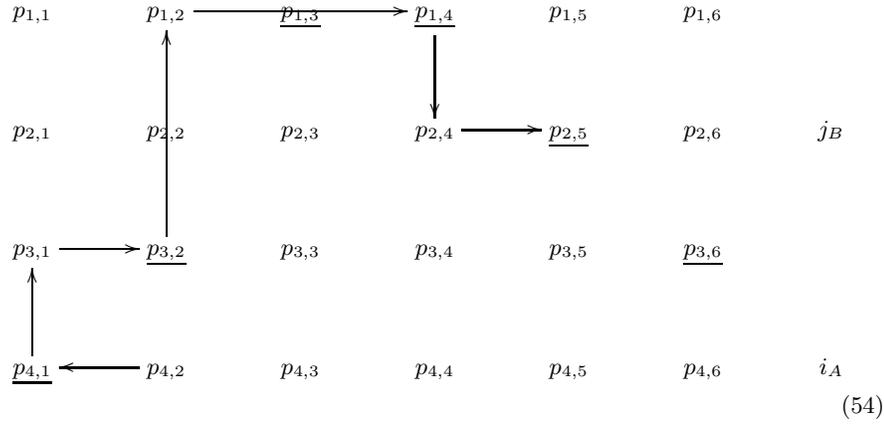
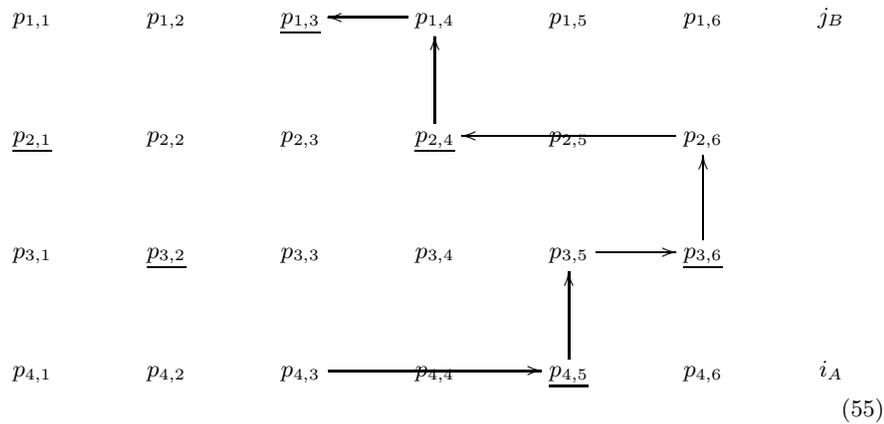
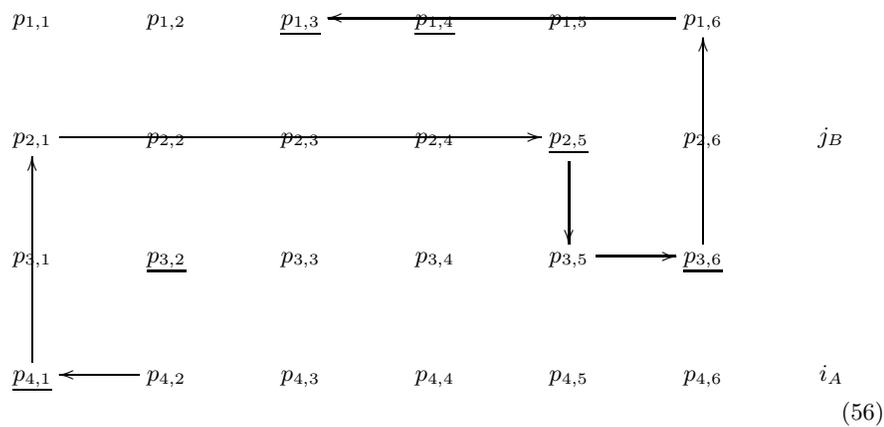
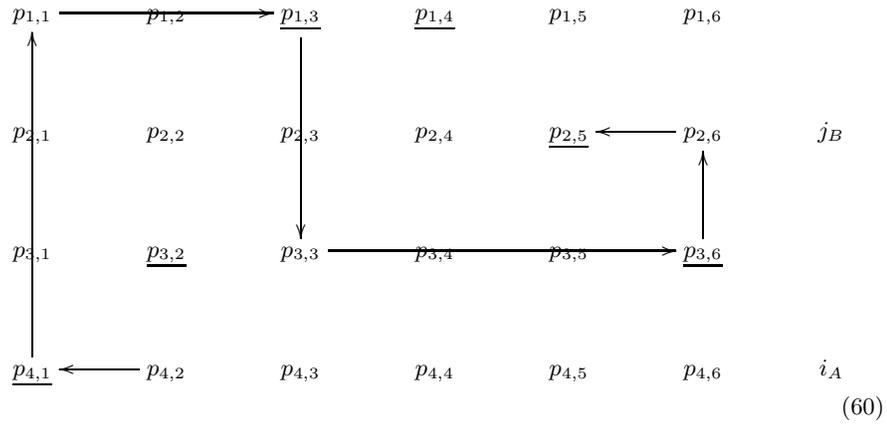
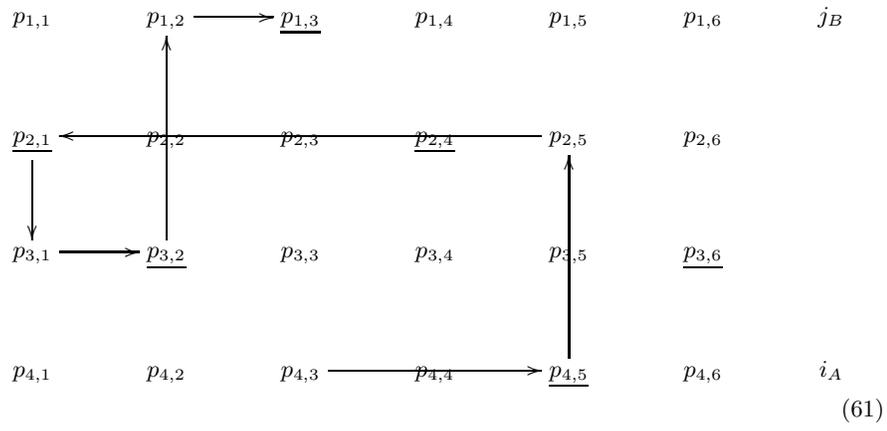
\end{document}